%% file: pap.tex
\documentclass{elsarticle}

\usepackage{bm}
\usepackage{fullpage}
\usepackage{amsfonts}
\usepackage{graphicx}
\usepackage{amsthm}
\usepackage{amssymb}
\usepackage{paralist}
\usepackage{mathtools}
\usepackage[section]{placeins}
\usepackage{algorithm}
\usepackage{algorithmic}
\usepackage{color}
\usepackage{moreverb}
\usepackage[colorlinks, bookmarksopen, bookmarksnumbered, citecolor=red,
            urlcolor=red]{hyperref}
\usepackage{subcaption}
\usepackage[table]{xcolor}

\usepackage{tikz}
\usepackage{pgfplots}
\usetikzlibrary{pgfplots.groupplots}
\usepackage{pgfplotstable, booktabs}

\input{mycommands.tex}

\begin{document}

\title{An Adjoint Method for a High-Order Discretization of Deforming Domain
       Conservation Laws for Optimization of Flow Problems}

\author[rvt1]{M.~J.~Zahr\fnref{fn1}\corref{cor1}}
\ead{mzahr@stanford.edu}

\author[rvt2]{P.-O.~Persson\fnref{fn2}}
\ead{persson@berkeley.edu}

\address[rvt1]{Institute for Computational and Mathematical Engineering,
               Stanford University, Stanford, CA 94035.}
\address[rvt2]{Department of Mathematics, University of California, Berkeley,
               Berkeley, CA 94720-3840.}
\cortext[cor1]{Corresponding author}
\fntext[fn1]{Graduate Student, Institute for Computational and Mathematical
             Engineering, Stanford University}
\fntext[fn2]{Associate Professor, Department of Mathematics, University of
             California, Berkeley.}

\begin{abstract}
The fully discrete adjoint equations and the corresponding adjoint method are
derived for a globally high-order accurate discretization of conservation laws
on parametrized, deforming domains. The conservation law on the deforming domain
is transformed into one on a fixed reference domain by the introduction of a
time-dependent mapping that encapsulates the domain deformation and
parametrization, resulting in an Arbitrary Lagrangian-Eulerian form of the
governing equations. A high-order discontinuous Galerkin method is used to
discretize the transformed equation
in space and a high-order diagonally implicit Runge-Kutta scheme is used for
the temporal discretization. Quantities of interest that take the form of
space-time integrals are discretized in a solver-consistent manner. The
corresponding fully discrete adjoint method is used to compute \emph{exact}
gradients of quantities of interest along the manifold of solutions of the
fully discrete conservation law. These quantities of interest and their
gradients are used in the context of gradient-based PDE-constrained
optimization.

The adjoint method is used to solve two optimal shape and control problems
governed by the isentropic, compressible Navier-Stokes equations. The first
optimization problem seeks the energetically optimal trajectory of a 2D
airfoil given a required initial and final spatial position. The optimization
solver, driven by gradients computed via the adjoint method, reduced the
total energy required to complete the specified mission nearly an order of
magnitude. The second optimization problem seeks the energetically optimal
flapping motion and time-morphed geometry of a 2D airfoil given an equality
constraint on the $x$-directed impulse generated on the airfoil. The
optimization solver satisfied the impulse constraint to greater than $8$
digits of accuracy and reduced the required energy between a factor of $2$
and $10$, depending on the value of the impulse constraint, as compared to the
nominal configuration.
\end{abstract}

\maketitle

\section{Introduction}\label{sec:intro}
Optimization problems constrained by Partial Differential Equations (PDEs)
commonly arise in engineering practice, particularly in the context of design
or control of physics-based systems. A majority of the research in
PDE-constrained optimization has been focused on \emph{steady} or \emph{static}
PDEs, with a large body of literature detailing many aspects of the subject,
including
continuous and discrete adjoint methods
\cite{newman1999overview, nadarajah2000comparison, giles2003algorithm,
      mavriplis2007discrete, mader2008adjoint},
parallel implementations \cite{orozco1992massively, reuther1999aerodynamic},
one-shot or infeasible path methods
\cite{orozco1992massively, ghattas1997optimal},
and generalized reduced gradient or feasible path methods
\cite{newman1999overview, he1997computational}. 
This emphasis on steady problems is largely due to the fact that
\begin{inparaenum}[(a)]
  \item static analysis is sufficient for a large class of problems of
        interest and
  \item unsteady analysis is expensive to perform in a many-query
        setting, such as optimization \cite{yamaleev2008adjoint}.
\end{inparaenum}
However, there is a large class of problems where steady analysis
is insufficient, such as problems that are inherently dynamic and
problems where a steady-state solution does not exist or cannot
be found reliably with numerical methods.
Flapping flight is an example of the first type, a fundamentally unsteady
problem that has become increasingly relevant due to its application to
Micro-Aerial Vehicles (MAVs) \cite{platzer2008flapping}.  Systems with
chaotic solutions, such as those encountered in turbulent flows, are an
example of the second type of problems where steady analysis breaks down.
Design and control of these types of systems calls for
\emph{time-dependent} PDE-constrained optimization of the form
\begin{equation} \label{opt:uns-pde-constr-cont}
  \begin{aligned}
    & \underset{\Ubm,~\mubold}{\text{minimize}}
    & & \int_0^T \int_{\Gammabold}  j(\Ubm, \mubold, t)\,dS \,dt \\
    & \text{subject to}
    & & \int_0^T \int_{\Gammabold} 
                                   \cbm(\Ubm, \mubold, t)\,dS \,dt \leq 0 \\
    & & & \pder{\Ubm}{t} + \nabla \cdot \Fbm(\Ubm, \nabla\Ubm) = 0
  \end{aligned}
\end{equation}
where the last constraint corresponds to a system of conservation laws with
solution $\Ubm$, parametrized by $\mubold$; the objective and constraint
functions of the optimization take the form of space-time integrals of
pointwise, instantaneous quantities of interest $j$ and $\cbm$ over the
surface of the body $\Gammabold$.

In this work, the large computational cost associated with time-dependent
PDE-constrained optimization will be addressed by two means. The first is
the development of a \emph{globally high-order} numerical discretization of
conservation laws on deforming domains%
\footnote{As with all works on high-order methods, high-order
          accuracy relies on sufficient regularity in the solution.}.
For many important problems, high-order
methods have been shown to require fewer spatial degrees of freedom
\cite{wang2013high} and time steps \cite{mani2008unsteady, zahr2013performance}
for a given level of accuracy compared to low-order methods. Highly accurate
quantities of interest, usually the time-average of a relevant surface- or
volume-integrated quantity, is paramount, at least at convergence, since they
drive the optimization trajectory through the objective function and
constraints.  Large errors in quantities of interest will cause the optimization
procedure to be driven by discretization errors causing termination at a
suboptimal design/control. The second approach to reduce the computational
impact of time-dependent optimization is the use of gradient-based optimization
techniques due to their rapid convergence properties, particularly when compared
to derivative-free alternatives.

An efficient technique for computing derivatives of optimization functionals,
required by gradient-based optimization solvers, is the \emph{adjoint method}.
It has proven its utility in the context of
output-based mesh adaptivity and gradient-based PDE-constrained optimization
as only a single linearized dual solve is required to compute the gradient
of a single quantity of interest with respect to any number of parameters.
In the context of partial differential equations, the adjoint equations can
be derived at either the continuous, semi-discrete, or fully discrete level.
The fully discrete adjoint method will be the focus of this work as it ensures
discrete consistency  \cite{mader2008adjoint, yamaleev2008adjoint} of computed
gradients, i.e. the gradient of the discrete solution, including discretization
errors, is computed. Discrete consistency is beneficial in the context of
gradient-based optimization as inconsistent gradients may cause convergence
of black-box optimizers to be slowed or hindered \cite{gill1981practical},
unless specialized optimization algorithms are employed that handle gradient
inexactness \cite{heinkenschloss1995analysis}.


In this work, a \emph{globally high-order} numerical discretization of
general systems of conservation laws, defined on deforming domains, is
introduced and the corresponding fully discrete adjoint equations derived. The
goal is to harness the advantages of high-order methods in the context of
\emph{gradient-based} optimization.  The solution of the adjoint equations --
the dual solution -- will be used to construct exact gradients of fully
discrete quantities of interest. A Discontinuous Galerkin
Arbitrary Lagrangian-Eulerian (DG-ALE) method \cite{persson2009dgdeform} is
used for the high-order spatial discretization.  Previous work on the adjoint
method for conservation laws on deforming domains predominantly considers
a Finite Volume (FV) spatial discretization
\cite{newman1999overview, nadarajah2000comparison, giles2003algorithm,
      mavriplis2007discrete, nadarajah2007optimum, mani2008unsteady,
      mader2008adjoint, nielsen2010discrete},
and recently extended to DG-ALE schemes
\cite{fidkowski2010entropy, fidkowski2011output, van2013adjoint}. The DG-ALE
discretization is chosen rather than FV due to its stable, high-order
discretization of convective fluxes. The Geometric Conservation Law (GCL) is
satisfied in the DG-ALE scheme through the introduction of an element-level
auxiliary equation. The fully discrete adjoint equations derived in this work
fully incorporate this GCL augmentation \cite{kast2013output}, ensuring
discrete consistency is maintained.

High-order temporal discretization will be achieved using a
Diagonally Implicit Runge-Kutta (DIRK) \cite{alexander77dirk} method, marking
a departure from previous work on unsteady adjoints, which has mostly
considered temporal discretization via Backward Differentiation Formulas (BDF)
\cite{he1997computational, nadarajah2007optimum, mani2008unsteady, nielsen2010discrete, van2013adjoint, mishra2015time}, with some work on space-time DG
discretizations \cite{kast2013output}.
Apart from being limited to second-order accuracy, if A-stability is required,
high-order BDF schemes require special techniques for initialization
\cite{zahr2013performance}. While DIRK schemes require additional work to
achieve high-order convergence in the form of additional nonlinear solves at a
given timestep, this investment has been shown to be worthwhile
\cite{zahr2013performance}.  The fully discrete, time-dependent adjoint
equations corresponding to a Runge-Kutta temporal discretization has been
studied in the context of ODEs \cite{sandu2006properties} and optimal
control \cite{joslin1997self, hager2000rk}. A difficulty that arises when
considering Runge-Kutta discretizations is the stages encountered during the
primal (forward) and dual (reverse) solves will not align if the continuous or
semi-discrete adjoint method is employed. This difficulty cannot arise in the
fully discrete formulation as only the terms computed in the primal solve are
used in the derivation of the discrete adjoint equations.

High-order discretization of integrated quantities of interest
will be done in a solver-consistent manner, that is, spatial
integrals will be evaluated via integration of the finite element shape
functions used for the spatial discretization and temporal integrals evaluated
via the DIRK scheme from the temporal discretization. This ensures the
discretization order of quantities of interest exactly matches the PDE
discretization. This marks a departure from the traditional method of simply
using the trapezoidal rule for temporal integration of quantities of interest
\cite{yamaleev2008adjoint, mani2008unsteady, nielsen2010discrete,
van2013adjoint}, which is limited to second-order accuracy. Since the
fully discrete adjoint method is used, discretization of the quantities of
interest is accounted for in the adjoint equations, which is the final
component in ensuring discrete consistency of their gradients.

This notion of a solver-consistent discretization of quantities of interest has
an immediate connection to adjoint-consistency for stationary
problems \cite{arnold2002unified, hartmann2007adjoint}. A necessary condition
for adjoint-consistency is that the quantity of interest is discretized with
the same scheme as the governing equations and boundary conditions
\cite{hartmann2015generalized}; however, this does not refer to the practical
issue of evaluating the integrals that arise in the finite-dimensional
primal or adjoint residual. Solver-consistency ensures the same approximation
is used to evaluate the integrals arising the primal residual and quantity
of interest, and therefore the adjoint residual. A key property of
adjoint-consistent discretizations is they possess optimal convergence rates
in the $L^2$-norm and in quantities of interest
\cite{houston2001hp, harriman2002hp, harriman2004importance}.

With the high-order numerical discretization in place, and the corresponding
adjoint method developed, the proposed methodology is employed to solve
deforming-domain PDE-constrained optimization problems, such as optimal
control and shape optimization for fluid flows.
Existing techniques for parametrizing the domain deformation include
\begin{inparaenum}[(a)]
  \item boundary-driven deformation where the deformation of a boundary is
        parametrized using existing techniques
        \cite{imam1982three, samareh1999survey} and the domain is deformed
        using a structural analogy
        \cite{degand2002three, mani2008unsteady, nielsen2010discrete}, and
  \item analytical expressions for direct domain deformation
        \cite{persson2009dgdeform, fidkowski2011output}.
\end{inparaenum}
The latter approach is adopted in this work. The proposed methodology,
including the high-order primal discretization, fully discrete adjoint method,
and domain parametrization, will be demonstrated on relevant optimal control and
shape optimization problems. An important component will be the incorporation of
\emph{solution-based} constraints, which has previously been done only
through heuristic penalization of the objective function \cite{van2013adjoint}.

The remainder of this document is organized as follows.
Section~\ref{sec:govern} introduces the general form of a system of
conservation laws defined on deforming domains, as well as the globally
high-order discretization using DG-in-space, DIRK-in-time and a
solver-consistent discretization of output integrals.
Section~\ref{sec:adj} derives the fully discrete adjoint method corresponding
to the high-order discretization from Section~\ref{sec:govern}, with
relevant implementational details discussed in Section~\ref{sec:implement}.
Section~\ref{sec:app} couples the adjoint framework to state-of-the-art
optimization software to solve a realistic optimization problem of
determining energetically optimal flapping motions of an airfoil. An alternate
derivation of the adjoint equations is provided in \ref{sec:app-a} by
interpreting the dual variables as Lagrange multipliers of an auxiliary
PDE-constrained optimization problem.

\section{Governing Equations and Discretization}\label{sec:govern}
This section is devoted to the treatment of conservation laws on a
\emph{parametrized, deforming domain} using an Arbitrary Lagrangian-Eulerian
(ALE) description of the governing equations and exposition of a globally
high-order numerical discretization of the ALE form of the system of
conservation laws. Subsequently, Section~\ref{sec:adj} will develop the
corresponding fully discrete adjoint equations and the adjoint method for
constructing gradients of quantities of interest.

The methods introduced in this work
are not necessarily limited to Partial Differential Equations (PDE) that can be
written as conservation laws (\ref{eqn:claw-phys}).  In
Section~\ref{subsec:govern-spatial}, the chosen spatial discretization
(discontinuous Galerkin Arbitrary Lagrangian-Eulerian method) is
applied to the PDE, resulting in a system of first-order Ordinary
Differential Equations (ODE), which is the point of departure for all
adjoint-related derivations.  Time-dependent PDEs that are not conservation
laws can be written similarly at the semi-discrete level after application
of an appropriate spatial discretization, e.g. a continuous finite element
method for parabolic PDEs.  In this work, the scope is limited to first-order
temporal systems, or those which are recast as such.

\subsection{System of Conservation Laws on Deforming Domain:
            Arbitrary Lagrangian-Eulerian Description}
Consider a general system of conservation laws, defined on a parametrized,
deforming domain, $v(\mubold, t)$, written at the continuous level as
\begin{equation}\label{eqn:claw-phys}
  \pder{\Ubm}{t} + \nabla\cdot\Fbm(\Ubm,~\nabla\Ubm) = 0
                                                  \quad \text{in}~~v(\mubold, t)
\end{equation}
where the physical flux is decomposed into an inviscid and a viscous part
$\Fbm(\Ubm,~\nabla\Ubm) = \Fbm^{inv}(\Ubm) + \Fbm^{vis}(\Ubm,~\nabla\Ubm)$,
$\Ubm(\xbm, \mubold, t)$ is the solution of the system of conservation laws,
$t \in (0, T)$ represents time, and $\mubold \in \Rbb^{N_\mubold}$ is a
vector of parameters. This work will focus on the case where the \emph{domain}
is parametrized by $\mubold$, although extension to other types of parameters,
e.g. constants defining the conservation law, is straightforward.

The conservation law on a deforming domain is transformed into a conservation
law on a \emph{fixed} reference domain through the introduction of a
time-dependent mapping between the physical and reference domains, resulting
in an Arbitrary Lagrangian-Eulerian description of the governing equations.

Denote the physical domain by $v(\mubold, t) \subset \Rbb^{n_{sd}}$ and the
fixed, reference domain by $V \subset \Rbb^{n_{sd}}$,
where $n_{sd}$ is the number of spatial dimensions. At each time $t$,
let $\Gcal$ be a time-dependent diffeomorphism between the reference domain
and physical domain: $\xbm(\Xbm, \mubold, t) = \Gcal(\Xbm, \mubold, t)$,
where $\Xbm \in V$ is a point in the reference domain and
$\xbm(\Xbm, \mubold, t) \in v(\mubold, t)$ is the corresponding point
in the physical domain at time $t$ and parameter configuration $\mubold$.

\begin{figure}
  \centering
  \includegraphics[width=2.5in]{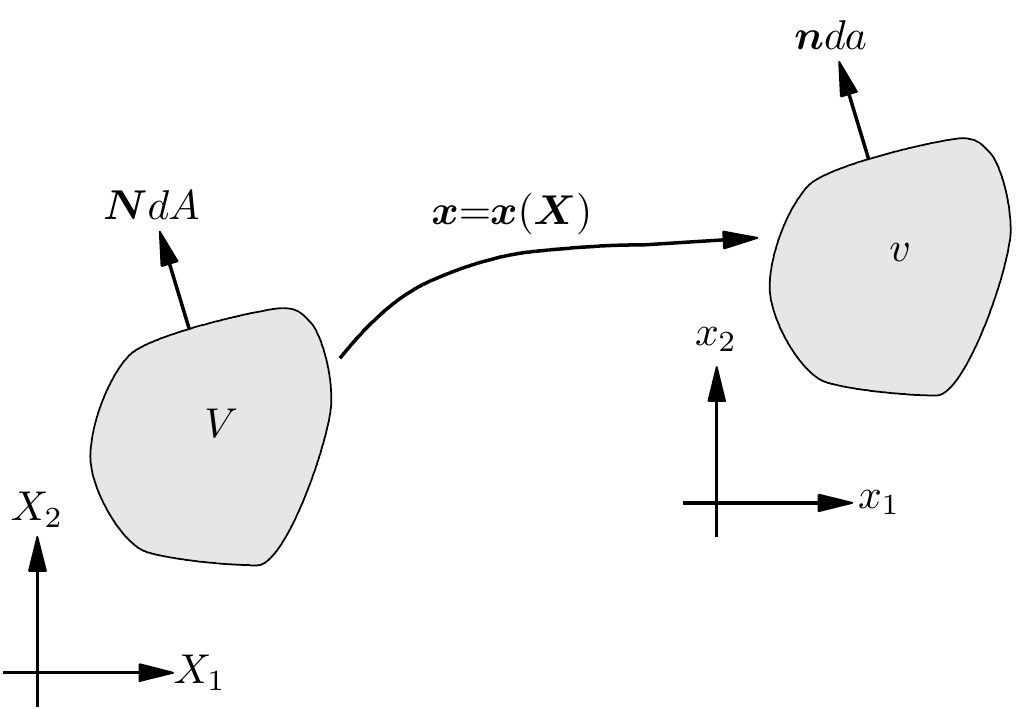}
  \caption{Time-dependent mapping between reference and physical domains.}
\end{figure}

The transformed system of conservation laws from (\ref{eqn:claw-phys}),
under the mapping $\Gcal$, defined on the reference domain takes the form
\begin{equation} \label{eqn:ns-cons-ref}
  \left.\pder{\Ubm_\Xbm}{t}\right|_\Xbm +
  \nabla_\Xbm\cdot\Fbm_\Xbm(\Ubm_\Xbm,~\nabla_\Xbm\Ubm_\Xbm) = 0
\end{equation}
where $\nabla_\Xbm$ denotes spatial derivatives with respect to the reference
variables, $\Xbm$.  The transformed state vector, $\Ubm_\Xbm$, and its
corresponding spatial gradient with respect to the reference configuration
take the form
\begin{equation} \label{eqn:transform1-nogcl}
\Ubm_\Xbm = g\Ubm, \qquad
\nabla_\Xbm \Ubm_\Xbm = g^{-1}\Ubm_\Xbm\pder{g}{\Xbm} +
                        g\nabla\Ubm\cdot\Gbm,
\end{equation}
where $\Gbm = \nabla_\Xbm \Gcal$, $g = \text{det}(\Gbm)$,
$\displaystyle{\vbm_\Gbm = \pder{\xbm}{t} = \pder{\Gcal}{t}}$,
and the arguments have been dropped, for brevity. The transformed fluxes are
\begin{equation} \label{eqn:transform2-nogcl}
  \begin{aligned}
    \Fbm_\Xbm(\Ubm_\Xbm, \nabla_\Xbm\Ubm_\Xbm) &= \Fbm_\Xbm^{inv}(\Ubm_\Xbm) +
                     \Fbm_\Xbm^{vis}(\Ubm_\Xbm, \nabla_\Xbm\Ubm_\Xbm), \qquad \\
    \Fbm_\Xbm^{inv}(\Ubm_\Xbm) &= g\Fbm^{inv}(g^{-1}\Ubm_\Xbm)\Gbm^{-T} -
                                 \Ubm_\Xbm \otimes \Gbm^{-1}\vbm_\Gbm, \qquad \\
    \Fbm_\Xbm^{vis}(\Ubm_\Xbm, \nabla_\Xbm\Ubm_\Xbm) &=
     g\Fbm^{vis}\left(g^{-1}\Ubm_\Xbm, g^{-1}\left[\nabla_\Xbm\Ubm_\Xbm
                - g^{-1}\Ubm_\Xbm\pder{g}{\Xbm}\right]\Gbm^{-1}\right)\Gbm^{-T}.
  \end{aligned}
\end{equation}
For details regarding the derivation of the transformed equations, the reader
is referred to \cite{persson2009dgdeform}.

When integrated using inexact numerical schemes,
this ALE formulation does not satisfy the Geometric
Conservation Law (GCL) \cite{farhat2001discrete,persson2009dgdeform}.  This
is overcome by introduction of an auxiliary variable $\bar{g}$, defined as
the solution of
\begin{equation} \label{eqn:gbar}
  \pder{\bar{g}}{t} - \nabla_\Xbm\cdot\left(g\Gbm^{-1}\vbm_\Gbm\right) = 0.
\end{equation}
The auxiliary variable, $\bar{g}$ is used to modify the \emph{transformed}
conservation law according to
\begin{equation} \label{eqn:ns-cons-ref-gcl}
  \left.\pder{\Ubm_{\bar\Xbm}}{t}\right|_\Xbm +
  \nabla_\Xbm\cdot\Fbm_{\bar\Xbm}(\Ubm_{\bar\Xbm},~\nabla_\Xbm\Ubm_{\bar\Xbm})
  = 0
\end{equation}
where the GCL-transformed state variables are
\begin{equation} \label{eqn:transform1-gcl}
\Ubm_{\bar\Xbm} = \bar{g}\Ubm, \qquad
\nabla_\Xbm \Ubm_{\bar\Xbm} = \bar{g}^{-1}\Ubm_{\bar\Xbm}\pder{\bar{g}}{\Xbm} +
                              \bar{g}\nabla\Ubm \cdot \Gbm
\end{equation}
and the corresponding fluxes
\begin{equation} \label{eqn:transform2-gcl}
  \begin{aligned}
    \Fbm_{\bar\Xbm}(\Ubm_{\bar\Xbm}, \nabla_\Xbm\Ubm_{\bar\Xbm}) &=
               \Fbm_{\bar\Xbm}^{inv}(\Ubm_{\bar\Xbm}) +
               \Fbm_{\bar\Xbm}^{vis}(\Ubm_{\bar\Xbm},
                     \nabla_\Xbm\Ubm_{\bar\Xbm}), \qquad \\
    \Fbm_{\bar\Xbm}^{inv}(\Ubm_{\bar\Xbm}) &=
                             g\Fbm^{inv}(\bar{g}^{-1}\Ubm_{\bar\Xbm})\Gbm^{-T} -
                           \Ubm_{\bar\Xbm} \otimes \Gbm^{-1}\vbm_\Gbm, \qquad \\
    \Fbm_{\bar\Xbm}^{vis}(\Ubm_{\bar\Xbm}, \nabla_\Xbm\Ubm_{\bar\Xbm}) &=
     g\Fbm^{vis}\left(\bar{g}^{-1}\Ubm_{\bar\Xbm},
                      \bar{g}^{-1}\left[\nabla_\Xbm\Ubm_{\bar\Xbm}
                 - \bar{g}^{-1}\Ubm_{\bar\Xbm}\pder{\bar{g}}{\Xbm}\right]
                                                      \Gbm^{-1}\right)\Gbm^{-T}.
  \end{aligned}
\end{equation}
It was shown in \cite{persson2009dgdeform} that the transformed equations
(\ref{eqn:ns-cons-ref-gcl}) satisfy the GCL. In the next section, the ALE
description of the governing equations (\ref{eqn:ns-cons-ref}) and
(\ref{eqn:ns-cons-ref-gcl}) will be converted to first-order form and
discretized via a high-order discontinuous Galerkin method.

\subsection{Spatial Discretization: Arbitrary Lagrangian-Eulerian Discontinuous
            Galerkin Method}\label{subsec:govern-spatial}

The ALE description of the conservation law without GCL augmentation will
be considered first. To proceed, the second-order system of partial differential
equations in (\ref{eqn:ns-cons-ref}) is converted to first-order form
\begin{equation} \label{eqn:ns-cons-ref-first}
  \begin{aligned}
    \left.\pder{\Ubm_\Xbm}{t}\right|_\Xbm +
        \nabla_\Xbm\cdot\Fbm_\Xbm(\Ubm_\Xbm,~\Qbm_\Xbm) &= 0 \\
    \Qbm_\Xbm - \nabla_\Xbm\Ubm_\Xbm &= 0,
  \end{aligned}
\end{equation}
where $\Qbm_\Xbm$ is introduced as an auxiliary variable to represent the
spatial gradient of the $\Ubm_\Xbm$. Equation (\ref{eqn:ns-cons-ref-first})
is discretized using a standard nodal discontinuous Galerkin finite element
method \cite{cockburn01rkdg}, which after local elimination of the auxiliary
variables $Q_\Xbm$ leads to the following system of ODEs
\begin{equation}\label{eqn:semidisc-nogcl}
  \Mbb_\Xbm \pder{\ubm_\Xbm}{t} = \rbm_{\ubm_\Xbm}(\ubm_\Xbm, \mubold, t),
\end{equation}
where $\Mbb_\Xbm$ is the block-diagonal, symmetric, \emph{fixed} mass matrix,
$\ubm_\Xbm$ is the vectorization of $\Ubm_\Xbm$ at all nodes in the high-order
mesh, and $\rbm_{\ubm_\Xbm}$ is the nonlinear function defining the DG
discretization of the inviscid and viscous fluxes.

The GCL augmentation is treated identically, i.e. conversion to first-order
form and subsequent application of the discontinuous Galerkin finite element
method, where $\Ubm_{\bar\Xbm}$ is taken as the state variable. The result is
a system of ODEs corresponding to a high-order ALE scheme that satisfies the GCL
\begin{equation}\label{eqn:semidisc-gcl1}
  \begin{aligned}
    \Mbb_{\bar\gbold} \pder{\bar\gbold}{t} &= \rbm_{\bar\gbold}(\mubold, t)\\
    \Mbb_\Xbm \pder{\ubm_{\bar\Xbm}}{t} &= \rbm_{\ubm_{\bar\Xbm}}(
                                                 \ubm_{\bar\Xbm}, \bar\gbold,
                                                 \mubold, t)
  \end{aligned}
\end{equation}
where each term is defined according to their counterparts in
(\ref{eqn:semidisc-nogcl}). From the conservation law defining $\bar{g}$
(\ref{eqn:gbar}), the corresponding flux is continuous, implying the physical
flux $g\Gbm^{-1}\vbm_\Gbm$ can be used as the numerical flux.  This implies no
information is required from neighboring elements and (\ref{eqn:gbar}) can be
solved at the element level, i.e. statically condensed. Furthermore, the
$\bar{\gbold}$ residual, $\rbm_{\bar{\gbold}}$, does not depend on
$\bar{\gbold}$ itself since the physical flux $g\Gbm^{-1}\vbm_\Gbm$ is
independent of $\bar{g}$.

Since the equation for $\bar\gbold$ does not depend on $\ubm_{\bar\Xbm}$,
it can be solved independently of the equation for $\ubm_{\bar\Xbm}$. This
enables $\bar\gbold$ to be considered an implicit function of $\mubold$, i.e.
$\bar\gbold = \bar\gbold(\mubold, t)$ through application of the implicit
function theorem. Then, (\ref{eqn:semidisc-gcl1}) reduces to
\begin{equation} \label{eqn:semidisc-gcl}
  \Mbb_\Xbm \pder{\ubm_{\bar\Xbm}}{t} = \rbm_{\ubm_{\bar\Xbm}}(
                                              \ubm_{\bar\Xbm},
                                              \bar\gbold(\mubold, t),
                                              \mubold, t).
\end{equation}
Equations (\ref{eqn:semidisc-nogcl}) and (\ref{eqn:semidisc-gcl}) are
abstracted into the following system of ODEs
\begin{equation} \label{eqn:semidisc}
  \Mbb\pder{\ubm}{t} = \rbm(\ubm, \mubold, t),
\end{equation}
for convenience in the derivation of the fully discrete adjoint equations.
Evaluation of the residual, $\rbm$, in (\ref{eqn:semidisc}) at 
parameter $\mubold$ and time $t$ requires evaluation of the mapping,
$\xbm(\mubold, t)$ and $\dot\xbm(\mubold, t)$, and $\bar\gbold(\mubold, t)$,
if GCL augmentation is employed. The implicit dependence of $\bar\gbold$ on
$\mubold$ requires special treatment when computing derivatives with respect
to $\mubold$, which will be required in the adjoint method
(Section~\ref{sec:adj}). Discussion of treatment of such terms will be
deferred to Section~\ref{subsec:implement-partial-gcl}.

A convenient property of this DG-ALE scheme is that all computations are
performed on the reference domain which is \emph{independent} of time and
parameter. An implication of this is that the mass matrix of the ODE
(\ref{eqn:semidisc}) is also time- and parameter-independent.  This property
simplifies all adjoint computations introduced in Section~\ref{sec:adj} as
terms involving
$\displaystyle{\pder{\Mbb}{\ubm}}$ and $\displaystyle{\pder{\Mbb}{\mubold}}$
are identically zero.
This, in turn, simplifies the implementation of the adjoint method
and translates to computational savings since contractions with these
third-order tensor are not required; see \cite{he1997computational} for a
discretization with parameter-dependent mass matrices and the corresponding
adjoint derivation. In subsequent sections, it will be assumed that the
mass matrix is time- and parameter-independent.

The DG-ALE scheme outlined in this section constitutes a \emph{spatial}
discretization, which yields a system of Ordinary Differential Equations (ODEs)
when applied to the PDE in (\ref{eqn:claw-phys}).  The semi-discrete form of the
conservation law is the point of departure for the remainder of this document.
The subsequent development applies to any system of ODEs of the form
(\ref{eqn:semidisc}) without relying on the specific spatial discretization
scheme employed. The DG-ALE scheme was chosen to provide a high-order, stable
spatial discretization of the conservation law (\ref{eqn:claw-phys}).
With the high-order spatial discretization of the state equations introduced,
focus is shifted toward high-order temporal discretization to yield a
globally high-order accurate numerical scheme.



\subsection{Temporal Discretization: Diagonally Implicit Runge-Kutta}
\label{subsec:govern-temporal}
The two prevailing classes of high-order implicit temporal integration
schemes are:
\begin{inparaenum}[(a)]
  \item Backward Differentiation Formulas (BDF) and
  \item Implicit Runge-Kutta (IRK).
\end{inparaenum}
BDF schemes are popular since high-order accuracy can be achieved at the cost
of the solution of a single nonlinear system of equations of size $N_\ubm$ at
each time step.  However, they suffer from initialization issues and are
limited to second-order, if A-stability is required.  In contrast, IRK schemes
are single-step methods that can be A-stable and arbitrarily high-order,
at the cost of solving an \emph{enlarged} nonlinear system of equations of
size $s \cdot N_\ubm$, for an $s$-stage scheme.  For practical problems,
this can be prohibitively expensive, in terms of memory and CPU time.

\begin{table}
  \centering
  \begin{tabular} { l | c c c c}
     $c_1$ & $a_{11}$  &  &  & \\
     $c_2$ & $a_{21}$ & $a_{22}$ & & \\
     $\vdots$ &  $\vdots$ & $\vdots$ & $\ddots$ & \\
     $c_s$ & $a_{s1}$ & $a_{s2}$ & $\cdots$ & $a_{ss}$ \\ \hline
     & $b_1$ & $b_2$ & $\cdots$ & $b_s$
  \end{tabular} 
  \caption{Butcher Tableau for $s$-stage diagonally implicit Runge-Kutta scheme}
  \label{tab:dirk}
\end{table}

A particular subclass of the IRK schemes, known as Diagonally Implicit
Runge-Kutta (DIRK) schemes \cite{alexander77dirk}, are capable of achieving
high-order accuracy with the desired stability properties, without requiring
the solution of an enlarged system of equations.  The DIRK schemes are
defined by a \emph{lower triangular} Butcher tableau (Table~\ref{tab:dirk})
and take the following form when applied to (\ref{eqn:semidisc})
\begin{equation} \label{eqn:dirk}
  \begin{aligned}
    \ubm^{(0)} &= \ubm_0(\mubold) \\
    \ubm^{(n)} &= \ubm^{(n-1)} + \sum_{i = 1}^s b_i\kbm^{(n)}_i \\
    \Mbb\kbm^{(n)}_i &= \Delta t_n\rbm\left(\ubm_i^{(n)},~\mubold,~
                                            t_{n-1} + c_i\Delta t_n\right),
  \end{aligned}
\end{equation}
for $n = 1, \dots, N_t$ and $i = 1, \dots, s$, where $N_t$ are the number of
time steps in the temporal discretization and $s$ is the number of stages
in the DIRK scheme.  The temporal domain, $[0,~T]$ is discretized into
$N_t$ segments with endpoints $\{t_0, t_1, \dots, t_{N_t}\}$, with the $n$th
segment having length $\Delta t_n = t_n - t_{n-1}$ for $n = 1, \dots, N_t$.
Additionally, in (\ref{eqn:dirk}), $\ubm_i^{(n)}$ is used to denote the
approximation of $\ubm^{(n)}$ at the $i$th stage of time step $n$
\begin{equation}\label{eqn:dirk-stage}
  \ubm_i^{(n)} = \ubm_i^{(n)}(\ubm^{(n-1)},~\kbm_1^{(n)},\dots,~\kbm_s^{(n)})
               = \ubm^{(n-1)} + \sum_{j = 1}^i a_{ij}\kbm^{(n)}_j.
\end{equation}
From (\ref{eqn:dirk}), a complete time step requires the solution of a
sequence of $s$ nonlinear systems of equation of size $N_\ubm$. 

Given this exposition on the spatio-temporal discretization of the
deforming domain conservation law, the expression for the quantity of
interest must be discretized. The next section introduces a spatio-temporal
discretization of the quantities of interest that is \emph{consistent} with
that of the conservation law itself.

\subsection{Solver-Consistent Discretization of Quantities of Interest}
In this section, the high-order discretization of quantities of interest --
space-time integrals of a nonlinear function of the solution of the
conservation law -- is considered. The output quantity takes the form
\begin{equation} \label{eqn:funcl}
 \Fcal(\Ubm, \mubold, t) = \int_{0}^t \int_\Gammabold w(\xbm, \tau)
                                        f(\Ubm, \mubold, \tau)~dS~d\tau.
\end{equation}
In the context of the optimization problem in (\ref{opt:uns-pde-constr-cont}),
$\Fcal$ corresponds to either the objective or a constraint function.

Discretization of quantities of interest will introduce an \emph{additional}
discretization error, on top of that in the approximation of $\Ubm$ itself. To
ensure the quantity of interest discretization does not dominate, thereby
lowering the global order of the scheme, it is necessary that its
discretization order \emph{matches} that of the discretization of the
governing equation.  Clearly, it is wasteful to discretize this to a higher
order than the state equation, using a similar argument.

For these reasons, discretization of (\ref{eqn:funcl}) will be done in a
\emph{solver-consistent} manner, i.e. the spatial and temporal discretization
used for the governing equation will also be used for the quantities of
interest. Define $f_h$ as the approximation of
$\displaystyle{\int_\Gammabold w(\xbm, t)f(\Ubm, \mubold, t)~dS}$
using the DG shape functions from the spatial discretization of the governing
equations.
Then, the solver-consistent spatial discretization of (\ref{eqn:funcl}) becomes
\begin{equation} \label{eqn:funcl-semidisc}
 \Fcal_h(\ubm, \mubold, t) = \int_{0}^t f_h(\ubm, \mubold, \tau)~d\tau,
\end{equation}
ensuring the spatial integration error in the quantity of interest exactly
matches that of the governing equations. Solver-consistent temporal
discretization requires the semi-discrete functional in
(\ref{eqn:funcl-semidisc}) be converted to an ODE, which is accomplished via
differentiation of (\ref{eqn:funcl-semidisc}) with respect to $t$
\begin{equation} \label{eqn:funcl-semidisc-ode}
 \dot\Fcal_h(\ubm, \mubold, t) = f_h(\ubm, \mubold, t).
\end{equation}
Augmenting the semi-discrete governing equations with this ODE
(\ref{eqn:funcl-semidisc-ode}) yields the system of ODEs
\begin{equation}\label{eqn:ode-aug}
  \begin{bmatrix} \Mbb & \zerobold \\ \zerobold & 1 \end{bmatrix}
  \begin{bmatrix} \dot\ubm \\ \dot\Fcal_h \end{bmatrix} =
  \begin{bmatrix} \rbm(\ubm, \mubold, t) \\ f_h(\ubm, \mubold, t) \end{bmatrix}.
\end{equation}
Application of the DIRK temporal discretization introduced in
Section~\ref{subsec:govern-temporal} yields the fully discrete governing
equations and corresponding solver-consistent discretization of the quantity of
interest (\ref{eqn:funcl})
\begin{equation} \label{eqn:dirk1}
  \begin{aligned}
    \ubm^{(n)} &= \ubm^{(n-1)} + \sum_{i=1}^s b_i \kbm_i^{(n)}\\
    \Fcal_h^{(n)} &= \Fcal_h^{(n-1)} + \sum_{i=1}^s b_i
           f_h\left(\ubm_i^{(n)},~\mubold,~t_{n-1}+c_i\Delta t_n\right) \\
    \Mbb\kbm_i^{(n)} &= \Delta t_n\rbm\left(\ubm_i^{(n)},~\mubold,~
                                            t_{n-1} + c_i\Delta t_n\right),
  \end{aligned}
\end{equation}
for $n = 1, \dots, N_t$, $i = 1, \dots, s$, and $\ubm_i^{(n)}$ is defined in
(\ref{eqn:dirk-stage}).  Finally, the functional in (\ref{eqn:funcl}) is
evaluated at time $t = T$ to yield the solver-consistent approximation
of $\Fcal(\ubm, \mubold, T)$
\begin{equation} \label{eqn:funcl-disc}
  F(\ubm^{(0)}, \dots, \ubm^{(N_t)}, \kbm_1^{(1)}, \dots, \kbm_s^{(N_t)}) =
 \Fcal_h^{(N_t)} \approx \Fcal(\ubm, \mubold, T).
\end{equation}

\begin{Rem}
In addition to the standard choice of $w(\xbm, t) = 1$, other common choices
include
\begin{inparaenum}[(a)]
 \item instantaneous, spatially integrated quantities, i.e.
       $w(\xbm, t) = \delta(t - t^*)$,
 \item pointwise, temporally integrated quantities, i.e.
       $w(\xbm, t) = \delta(\norm{\xbm - \xbm^*})$, and
 \item pointwise, instantaneous quantities, i.e.
       $w(\xbm, t) = \delta(\norm{\xbm - \xbm^*})\delta(t - t^*)$.
\end{inparaenum}
where $\xbm^* \in \Gammabold$ and $t^* \in (0, T)$.
For these specific choices of $w(\xbm, t)$, simplifications to the above
discretization can be leveraged since the spatial and/or temporal integrals
reduce to an evaluation of $f$ at $\xbm = \xbm^*$ and/or $t = t^*$.
\end{Rem}

While the spatially solver-consistent discretization of quantities of interest
is widely used, particularly in the context of finite element methods,
temporal discretization is commonly done via low-order quadrature rules,
usually the trapezoidal rule \cite{van2013adjoint, mavriplis2007discrete,
mani2008unsteady, yamaleev2008adjoint, jones2013adjoint}.
The main advantage of this solver-consistent discretization
is the asymptotic discretization order of the governing equation
and quantity of interest are guaranteed to exactly match,
which ensures there is no wasted error in ``over-integrating'' one of
the terms.
Furthermore, solver-consistency extends the notion of adjoint-consistency
that requires the discrete adjoint equations represent a consistent
discretization of the continuous adjoint equations, but does not necessarily
address the practical issue of evaluating the integrals that arise in each
term of the finite-dimensional primal and adjoint residual. Solver-consistency
ensures the integrals arising in the primal residual and quantity of
interest are approximated identically, which in turn implies that each term in
the adjoint equations uses the same integral approximation.
The solver-consistent discretization also has the advantage of
a natural and convenient implementation given the spatial and temporal
discretization implementation.  Finally, this method has the additional
convenience of keeping a high-order accurate ``current'' value of the integral, 
i.e. at time step $n$,
$\displaystyle{\Fcal_h^{(n)} \approx \int_{0}^{t_n} f_h(\tau)\,d\tau}$ to
high-order accuracy. This property does not hold for high-order numerical
quadrature since
$\displaystyle{\int_{0}^{t_n } f_h(\tau)\,d\tau}$ will involve $\ubm^{(n+j)}$,
where $j \geq 1$ depends on the quadrature rule used.

\section{Fully Discrete, Time-Dependent Adjoint Equations}\label{sec:adj}
The purpose of this section is to derive an expression for the total derivative
of the discrete quantity of interest $F$ in (\ref{eqn:funcl-disc}), which can
be expanded as
\begin{equation} \label{eqn:funcl-disc-chain}
  \oder{F}{\mubold} = \pder{F}{\mubold} +
                      \sum_{n=0}^{N_t} \pder{F}{\ubm^{(n)}}
                                       \pder{\ubm^{(n)}}{\mubold} + 
                      \sum_{n=1}^{N_t} \sum_{i=1}^s\pder{F}{\kbm_i^{(n)}}
                                       \pder{\kbm_i^{(n)}}{\mubold},
\end{equation}
\emph{that does not depend on the sensitivities of the state
      variables},
$\displaystyle{\pder{\ubm^{(n)}}{\mubold}}$ and
$\displaystyle{\pder{\kbm_i^{(n)}}{\mubold}}$.
Each of the $N_\mubold$ state
variable sensitivities is the solution of a linear evolution equation of the
same dimension and number of steps as the primal equation (\ref{eqn:dirk}),
rendering these quantities intractable to compute when $N_\mubold$ is large.
Elimination of the state variable sensitivities from
(\ref{eqn:funcl-disc-chain}) is accomplished through introduction of the
adjoint equations corresponding to the functional $F$, and the corresponding
dual variables.  From the derivation of the adjoint equation in
Section~\ref{subsec:adj-deriv}, an expression for the reconstruction of the
gradient of $F$, independent of the state variables sensitivities, follows
naturally. At this point, it is emphasized that $F$ represents \emph{any}
quantity of interest whose gradient is desired, such as the optimization
objective function or a constraint. This section concludes with a discussion
of the advantages of the fully discrete framework in the setting of the
high-order numerical scheme.

Before proceeding to the derivation of the adjoint method, the following
definitions are introduced for the Runge-Kutta stage equations and
state updates
\begin{equation} \label{eqn:abs-govern-disc}
  \begin{aligned}
    \tilde\rbm^{(0)}(\ubm^{(0)},~\mubold)
    =& ~\ubm^{(0)} - \ubm_0(\mubold) = 0\\
    \tilde\rbm^{(n)}(\ubm^{(n-1)},~\ubm^{(n)},~
                        \kbm_1^{(n)}, \dots, \kbm_s^{(n)},~
                        \mubold)
    =& ~\ubm^{(n)}-\ubm^{(n-1)}-\sum_{i=1}^sb_i\kbm_i^{(i)} = 0\\
    \Rbm_i^{(n)}(\ubm^{(n-1)},
                  \kbm_1^{(n)}, \dots, \kbm_i^{(n)},
                  \mubold)
    =& ~\Mbb\kbm_i^{(n)} - \Delta t_n\rbm\left(\ubm_i^{(n)},~\mubold,~
       t_{n-1}+c_i\Delta t_n\right) = 0
  \end{aligned}
\end{equation}
for $n = 1, \dots, N_t$ and $i = 1, \dots, s$.  Differentiation of these
expressions with respect to $\mubold$ gives rise to the \emph{fully 
discrete sensitivity equations}
\begin{equation} \label{eqn:uns-disc-sens}
  \begin{aligned}
    \pder{\tilde\rbm^{(0)}}{\mubold} +
    \pder{\tilde\rbm^{(0)}}{\ubm^{(0)}}\pder{\ubm^{(0)}}{\mubold} &= 0 \\
    \pder{\tilde\rbm^{(n)}}{\mubold} +
    \pder{\tilde\rbm^{(n)}}{\ubm^{(n)}}\pder{\ubm^{(n)}}{\mubold} +
    \pder{\tilde\rbm^{(n)}}{\ubm^{(n-1)}}\pder{\ubm^{(n-1)}}{\mubold} +
    \sum_{p=1}^s \pder{\tilde\rbm^{(n)}}{\kbm_p^{(n)}}
                                           \pder{\kbm_p^{(n)}}{\mubold} &= 0 \\
  \pder{\Rbm_i^{(n)}}{\mubold} +
  \pder{\Rbm_i^{(n)}}{\ubm^{(n-1)}}\pder{\ubm^{(n-1)}}{\mubold} +
  \sum_{j=1}^i\pder{\Rbm_i^{(n)}}{\kbm_j^{(n)}}\pder{\kbm_j^{(n)}}{\mubold} &= 0
  \end{aligned}
\end{equation}
where $n = 1, \dots, N_t$, $i = 1, \dots, s$, and arguments have been dropped.

\subsection{Derivation} \label{subsec:adj-deriv}
The derivation of the fully discrete adjoint equations corresponding to the
quantity of interest, $F$, begins with the introduction of test variables
\begin{equation}
  \lambdabold^{(0)},~\lambdabold^{(n)},~\kappabold_i^{(n)} \in \Rbb^{N_\ubm}
\end{equation}
for $n = 1, \dots, N_t$ and $i = 1, \dots, s$. To eliminate the state
sensitivities from the expression for $\displaystyle{\oder{F}{\mubold}}$ in
(\ref{eqn:funcl-disc-chain}), multiply the sensitivity equations
(\ref{eqn:uns-disc-sens}) by the test variables, integrate (i.e. sum in
the discrete setting) over the time domain, and subtract from the expression for
the gradient in (\ref{eqn:funcl-disc-chain}) to obtain
\begin{equation} \label{eqn:uns-disc-adj-derive-1}
  \begin{aligned}
    \oder{F}{\mubold} = \pder{F}{\mubold} &+
                        \sum_{n=0}^{N_t} \pder{F}{\ubm^{(n)}}
                                         \pder{\ubm^{(n)}}{\mubold} + 
                        \sum_{n=1}^{N_t} \sum_{i=1}^s\pder{F}{\kbm_i^{(n)}}
                                         \pder{\kbm_i^{(n)}}{\mubold} -
                 {\lambdabold^{(0)}}^T\left[\pder{\tilde\rbm^{(0)}}{\mubold}
          + \pder{\tilde\rbm^{(0)}}{\ubm^{(0)}}\pder{\ubm^{(0)}}{\mubold}\right]
          \\ &- \sum_{n=1}^{N_t} {\lambdabold^{(n)}}^T\left[
            \pder{\tilde\rbm^{(n)}}{\mubold} +
            \pder{\tilde\rbm^{(n)}}{\ubm^{(n)}}\pder{\ubm^{(n)}}{\mubold} +
            \pder{\tilde\rbm^{(n)}}{\ubm^{(n-1)}}\pder{\ubm^{(n-1)}}{\mubold} +
            \sum_{p=1}^s \pder{\tilde\rbm^{(n)}}{\kbm_p^{(n)}}
                                             \pder{\kbm_p^{(n)}}{\mubold}\right]
          \\ &- \sum_{n=1}^{N_t} \sum_{i=1}^s {\kappabold_i^{(n)}}^T
        \left[\pder{\Rbm_i^{(n)}}{\mubold} +
             \pder{\Rbm_i^{(n)}}{\ubm^{(n-1)}}\pder{\ubm^{(n-1)}}{\mubold} +
             \sum_{j=1}^i\pder{\Rbm_i^{(n)}}{\kbm_j^{(n)}}
                         \pder{\kbm_j^{(n)}}{\mubold}\right].
   \end{aligned}
\end{equation}
The right side of the equality in (\ref{eqn:uns-disc-adj-derive-1}) is an
equivalent expression for $\displaystyle{\oder{F}{\mubold}}$
\emph{for any value of the test variables} since
the terms in the brackets are zero, i.e. the sensitivity equations.
Re-arrangement of terms in (\ref{eqn:uns-disc-adj-derive-1}) leads to the
following expression for $\displaystyle{\oder{F}{\mubold}}$, where the
state variable sensitivities have been isolated
\begin{equation} \label{eqn:uns-disc-adj-derive-2}
  \begin{aligned}
    \oder{F}{\mubold} = \pder{F}{\mubold} &+
   \left[\pder{F}{\ubm^{(N_t)}} - {\lambdabold^{(N_t)}}^T
    \pder{\tilde\rbm^{(N_t)}}{\ubm^{(N_t)}}\right]\pder{\ubm^{(N_t)}}{\mubold} -
        \sum_{n=0}^{N_t}{\lambdabold^{(n)}}^T\pder{\tilde\rbm^{(n)}}{\mubold} - 
        \sum_{n=1}^{N_t}\sum_{p=1}^s {\kappabold_p^{(n)}}^T
                                                  \pder{\Rbm_p^{(n)}}{\mubold}\\
         &+ \sum_{n=1}^{N_t}\left[\pder{F}{\ubm^{(n-1)}} -
                 {\lambdabold^{(n-1)}}^T\pder{\tilde\rbm^{(n-1)}}{\ubm^{(n-1)}}-
                 {\lambdabold^{(n)}}^T\pder{\tilde\rbm^{(n)}}{\ubm^{(n-1)}}-
                        \sum_{i=1}^s{\kappabold_i^{(n)}}^T
                                     \pder{\Rbm_i^{(n)}}{\ubm^{(n-1)}}\right]
                        \pder{\ubm^{(n-1)}}{\mubold} \\
         &+ \sum_{n=1}^{N_t} \sum_{p=1}^s \left[\pder{F}{\kbm_p^{(n)}} -
               {\lambdabold^{(n)}}^T\pder{\tilde\rbm^{(n)}}{\kbm_p^{(n)}}-
                \sum_{i=p}^s {\kappabold_i^{(n)}}^T
                                   \pder{\Rbm_i^{(n)}}{\kbm_p^{(n)}}
                                   \right]\pder{\kbm_p^{(n)}}{\mubold}.
   \end{aligned}
\end{equation}
The dual variables, $\lambdabold^{(n)}$ and $\kappabold_i^{(n)}$, which have
remained arbitrary to this point, are chosen as the solution to the
following equations
\begin{equation} \label{eqn:uns-disc-adj}
  \begin{aligned}
    \pder{\tilde\rbm^{(N_t)}}{\ubm^{(N_t)}}^T\lambdabold^{(N_t)} &=
                                                     \pder{F}{\ubm^{(N_t)}}^T \\
     \pder{\tilde\rbm^{(n)}}{\ubm^{(n-1)}}^T\lambdabold^{(n)}
   + \pder{\tilde\rbm^{(n-1)}}{\ubm^{(n-1)}}^T\lambdabold^{(n-1)} &=
     \pder{F}{\ubm^{(n-1)}}^T - \sum_{i=1}^s \pder{\Rbm_i^{(n)}}{\ubm^{(n-1)}}^T
                                                            \kappabold_i^{(n)}\\
     \sum_{j=i}^s \pder{\Rbm_j^{(n)}}{\kbm_i^{(n)}}^T\kappabold_j^{(n)} &=
      \pder{F}{\kbm_i^{(n)}}^T - \pder{\tilde\rbm^{(n)}}{\kbm_i^{(n)}}^T
                                                               \lambdabold^{(n)}
  \end{aligned}
\end{equation}
for $n = 1, \dots, N_t$ and $i = 1, \dots, s$.  These are the \emph{fully 
discrete adjoint equations} corresponding to the primal evolution equations in
(\ref{eqn:abs-govern-disc}) and quantity of interest $F$. Defining the dual
variables as the solution of the adjoint equations in
(\ref{eqn:uns-disc-adj}), the expression for
$\displaystyle{\oder{F}{\mubold}}$ in (\ref{eqn:uns-disc-adj-derive-2})
reduces to
\begin{equation} \label{eqn:funcl-grad-nosens}
  \oder{F}{\mubold} = \pder{F}{\mubold} -
        \sum_{n=0}^{N_t}{\lambdabold^{(n)}}^T\pder{\tilde\rbm^{(n)}}{\mubold} - 
        \sum_{n=1}^{N_t}\sum_{p=1}^s {\kappabold_p^{(n)}}^T
                                                  \pder{\Rbm_p^{(n)}}{\mubold},
\end{equation}
which is \emph{independent} of the state sensitivities.  Finally, elimination
of the auxiliary terms, $\tilde\rbm^{(n)}$ and $\Rbm_i^{(n)}$, in
equations (\ref{eqn:uns-disc-adj-derive-2}) and
(\ref{eqn:uns-disc-adj}) through differentiation of their
expressions in (\ref{eqn:abs-govern-disc}) gives rise to the adjoint equations
\begin{equation} \label{eqn:uns-disc-adj-dirk}
  \begin{aligned}
    \lambdabold^{(N_t)} &= \pder{F}{\ubm^{(N_t)}}^T \\
    \lambdabold^{(n-1)} &= \lambdabold^{(n)} + \pder{F}{\ubm^{(n-1)}}^T +
      \sum_{i=1}^s \Delta t_n\pder{\rbm}{\ubm}\left(\ubm_i^{(n)},~\mubold,~
                             t_{n-1}+c_i\Delta t_n\right)^T\kappabold_i^{(n)} \\
    \Mbb^T\kappabold_i^{(n)} &= \pder{F}{\kbm_i^{(n)}}^T +b_i\lambdabold^{(n)} +
                         \sum_{j=i}^s a_{ji}\Delta t_n\pder{\rbm}{\ubm}
                         \left(\ubm_j^{(n)},
                         ~\mubold,~t_{n-1}+c_j\Delta t_n\right)^T
                         \kappabold_j^{(n)}
  \end{aligned}
\end{equation}
for $n = 1, \dots, N_t$ and $i = 1, \dots, s$ and the expression for gradient
reconstruction, independent of state sensitivities,
\begin{equation} \label{eqn:funcl-grad-nosens-dirk}
  \oder{F}{\mubold} = \pder{F}{\mubold} +
                     {\lambdabold^{(0)}}^T\pder{\ubm_0}{\mubold} +
                     \sum_{n=1}^{N_t} \Delta t_n \sum_{i=1}^s
                      {\kappabold_i^{(n)}}^T\pder{\rbm}{\mubold}(\ubm_i^{(n)},~
                                                \mubold,~t_{n-1}+c_i\Delta t_n),
\end{equation}
specialized to the case of a DIRK temporal discretization. From inspection
of (\ref{eqn:funcl-grad-nosens-dirk}), it is clear that the initial condition
sensitivity $\displaystyle{\pder{\ubm_0}{\mubold}}$ is the only sensitivity
term required to reconstruct $\displaystyle{\oder{F}{\mubold}}$. The
presence of this term does not destroy the efficiency of the adjoint method
for two reasons:%
\begin{inparaenum}[(a)]
 \item only matrix-vector products with
       $\displaystyle{\pder{\ubm_0}{\mubold}^T}$ are required and
 \item the parametrization of the initial condition is either known
       analytically (uniform flow, zero freestream, independent of
       $\mubold$, etc) or is the solution of some nonlinear system of
       equations (most likely the steady-state equations).
\end{inparaenum}
In the first case, $\displaystyle{{\lambdabold^{(0)}}^T\pder{\ubm_0}{\mubold}}$
can be computed analytically once $\lambdabold^{(0)}$ is known. The next
section details efficient computation of
$\displaystyle{{\lambdabold^{(0)}}^T\pder{\ubm_0}{\mubold}}$ using the
adjoint method of the steady-state problem.

\subsection{Parametrization of Initial Condition}
Suppose the initial condition $\ubm_0(\mubold)$ is defined as the solution
of the nonlinear system of equations -- whose Jacobian is invertible at
$\ubm_0(\mubold)$ -- which is most likely the fully discrete steady-state
form of the governing equations
\begin{equation}
\Rbm(\ubm_0(\mubold), \mubold) = 0.
\end{equation}
Differentiating with respect to the parameter $\mubold$ leads to the expansion 
\begin{equation}
 \oder{\Rbm}{\mubold} = \pder{\Rbm}{\mubold} +
                        \pder{\Rbm}{\ubm_0}\pder{\ubm_0}{\mubold} = 0,
\end{equation}
where arguments have been dropped for brevity.
Assuming the Jacobian matrix is invertible, multiply the preceding equation
by the $\displaystyle{\lambdabold^{(0)}}$ and rearrange to obtain
\begin{equation}
 -{\lambdabold^{(0)}}^T\pder{\ubm_0}{\mubold} =
   \left[\pder{\Rbm}{\ubm_0}^{-T}\lambdabold^{(0)}\right]^T\pder{\Rbm}{\mubold}.
\end{equation}
This reveals the term
$\displaystyle{{\lambdabold^{(0)}}^T\pder{\ubm_0}{\mubold}}$
can be computed at the cost of one linear system solve of the form
$\displaystyle{\pder{\Rbm}{\ubm_0}^T\vbm = \lambdabold^{(0)}}$ and an
inner product $\displaystyle{\vbm^T\pder{\Rbm}{\mubold}}$. The only operation
whose cost scales with the size of $\mubold$ is the evaluation of
$\displaystyle{\pder{\Rbm}{\mubold}}$ and subsequent inner product. Given
this exposition on the fully discrete, time-dependent adjoint method and 
the discrete adjoint method for computing
$\displaystyle{{\lambdabold^{(0)}}^T\pder{\ubm_0}{\mubold}}$, a discussion
is provided detailing the advantages of the fully discrete framework when
computing gradients of output quantities before discussing implementation
details in Section~\ref{sec:implement}.

\subsection{Benefits of Fully Discrete Framework}
In the context of optimization, the
fully discrete adjoint method is advantageous compared to the continuous or
semi-discrete version as it is guaranteed that the resulting derivatives
will be consistent with the quantity of interest, $F$.
This emanates from the fact that in the fully discrete setting, the
\emph{discretization errors are also differentiated}.  This property is
practically relevant as convergence guarantees and convergence rates of
many black-box optimizers are heavily dependent on consistent gradients
of optimization functionals.

Additionally, when Runge-Kutta schemes are chosen for the temporal
discretization, the fully discrete framework is particularly advantageous since
the \emph{stages} are rarely invariant with respect to the direction of
time, that is to say,
\begin{equation} \label{eqn:rk-noninvar}
  \not\exists i, j \in \{1, \dots, s\} \quad \text{such that} \quad
  t_{n-1} + c_i\Delta t_n = t_n - c_j\Delta t_n,
\end{equation}
where $c$ is from the Butcher tableau.  Temporal invariance of a Runge-Kutta
scheme, as defined in (\ref{eqn:rk-noninvar}) is significant when computing
adjoint variables. During the primal solve, $\ubm$ will be
computed at $t_n$ for $n = 1, \dots, N$ and its stage values
at $t_{n-1} + c_i\Delta t_n$ for $n = 1, \dots, N$ and $i = 1, \dots, s$. 
If the same RK scheme is applied to integrate the \emph{semi-discrete}
adjoint equations backward in time, the primal solution will be required
at $t_n - c_i\Delta t_n$ for $n = 1, \dots, N$ and $i = 1, \dots, s$.  Due
to condition (\ref{eqn:rk-noninvar}), the solution to the primal problem was
not computed during the forward solve.  Obtaining the primal solution at this
time requires interpolation, which complicates the implementation,
degrades the accuracy of the computed adjoint variables, and destroys discrete
consistency of the computed gradients. This issue does not arise in the fully
discrete setting as only terms computed during the primal solve appear in the
adjoint equations, by construction.

The next section is devoted to detailing an efficient and modular
implementation of the fully discrete adjoint method on deforming domains.

\section{Implementation} \label{sec:implement}
Implementation of the fully discrete adjoint method introduced in
Section~\ref{sec:adj} relies on the computation of the following terms
from the spatial discretization
\begin{equation} \label{eqn:spat-quant1}
  \Mbb, \rbm, \pder{\rbm}{\ubm}, \pder{\rbm}{\ubm}^T, \pder{\rbm}{\mubold},
  f_h, \pder{f_h}{\ubm}, \pder{f_h}{\mubold}.
\end{equation}
Here, $\Mbb$ is the mass matrix of the semi-discrete conservation law,
and $\rbm$ is the spatial residual vector with
derivatives $\displaystyle{\pder{\rbm}{\ubm}}$ (Jacobian) and
$\displaystyle{\pder{\rbm}{\mubold}}$. As in the previous section, $f_h$ is
the discretization of the spatial integral of the output quantity of interest
with derivatives $\displaystyle{\pder{f_h}{\ubm}}$ and
$\displaystyle{\pder{f_h}{\mubold}}$.
The mass matrix, spatial flux, Jacobian of spatial flux, and output quantity
are standard terms required by an implicit solver and will not be considered
further.  The Jacobian transpose is explicitly mentioned as additional
implementational effort is required when performing parallel matrix
transposition. The derivatives with respect to $\mubold$ are rarely required
outside adjoint method computations and will be considered further in
Section~\ref{subsec:implement-partial}. As indicated in
Section~\ref{subsec:govern-spatial}, all relevant derivatives of the
mass matrix are zero since it is independent of time, parameter, and
state variable, which is an artifact of the transformation to a fixed
reference domain.

The parallel implementation of all semi-discrete quantities in
(\ref{eqn:spat-quant1}) is performed using domain decomposition, where each
processor contains a subset of the elements in the mesh, including a halo of
elements to be communicated with neighbors \cite{persson09parallel}. Linear
systems of the form
\begin{equation*}
 \pder{\rbm}{\ubm}\xbold=\bbold \qquad\qquad \pder{\rbm}{\ubm}^T\xbold=\bbold
\end{equation*}
are solved in parallel using a GMRES solver with a block Incomplete-LU (ILU)
preconditioner.

Given the availability of all terms in (\ref{eqn:spat-quant1}), the solution
of the primal problem and integration of the output quantity $F$
is given in Algorithm~\ref{alg:primal}. The solution of the corresponding
fully discrete adjoint equation, and reconstruction of the gradient
of $F$, is given in Algorithm~\ref{alg:dual}.
\begin{algorithm}[htbp]
  \caption{Primal Solution: Functional Evaluation} \label{alg:primal}
  \begin{algorithmic}[1]
    \REQUIRE Initial condition, $\ubm^{(0)}$;
             parameter configuration, $\mubold$
    \ENSURE Integrated output quantity, $F = \Fcal_h^{(N_t)}$, and
            primal state quantities, $\ubm^{(n)}$ and $\kbm_i^{(n)}$ for
            $n = 1, \dots, N_t$ and $i = 1, \dots, s$
    \STATE Initialize: $\Fcal_h^{(0)} = 0$
    \FOR{$n = 1, \dots, N_t$}
      \FOR{$i = 1, \dots, s$}
        \STATE Solve (\ref{eqn:dirk}) for $\kbm_i^{(n)}$
        $$\Mbb\kbm_i^{(n)} = \Delta t\rbm\left(\ubm_i^{(n)},
                                         \mubold, t_{n-1} + c_i\Delta t\right)$$
        where $\ubm_i^{(n)} = \ubm^{(n-1)} + \sum_{j=1}^i a_{ij}\kbm_j^{(n)}$
        \STATE Write $\kbm_i^{(n)}$ to disk
      \ENDFOR
      \STATE Update $\ubm$ according to (\ref{eqn:dirk})
      $$\ubm^{(n)} = \ubm^{(n-1)} + \sum_{i=1}^s b_i \kbm_i^{(n)}$$
      \STATE Update $\Fcal_h$ according to (\ref{eqn:dirk1})
      $$\Fcal_h^{(n)} = \Fcal_h^{(n-1)} + \sum_{i=1}^s b_i f\left(\ubm_i^{(n)},
                                       \mubold, t_{n-1} + c_i\Delta t_n\right)$$

      \STATE Write $\ubm^{(n)}$ to disk
    \ENDFOR
  \end{algorithmic}
\end{algorithm}

\begin{algorithm}[htbp]
  \caption{Dual Solution: Gradient Evaluation} \label{alg:dual}
  \begin{algorithmic}[1]
    \REQUIRE Primal state quantities, $\ubm^{(n)}$ and $\kbm_i^{(n)}$ for
             $n = 1, \dots, N_t$ and $i = 1, \dots, s$; initial condition
             sensitivity, $\displaystyle{\pder{\ubm^{(0)}}{\mubold}}$; parameter
             configuration, $\mubold$
    \ENSURE Gradient of integrated output quantity,
            $\displaystyle{\oder{F}{\mubold}}$, and
            dual state quantities, $\lambdabold^{(n)}$ and
            $\kappabold_i^{(n)}$ for $n = 1, \dots, N_t$ and $i = 1, \dots, s$
    \STATE Read primal solution $\ubm^{(N_t)}$ from disk
    \STATE $\displaystyle{\lambdabold^{(N_t)} = \pder{F}{\ubm^{(N_t)}}^T}$
    \STATE Initial gradient of $F$ with partial derivative and initial
           condition sensitivity
    $$\oder{F}{\mubold} = \pder{F}{\mubold} +
                                   {\lambdabold^{(0)}}^T\pder{\ubm_0}{\mubold}$$
    \FOR{$n = N_t, \dots, 1$}
      \STATE Read primal solution $\ubm^{(n-1)}$ from disk
      \FOR{$i = s, \dots, 1$}
        \STATE Read primal solution $\kbm_i^{(n)}$ from disk
        \STATE Solve (\ref{eqn:uns-disc-adj-dirk}) for $\kappabold_i^{(n)}$
               $$\Mbb^T\kappabold_i^{(n)} = \pder{F}{\kbm_i^{(n)}}^T +
                         b_i\lambdabold^{(n)} +
               \sum_{j=i}^s a_{ji}\Delta t_n
               \pder{\rbm}{\ubm}(\ubm_j^{(n)}, \mubold, t_{n-1}+c_j\Delta t_n)^T
               \kappabold_j^{(n)}$$
      \STATE Update $\displaystyle{\oder{F}{\mubold}}$ according to
             (\ref{eqn:funcl-grad-nosens-dirk})
      $$\oder{F}{\mubold} = \oder{F}{\mubold} +
        \Delta t_n {\kappabold_i^{(n)}}^T\pder{\rbm}{\mubold}(\ubm_i^{(n)},~
                                               \mubold,~t_{n-1}+c_i\Delta t_n)$$
      \ENDFOR
      \STATE Update $\lambdabold$ according to (\ref{eqn:uns-disc-adj-dirk})
      $$\lambdabold^{(n-1)} = \lambdabold^{(n)} + \pder{F}{\ubm^{(n-1)}}^T
                        + \sum_{i=1}^s \Delta t_n\pder{\rbm}{\ubm}(\ubm_i^{(n)},
                                           \mubold, t_i + c_i\Delta t_n)^T
                                           \kappabold_i^{(n)}$$
    \ENDFOR
  \end{algorithmic}
\end{algorithm}

A well-documented implementational issue corresponding to the unsteady adjoint
method pertains to storage and I/O demands. 
The adjoint equations are solved backward in time and require
the solution of the primal problem at each of the corresponding steps/stages.
Therefore, the adjoint computations cannot begin until all primal states have
been computed.  Additionally, this implies all primal states must be stored
since they will be required in \emph{reverse} order during the adjoint
computation.  For most problems, storing all primal states in memory will
be infeasible, requiring disk I/O, which must be performed in parallel to
ensure parallel scaling is maintained.
There have been a number of strategies to minimize the required I/O operations,
such as local-in-time adjoint strategies \cite{yamaleev2010local} and
checkpointing \cite{charpentier2001checkpointing, heimbach2005efficient,
                    heuveline2006online}. 
For the DG-ALE method in this work, the cost of I/O was not significant
compared to the cost of assembly and solving the linearized system of equations.

In this work, the 3DG software \cite{peraire2008compact} was used for the
high-order DG-ALE scheme. The temporal discretization and unsteady adjoint
method were implemented in the Model Order Reduction Testbed (MORTestbed)
\cite{zahr2010mortestbed, zahr2015progressive} code-base, which was used to
wrap 3DG such that all data structures, and thus all parallel capabilities,
were inherited.

\subsection{Partial Derivatives of Residuals and Output Quantities}
\label{subsec:implement-partial}
This section details computation of partial derivatives of the residual,
$\rbm$, and the output quantity, $f_h$, with respect to the parameter $\mubold$.
The DG-ALE discretizations of Section~\ref{subsec:govern-spatial},
with and without GCL augmentation, are considered separately as the implicit
dependence of $\bar\gbold$ on $\mubold$ requires special treatment.

\subsubsection{Without GCL Augmentation}
\label{subsec:implement-partial-nogcl}
When the GCL augmentation is not considered, the dependence of $\rbm$ and $f_h$
on the parameter $\mubold$ is solely due to the domain parametrization.
Therefore, the following expansion of the partial derivatives with respect
to $\mubold$ is exploited
\begin{equation} \label{eqn:partial-nogcl}
  \pder{\rbm}{\mubold} = \pder{\rbm}{\xbm}\pder{\xbm}{\mubold} + 
                         \pder{\rbm}{\dot{\xbm}}\pder{\dot{\xbm}}{\mubold}
  \qquad\qquad
  \pder{f_h}{\mubold} = \pder{f_h}{\xbm}\pder{\xbm}{\mubold} + 
                        \pder{f_h}{\dot{\xbm}}\pder{\dot{\xbm}}{\mubold}
\end{equation}
where $\displaystyle{\pder{\xbm}{\mubold}}$ and
$\displaystyle{\pder{\dot{\xbm}}{\mubold}}$ are determined solely from the
domain parametrization in Section~\ref{subsec:implement-domain} and the terms
\begin{equation} \label{eqn:partial-terms-nogcl}
  \pder{\rbm}{\xbm}, \pder{\rbm}{\dot{\xbm}},
  \pder{f_h}{\xbm} , \pder{f_h}{\dot{\xbm}}
\end{equation}
are determined from the form of the governing equations and spatial
discretization outlined in Section~\ref{sec:govern}. From the expressions in
(\ref{eqn:partial-nogcl}), the terms in (\ref{eqn:partial-terms-nogcl}) are not explicitly required in matrix form, rather matrix-vector products with
$\displaystyle{\pder{\xbm}{\mubold}}$ and
$\displaystyle{\pder{\dot\xbm}{\mubold}}$ from
Section~\ref{subsec:implement-domain} are required.

\subsubsection{With GCL Augmentation}
\label{subsec:implement-partial-gcl}
For the DG-ALE scheme with GCL augmentation, the dependence of $\rbm$ and $f$
on the parameter $\mubold$ arises from two sources, the domain parametrization
and the implicit dependence of $\bar\gbold$ on $\mubold$.
Therefore, the chain rule expansions in (\ref{eqn:partial-nogcl}) must include
an additional term to account for the dependence of $\bar\gbold$ on $\mubold$
\begin{equation}
  \pder{\rbm}{\mubold} = \pder{\rbm}{\xbm}\pder{\xbm}{\mubold} + 
                         \pder{\rbm}{\dot{\xbm}}\pder{\dot{\xbm}}{\mubold} +
                         \pder{\rbm}{\bar\gbold}\pder{\bar\gbold}{\mubold}
  \qquad\qquad
  \pder{f}{\mubold} = \pder{f}{\xbm}\pder{\xbm}{\mubold} + 
                      \pder{f}{\dot{\xbm}}\pder{\dot{\xbm}}{\mubold} +
                      \pder{f}{\bar\gbold}\pder{\bar\gbold}{\mubold}.
\end{equation}
Similar to the previous section, the terms 
$\displaystyle{\pder{\xbm}{\mubold}}$ and
$\displaystyle{\pder{\dot{\xbm}}{\mubold}}$ are determined solely from the
domain parametrization in Section~\ref{subsec:implement-domain} and
\begin{equation}
  \pder{\rbm}{\xbm}, \pder{\rbm}{\dot{\xbm}}, \pder{\rbm}{\bar\gbold},
  \pder{f}{\xbm}, \pder{f}{\dot{\xbm}}, \pder{f}{\bar\gbold}
\end{equation}
are determined from the form of the governing equations and spatial
discretization in Section~\ref{sec:govern}. The only remaining term
$\displaystyle{\pder{\bar\gbold}{\mubold}}$ is defined as the solution of
the following ODE
\begin{equation}\label{eqn:gbar-sens}
  \Mbb_{\bar\gbold} \pder{}{t}\left(\pder{\bar\gbold}{\mubold}\right) =
                           \pder{\rbm_{\bar\gbold}}{\mubold} +
                 \pder{\rbm_{\bar\gbold}}{\bar\gbold}\pder{\bar\gbold}{\mubold}
  = \pder{\rbm_{\bar\gbold}}{\mubold},
\end{equation}
obtained by direct differentiation of (\ref{eqn:semidisc-gcl1}). The last
equality uses the fact that $\rbm_{\bar\gbold}$ is independent of $\bar\gbold$,
which can be deduced from examination of the governing equation for $\bar{g}$
(\ref{eqn:gbar}). Equation (\ref{eqn:gbar-sens}) is discretized with the same
DIRK scheme used for the temporal discretization of the state equation.

\begin{Rem}
The special treatment of $\bar\gbold$ detailed in this section, including
integration of the sensitivity equations (\ref{eqn:gbar-sens}), can be
avoided by considering the ODEs in (\ref{eqn:semidisc-gcl1}) directly
without leveraging the fact that the $\bar\gbold$ equation is independent
of $\ubm_{\bar\Xbm}$. This implies the state vector will contain an
additional unknown for $\bar{g}$ for each DG node. This increases the cost
of a primal and dual solve, but simplifies the adjoint derivation and
implementation, allowing the form in
Section~\ref{subsec:implement-partial-nogcl} to be used directly.
\end{Rem}

\subsection{Time-Dependent, Parametrized Domain Deformation}
\label{subsec:implement-domain}
A crucial component of the fully discrete adjoint method on deforming domains
is a time-dependent parametrization of the domain, amenable to parallel
implementation. A parallel implementation is required as domain deformation
will involve operations on the entire computational mesh and will be queried
at every stage of each time step of both the primal and dual solves, according
to Algorithms~\ref{alg:primal}~and~\ref{alg:dual}.
In this work, the domain parametrization is required to be sufficiently
general to handle shape deformation, as well as kinematic motion.
Additionally, the domain deformation must be sufficiently smooth to ensure
sufficient regularity of the transformed solution, and the spatial and temporal
derivatives must be analytically available for fast, accurate computation of the
deformation gradient, $\Gbm$, and velocity, $\vbm_\Xbm$, of the mapping,
$\Gcal$.

The domain deformation will be defined by the superposition of a 
rigid body motion and a spatially varying deformation. To avoid large
mesh velocities at the far-field, which could arise from rigid rotations of
the body, the blending maps of \cite{persson2009dgdeform} are used.
First, define a spatial configuration consisting of a rigid body motion
($\Qbm(\mubold, t)$, $\vbm(\mubold, t)$) and deformation
($\varphibold(\Xbm, \mubold, t)$) to the \emph{reference} domain
\begin{equation}
\Xbm' = \Qbm(\mubold, t)\Xbm + \vbm(\mubold, t) + \varphibold(\Xbm, \mubold, t),
\end{equation}
which completely defines the \emph{physical} motion of the body.
This physical configuration is blended with the reference configuration
according to
\begin{equation}
\xbm = (1-b(d(\Xbm))) \Xbm' + b(d(\Xbm))\Xbm
\end{equation}
where $d(\Xbm) = \norm{\Xbm - \Xbm_0} - R_0$ is the signed distance
from the origin $\Xbm_0$ to the circle of radius $R_0$ centered at $\Xbm_0$
and
\begin{equation}
b(s) =
\begin{cases}
  0, & \text{if}~s<0\\
  1, & \text{if}~s>R_1\\
  r(s/R_1), & \text{otherwise}
\end{cases}
\end{equation}
where $r(s) = 3s^2-2x^3$ for a cubic blending and $r(s) = 10s^3-15s^4+6s^5$
for a quintic blending. Spatial blending of this form ensures the desired
physical motion of the body, $\Xbm'$ is exactly achieved within a radius $R_0$
of the origin. Further, there is no deformation outside a
radius $R_0+R_1$ of the origin. In the annulus about the origin with inner
radius $R_0$ and outer radius $R_0+R_1$, the spatial configuration is blended
smoothly between these two spatial configurations.

The specific form of $\Qbm(\mubold, t)$, $\vbm(\mubold, t)$, and
$\varphibold(\Xbm, \mubold, t)$ is problem-specific and will be deferred to
Section~\ref{sec:app}. Assuming these terms are known analytically, the specific
form of $\displaystyle{\Gbm = \pder{\xbm}{\Xbm}}$,
$\displaystyle{\vbm_\Xbm = \dot\xbm = \pder{\xbm}{t}}$,
$\displaystyle{\pder{\xbm}{\mubold}}$, and
$\displaystyle{\pder{\dot\xbm}{\mubold}}$ can be easily computed.

\section{Applications}\label{sec:app}
In this section, the high-order numerical discretization of a system of
conservation laws and corresponding adjoint method is applied to the
isentropic compressible Navier-Stokes equations to solve optimal control and
shape optimization problems using gradient-based optimization techniques.
The compressible Navier-Stokes equations take the form
\begin{align}
\frac{\partial \rho}{\partial t}  + \frac{\partial}{\partial x_i}
(\rho u_i) &= 0, \label{ns1}\\
\frac{\partial}{\partial t} (\rho u_i) +
\frac{\partial}{\partial x_i} (\rho u_i u_j+ p)  &=
+\frac{\partial \tau_{ij}}{\partial x_j}
\quad\text{for }i=1,2,3, \label{ns2} \\
\frac{\partial}{\partial t} (\rho E) +
\frac{\partial}{\partial x_i} \left(u_j(\rho E+p)\right) &=
-\frac{\partial q_j}{\partial x_j}
+\frac{\partial}{\partial x_j}(u_j\tau_{ij}), \label{ns3}
\end{align}
in $v(\mubold, t)$ where $\rho$ is the fluid density, $u_1,u_2,u_3$ are the
velocity components, and $E$ is the total energy. The viscous stress tensor and
heat flux are given by
\begin{align}
\tau_{ij} = \mu
\left( \frac{\partial u_i}{\partial x_j} +
\frac{\partial u_j}{\partial x_i} -\frac23
\frac{\partial u_k}{\partial x_k} \delta_{ij} \right)
\qquad \text{ and } \qquad
q_j = -\frac{\mu}{\mathrm{Pr}} \frac{\partial}{\partial x_j}
\left( E+\frac{p}{\rho} -\frac12 u_k u_k \right).
\end{align}
Here, $\mu$ is the viscosity coefficient and $\mathrm{Pr = 0.72}$ is
the Prandtl number which we assume to be constant. For an ideal gas,
the pressure $p$ has the form
\begin{align}
p=(\gamma-1)\rho \left( E - \frac12 u_k u_k\right), \label{ns5}
\end{align}
where $\gamma$ is the adiabatic gas constant.  Furthermore, the entropy of the
system is assumed constant, which is tantamount to the flow being adiabatic
and reversible.  For a perfect gas, the entropy is defined as
\begin{equation}\label{eqn:entropy}
  s = p/\rho^\gamma.
\end{equation}
Using (\ref{eqn:entropy}) to explicitly relate the pressure and density, the
energy equation becomes redundant. This effectively reduces the square system
of PDEs of size $n_{sd}+2$ to one of size $n_{sd}+1$, where $n_{sd}$ is
the number of spatial dimensions.  It can be shown, under suitable
assumptions, that the solution of the isentropic approximation of the
Navier-Stokes equations converges to the solution of the incompressible
Navier-Stokes equations as the Mach number goes to 0
\cite{lin1995incompressible, desjardins1999incompressible, froehle2013high}.

The DG-ALE scheme introduced in Section~\ref{sec:govern} is used for the
spatial discretization of the system of conservation laws with polynomial order
$p = 3$ and a diagonally implicit Runge-Kutta scheme for the temporal
discretization. The DG-ALE scheme uses the Roe flux \cite{roe1981approximate}
for the inviscid numerical flux and the Compact DG flux
\cite{peraire2008compact} for the viscous numerical flux.  The Butcher tableau
for the three-stage, third-order DIRK scheme considered in this work is given
in Table~\ref{tab:dirk3}.

\begin{figure}[h]
  \centering
  \includegraphics[width=2.25in]{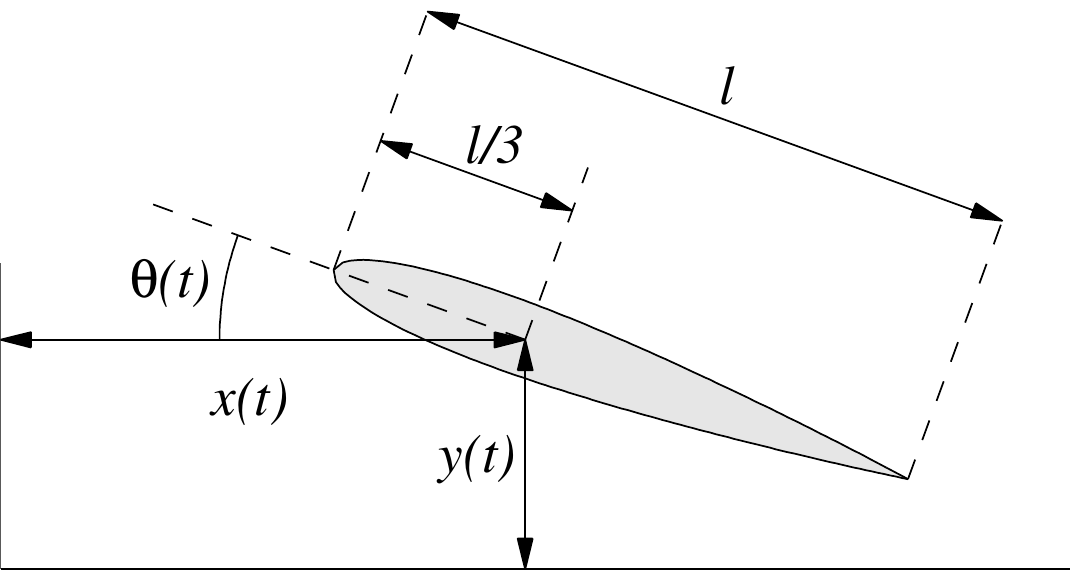}
  \caption{Airfoil}\label{fig:airfoil-traj}
\end{figure}

\begin{table}[h]
  \centering
  \begin{minipage}{0.52\textwidth}
  \centering
  \begin{tabular} { l | c c c}
     $\alpha$ & $\alpha$  &  &  \\
     $\frac{1+ \alpha}{2}$ & $\frac{1+ \alpha}{2}-\alpha$  & $\alpha$ &  \\
     $1$ &  $\gamma$ & $\omega$ & $\alpha$ \\ \hline
     & $\gamma$ & $\omega$ & $\alpha$
  \end{tabular}
  \caption{Butcher Tableau for 3-stage, 3rd order DIRK scheme
  \cite{alexander77dirk} \newline
  $\alpha =  0.435866521508459$, 
  $\gamma = -\frac{6\alpha^2-16\alpha+1}{4}$, 
  $\omega = \frac{6\alpha^2-20\alpha+5}{4}$.} \label{tab:dirk3}
  \end{minipage}
\end{table}
The instantaneous quantities of interest for a body, defined by the surface
$\Gammabold$, take the following form
\begin{equation} \label{eqn:funcls}
 \begin{aligned}
  \Fcal_x(\Ubm, \mubold, t) &= \int_\Gammabold \fbm(\Ubm, \mubold, t) \cdot
                                               \ebm_1~dS \qquad &
  \Fcal_y(\Ubm, \mubold, t) &= \int_\Gammabold \fbm(\Ubm, \mubold, t) \cdot
                                               \ebm_2~dS \\
  \Pcal(\Ubm, \mubold, t) &= \int_\Gammabold \fbm(\Ubm, \mubold, t) \cdot
                                             \dot\xbm~dS \qquad &
  \Pcal_x(\Ubm, \mubold, t) &= \int_\Gammabold
                                  \dot{x}\fbm(\Ubm, \mubold, t)\cdot\ebm_1~dS \\
  \Pcal_y(\Ubm, \mubold, t) &= \int_\Gammabold
                                  \dot{y}\fbm(\Ubm, \mubold, t)\cdot\ebm_2~dS
                         \qquad &
  \Pcal_\theta(\Ubm, \mubold, t) &= -\int_\Gammabold
                                               \dot\theta\fbm(\Ubm, \mubold, t)
                                               \times (\xbm - \xbm_0)~dS
 \end{aligned}
\end{equation}
where $\fbm \in \Rbb^{n_{sd}}$ is the force imparted by the fluid on the body,
$\ebm_i$ is the $i$th canonical basis vector in $\Rbb^{n_{sd}}$,
$\xbm$ and $\dot\xbm$ are the position and velocity of a point on the
surface $\Gammabold$, and $x$, $y$, $\theta$, $\dot{x}$, $\dot{y}$, $\dot\theta$
define the motion of the reference point, $\xbm_0$ (the $1/3$-chord of the
airfoil, in this case); see Figure~\ref{fig:airfoil-traj}. The $\Fcal_x$ and
$\Fcal_y$ terms correspond to the total $x$- and $y$-directed forces on the
body and $\Pcal$ is the total power exerted on the body by the fluid. The
total power $\Pcal$ is broken into its translational, $\Pcal_x$ and $\Pcal_y$,
and rotational, $\Pcal_\theta$, components. For a 2D rigid body motion, an
additive relationship among these terms holds
\begin{equation} \label{eqn:works}
 \Pcal(\Ubm, \mubold, t) = \Pcal_x(\Ubm, \mubold, t) +
                           \Pcal_y(\Ubm, \mubold, t) +
                           \Pcal_\theta(\Ubm, \mubold, t).
\end{equation}
The negative sign is included in the definition of $\Pcal_\theta$ due to the
clockwise definition of $\theta$ in Figure~\ref{fig:airfoil-traj}.
In the remainder of this document, a superscript $h$ will be used to denote the
high-order DG approximation to these spatial integrals that constitute the
instantaneous quantities of interest, e.g., $\Pcal^h(\ubm, \mubold, t)$ is the
high-order approximation of $\Pcal(\Ubm, \mubold, t)$, where $\ubm$ is the
semi-discrete approximation of $\Ubm$. Temporal integration of
the instantaneous quantities of interest leads to the integrated quantities of
interest
\begin{equation} \label{eqn:funcls}
 \begin{aligned}
  \Jcal_x(\Ubm, \mubold) &= \int_0^T\int_\Gammabold \fbm(\Ubm, \mubold, t) \cdot
                                                    \ebm_1~dS~dt \qquad &
  \Jcal_y(\Ubm, \mubold) &= \int_0^T\int_\Gammabold \fbm(\Ubm, \mubold, t) \cdot
                                                    \ebm_2~dS~dt \\
  \Wcal(\Ubm, \mubold) &= \int_0^T\int_\Gammabold \fbm(\Ubm, \mubold, t) \cdot
                                                  \dot\xbm~dS~dt \qquad &
  \Wcal_x(\Ubm, \mubold) &= \int_0^T\int_\Gammabold 
                               \dot{x}\fbm(\Ubm, \mubold, t)\cdot\ebm_1~dS~dt \\
  \Wcal_y(\Ubm, \mubold) &= \int_0^T\int_\Gammabold
                         \dot{y}\fbm(\Ubm, \mubold, t)\cdot\ebm_2~dS~dt \qquad &
  \Wcal_\theta(\Ubm, \mubold) &= -\int_0^T\int_\Gammabold
                               \dot\theta\fbm(\Ubm, \mubold, t)
                                                    \times (\xbm - \xbm_0)~dS~dt
 \end{aligned}
\end{equation}
which will be used as optimization functionals in subsequent sections. 
The terms $\Jcal_x$ and $\Jcal_y$ are the $x$- and $y$-directed impulse the
fluid exerts on the airfoil, respectively, $\Wcal$ is the total work done on
the airfoil by the fluid, and $\Wcal_x$, $\Wcal_y$, and $\Wcal_\theta$ are 
the translational and rotational components of the total work. The
fully discrete, high-order approximation of the integrated quantities of
interest (DG in space, DIRK in time) will be denoted with the corresponding
Roman symbol, e.g.,
$W(\ubm^{(0)}, \dots, \ubm^{(N_t)}, \kbm_1^{(n)}, \dots, \kbm_s^{(n)}, \mubold)$
is the fully discrete approximation of $\Wcal(\Ubm, \mubold)$.

\subsection{Energetically Optimal Trajectory of 2D Airfoil in Compressible,
            Viscous Flow} \label{subsec:app-pitch}
\begin{figure}[!htbp]
 \centering
 \begin{subfigure}{0.49\textwidth}
   \centering
   \input{tikz/pitch-traj-x.tikz}
 \end{subfigure} \hfill
 \begin{subfigure}{0.49\textwidth}
   \centering
   \input{tikz/pitch-traj-y.tikz}
 \end{subfigure} \\
 \begin{subfigure}{0.49\textwidth}
   \centering
   \input{tikz/pitch-traj-th.tikz}
 \end{subfigure}
 \caption{Trajectories of $x(t)$, $y(t)$, and $\theta(t)$ at initial guess
          (\ref{line:pitch-init}), solution of (\ref{opt:pitch-disc}) under
          parametrization PI (\ref{line:pitch-opt1}), and solution of
          (\ref{opt:pitch-disc}) under parametrization PII
          (\ref{line:pitch-opt2}) for ID = $7$.}
 \label{fig:pitch-traj-kine}
\end{figure}
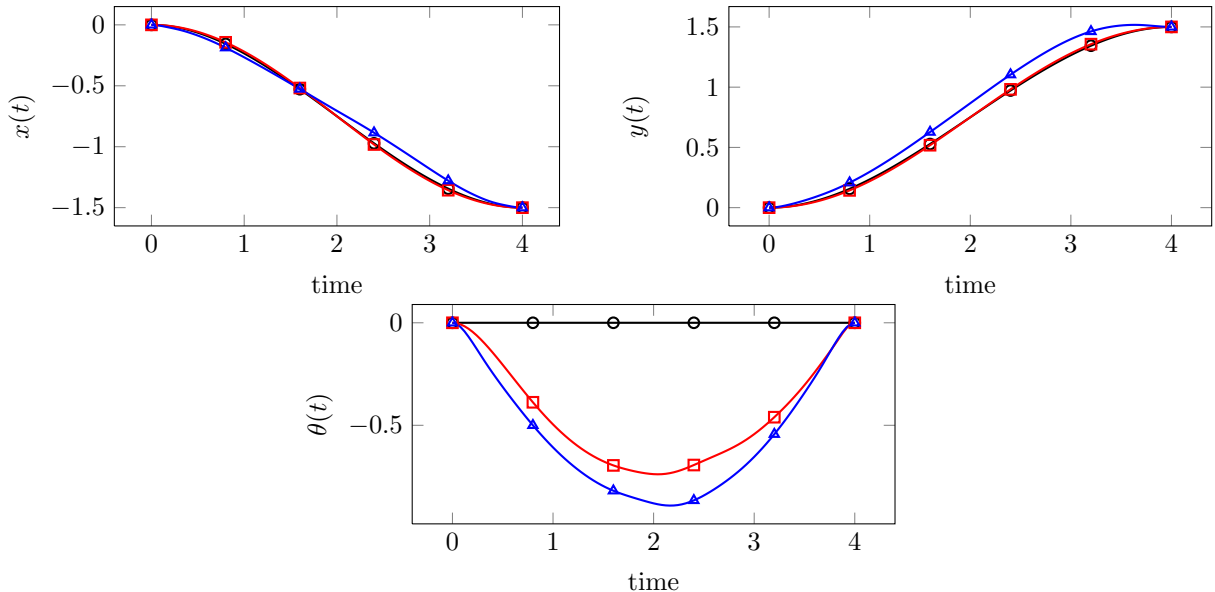

In this section, the high-order, time-dependent PDE-constrained optimization
framework introduced in this document is applied to find the
energetically optimal trajectory of a 2D NACA0012 airfoil with
chord length $l = 1$ and zero-thickness trailing edge. The governing
equations are the 2D compressible, isentropic Navier-Stokes equations.
The mission of the airfoil is to move a distance of $-1.5$ units horizontally
and $1.5$ units vertically in $T = 4$ units of time, with the restriction that
$\theta(0) = \theta(T) = 0$, i.e. the angle of attack at the initial and final
time is zero. Additionally, to ensure smoothness of the motion and avoid
non-physical transients,
$\dot{x}(0) = \dot{x}(T) = \dot{y}(0) = \dot{y}(T) = \dot\theta(0) =
 \dot\theta(T) = 0$ are enforced.
The goal is to determine the trajectory $x(t)$, $y(t)$, $\theta(t)$ of the
airfoil that minimizes the total energy required to complete the mission,
i.e.,
\begin{equation} \label{opt:pitch-cont}
  \begin{aligned}
    & \underset{\Ubm,~\mubold}{\text{minimize}}
    & & \Wcal(\Ubm, \mubold) \\
    & \text{subject to}
    & &   x(0) = \dot{x}(0) = \dot{x}(T) = 0,~x(T) = -1.5 \\
    & & & y(0) = \dot{y}(0) = \dot{y}(T) = 0,~y(T) =  1.5 \\
    & & & \theta(0) = \theta(T) = \dot\theta(0) = \dot\theta(T) = 0 \\
    & & & \pder{\Ubm}{t} + \nabla \cdot \Fbm(\Ubm, \nabla\Ubm) = 0 \quad
          \text{ in } v(\mubold, t).
  \end{aligned}
\end{equation}
The trajectory of the airfoil -- $x(t)$, $y(t)$, and $\theta(t)$ -- is
discretized via clamped cubic splines with $m_x+1$, $m_y+1$, and $m_\theta+1$
knots, respectively. The knots are uniformly spaced between $0$ and $T$ in
the $t$-dimension and the knot values are optimization parameters.
\begin{table}[!htbp]
 \centering
 \input{dat/pitch-setup.tab}
 \caption{Summary of parametrizations considered in
          Section~\ref{subsec:app-pitch}. The number of clamped cubic spline
          knots used to discretize $x(t)$, $y(t)$, and $\theta(t)$ are $m_x+1$,
          $m_y+1$, and $m_\theta+1$, respectively. PI freezes the rigid
          body translation ($m_x = m_y = 0$) and optimizes over only the
          rotation ($m_\theta \neq 0$). PII optimizes over all rigid body
          degrees of freedom ($m_x = m_y = m_\theta \neq 0$).}
\label{tab:pitch-setup}
\end{table}
Table~\ref{tab:pitch-setup} summarizes two parametrizations considered in this
section:
\begin{inparaenum}
 \item[(PI)]  the translational degrees of freedom -- $x(t)$ and $y(t)$ -- are
              frozen at their nominal value in Figure~\ref{fig:pitch-traj-kine}
              and the rotational degree of freedom -- $\theta(t)$ -- is
              parametrized with a $m_\theta+1$-knot clamped cubic spline and
 \item[(PII)] all rigid body modes are parametrized with clamped cubic splines.
\end{inparaenum}
The $7$ IDs in Table~\ref{tab:pitch-setup} correspond to levels of refinement
of the given parametrization with ID = 1 being the coarsest parametrization and
ID = 7 the finest. With this parametrization of the airfoil kinematics, spatial
and temporal discretization with the high-order scheme of
Section~\ref{sec:govern} leads to the  fully discrete version of the
optimization problem in (\ref{opt:pitch-cont})
\begin{equation} \label{opt:pitch-disc}
  \begin{aligned}
    & \underset{\substack{\ubm^{(0)},~\dots,~\ubm^{(N_t)} \in \Rbb^{N_\ubm},\\
                        \kbm_1^{(1)},~\dots,~\kbm_s^{(N_t)} \in \Rbb^{N_\ubm},\\
                        \mubold \in \Rbb^{N_\mubold}}}
                        {\text{minimize}}
    & & W(\ubm^{(0)},~\dots,~\ubm^{(N_t)},~\kbm_1^{(1)},~\dots,~\kbm_s^{(N_t)},~
          \mubold) \\
    & \text{subject to}
    & &   x(0) = \dot{x}(0) = \dot{x}(T) = 0,~x(T) = -1.5 \\
    & & & y(0) = \dot{y}(0) = \dot{y}(T) = 0,~y(T) =  1.5 \\
    & & & \theta(0) = \theta(T) = \dot\theta(0) = \dot\theta(T) = 0 \\
    & & &  \ubm^{(0)} = \ubm_0 \\
    & & & \ubm^{(n)} = \ubm^{(n-1)} + \sum_{i = 1}^s b_i\kbm^{(n)}_i \\
    & & & \Mbb\kbm^{(n)}_i = \Delta t_n\rbm\left(\ubm_i^{(n)},~\mubold,~
                                                 t_{n-1} + c_i\Delta t_n\right).
  \end{aligned}
\end{equation}

\begin{figure}[!htbp]
 \centering
 \input{tikz/pitch-findiff-rev.tikz}
 \caption{Verification of adjoint-based gradient with fourth-order
          centered finite difference approximation, for a range of finite
          intervals, $\tau$, for the total work $W$ -- the objective function in
          (\ref{opt:pitch-disc}) -- for parametrization PII
          (Table~\ref{tab:pitch-setup}). The computed gradient match the
          finite difference approximation to about $10$ digits of
          accuracy before round-off errors degrade the accuracy.}
 \label{fig:pitch-findiff}
\end{figure}
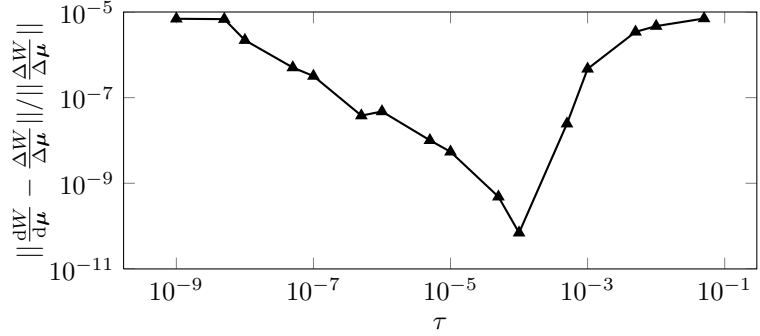

Before considering the optimization problem (\ref{opt:pitch-disc}), the proposed
adjoint method for computing gradients of quantities of interest on the manifold
of fully discrete, high-order solutions of the conservation law (\ref{eqn:dirk})
is verified against a fourth-order finite difference approximation. The finite
difference approximation to gradients on the aforementioned manifold requires
finding the solution of the fully discretized governing equations
\emph{at perturbations} about the nominal parameter configuration in
Figure~\ref{fig:pitch-traj-kine}. To mitigate round-offs errors as much as
possible in the finite difference computation, the number of time steps was
reduced to $10$ and only half the trajectory was simulated.
Figure~\ref{fig:pitch-findiff} shows the relative error between the gradients
computed via the adjoint method and this finite difference approximation for
a sweep of finite difference intervals, $\tau$. A relative error of nearly
$10^{-11}$ is observed for a finite difference step of $\tau = 10^{-4}$.
As expected, the error starts to increase after $\tau$ drops too small due
to the trade-off between finite difference accuracy and roundoff error.

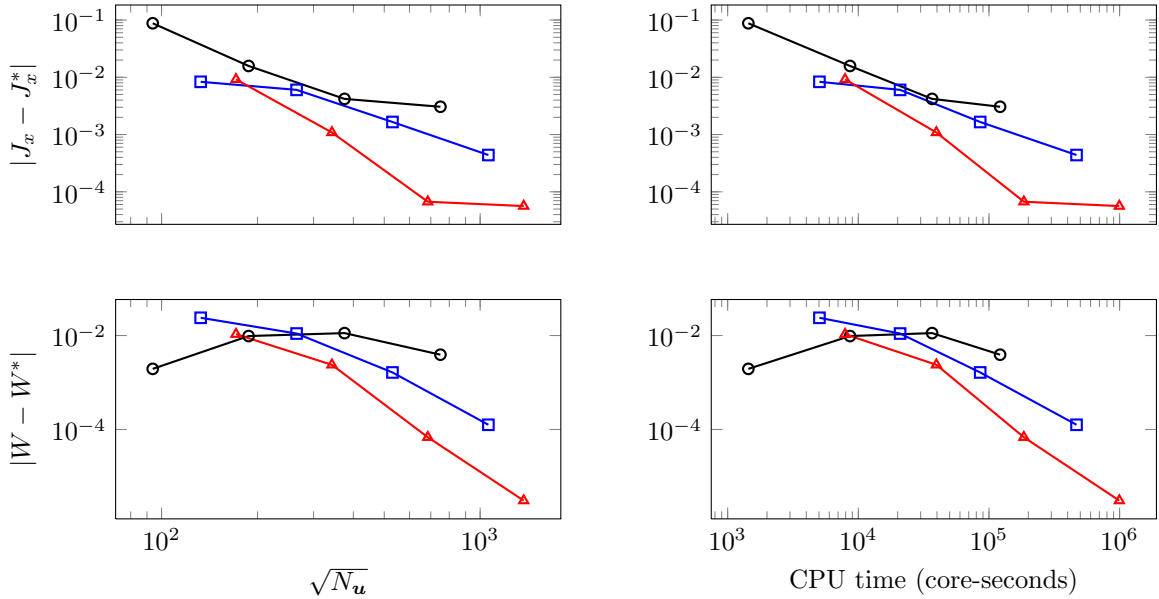
\begin{figure}
 \input{tikz/pitch-spatial-conv-qoi.tikz}
 \caption{Spatial convergence in quantities of interest for the DG-ALE-DIRK
          discretization for $p = 1$ (\ref{line:p1}), $p = 2$ (\ref{line:p2}),
          $p = 3$ (\ref{line:p3}). All schemes use the DIRK3 temporal
          discretization (Table~\ref{tab:dirk3}) with $1000$ time steps to
          mitigate temporal error. The reference values $W^*$ and $J_x^*$ are
          produced from a simulation consisting of $250368$ $p = 4$ elements
          and $1000$ DIRK3 time steps.}
 \label{fig:unnecessary-conv-plt-for-reviewer2}
\end{figure}

A brief study of the convergence behavior of the underlying DG scheme is
provided in Figure~\ref{fig:unnecessary-conv-plt-for-reviewer2} to justify our
use of $p = 3$. In this study, no attempt is made to quantify convergence rates
due to the sharp trailing edge of the NACA0012 that limits regularity of the
solution. From Figure~\ref{fig:unnecessary-conv-plt-for-reviewer2}, the higher
order methods outperform the lower order ones in terms of error -- in the
work, $W$, and $x$-impulse, $J_x$ -- for a given number of degrees of freedom
($N_\ubm$) and CPU time. All timings were obtained using 12 CPU cores.

With this verification of the adjoint-based gradients and high-order accuracy,
attention is turned to the optimization problem in (\ref{opt:pitch-disc}). The
optimization solver
used in this section is L-BFGS-B \citep{zhu1997algorithm}, a bound-constrained,
limited-memory BFGS algorithm. Figure~\ref{fig:pitch-traj-kine} contains the
initial guess for the optimization problem in (\ref{opt:pitch-disc}) as well
as its solution under both parametrization, PI and PII, at the finest level of
refinement (ID = 7). The initial guess for the optimization problem is a pure
translational motion with $\theta(t) = 0$. The solution under parametrization
PI freezes the translational motion at its nominal value and incorporates
rotational motion. The solution under parametrization PII increases the
amplitude of the rotation, flattens the trajectory of $x(t)$, and incorporates
an overshoot in $y(t)$ before settling to the required location, as compared
to the optimal solution corresponding to PI.

\begin{figure}[!htbp]
 \centering
 \input{tikz/pitch-traj-qoi.tikz}
 \caption{Time history of instantaneous quantities of interest ($x$-directed
          force -- $\Fcal_x^h(\ubm, \mubold, t)$, $y$-directed force  --
          $\Fcal_y^h(\ubm, \mubold, t)$, total power --
          $\Pcal^h(\ubm, \mubold, t)$, $x$-translational power --
          $\Pcal_x^h(\ubm, \mubold, t)$, $y$-translational power --
          $\Pcal_y^h(\ubm, \mubold, t)$, rotational power --
          $\Pcal_\theta^h(\ubm, \mubold, t)$) at initial guess
          (\ref{line:pitch-init}), solution of (\ref{opt:pitch-disc}) under
          parametrization PI (\ref{line:pitch-opt1}), and solution of
          (\ref{opt:pitch-disc}) under parametrization PII
          (\ref{line:pitch-opt2}) for ID = $7$.}
 \label{fig:pitch-traj-qoi}
\end{figure}
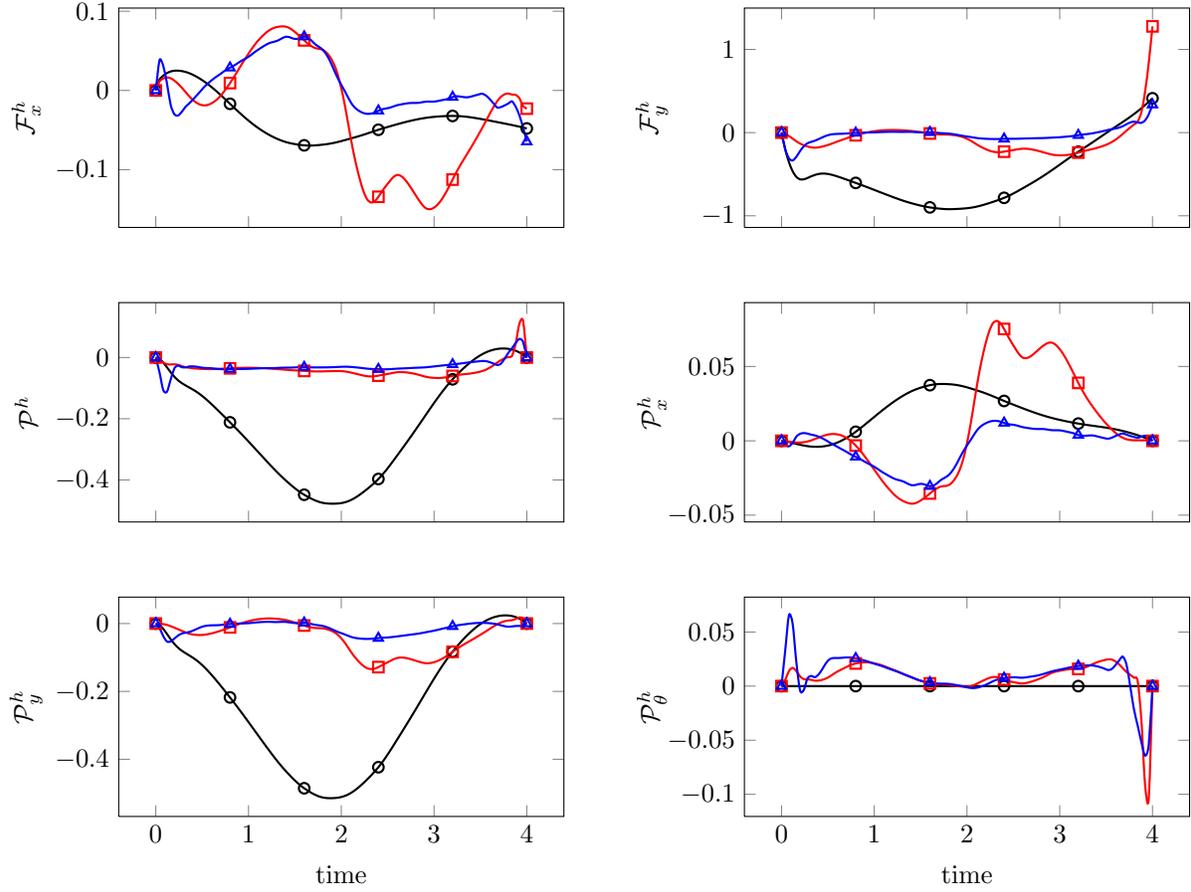

The instantaneous quantities of interest for the nominal trajectory and solution
of (\ref{opt:pitch-disc}) under parametrizations PI and PII are included in
Figure~\ref{fig:pitch-traj-qoi}. It is clear that the optimal solution under
both parametrizations result in a time history of the total power that is
uniformly closer to $0$ than that at the nominal trajectory, which is expected
since $W$ is the objective function. With the exception of the edges of the
time interval, the total power time history for the optimal solution under
parametrization PII is uniformly closer to $0$ than that of PI. The same
observation holds for the power due to the translational motion, $\Pcal_x^h$
and $\Pcal_y^h$. Whereas the total power corresponding to the nominal trajectory
is due solely to the translational motion (since there is no rotation), the
optimal solutions exchange large amounts of translational power for a small
amount of rotational power. These observations can also be verified in
Table~\ref{tab:pitch-summ} which summarizes the optimal values of the integrated
quantities of interest.

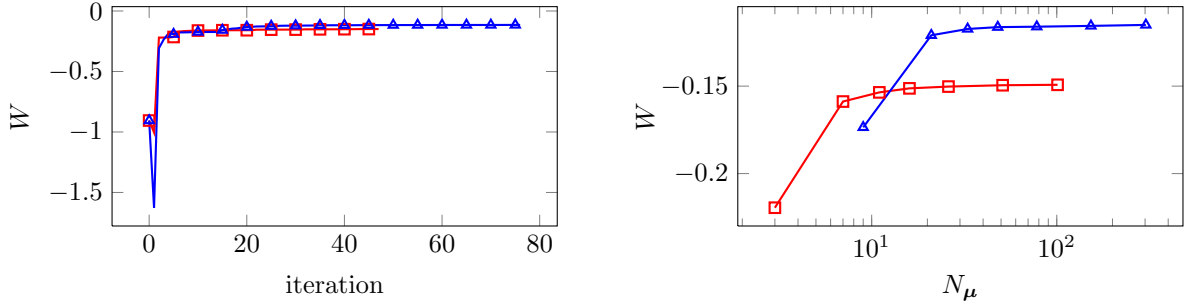
\begin{figure}[!htbp]
 \centering
 \input{tikz/pitch-conv.tikz}
 \caption{\emph{Left}: Convergence of total work $W$ with optimization iteration
          for parametrization PI (\ref{line:pitch-opt1}) and PII
          (\ref{line:pitch-opt2}) for ID = $7$. Both optimization problems
          converge to a motion with significantly lower required total work;
          PII finds a better motion than PI (in terms of total work) due to the
          enlarged search space, at the cost of additional iterations. Each
          optimization iteration requires a primal flow computation -- to
          evaluate the quantities of interest -- and its corresponding adjoint
          -- to evaluate the gradient of the quantity of interest. \emph{Right}:
          Convergence of optimal value of total work $W$ as parameter space is
          refined for parametrization PI (\ref{line:pitch-opt1}) and PII
          (\ref{line:pitch-opt2}). This implies convergence to an optimal,
          smooth trajectory that is not polluted by its discrete
          parametrization.}
 \label{fig:pitch-conv}
\end{figure}

The convergence of the total work, i.e., the objective function of the
optimization problem, with iterations of the optimization solver is summarized
in Figure~\ref{fig:pitch-conv} (left). Both parametrizations are included
and iterations are agglomerated over all IDs. The first iteration corresponds to
a steepest descent step, which causes an adverse jump in the objective value.
The following iterations make rapid progress toward the optimal solution,
which is slowed as convergence is approached. The solver requires additional
iterations to converge the solution corresponding to parametrization PII, which
is expected due to the larger parameter space.

Next, convergence of the total work as the parameter space is refined is
considered in Figure~\ref{fig:pitch-conv} (right) and
Table~\ref{tab:pitch-summ}. This implies the optimal trajectory among all
twice continuously differentiable functions is being approached. For both
parametrizations, the optimal value of the total work agrees to $3$ digits
between IDs $6$ and $7$ (roughly a factor of $2$ difference in dimension of
parameter spaces) and $2$ digits between IDs $3$ and $7$ (roughly a factor of
$10$ difference in dimension of parameter spaces).
This sweep of optimization problems was repeated with a low-order method
($p = 1$) with a comparable number of degrees of freedom ($N_\ubm$). As
expected, the same local optimum was obtained and the convergence behavior with
respect to refinement in $N_\mubold$ was very similar.

\begin{table}[!htbp]
 \centering
 \input{dat/pitch-summ.tab}
 \caption{Table summarizing integrated quantities of interest at optimal
          solution of (\ref{opt:pitch-disc}) for each parametrization (PI, PII)
          for each level of refinement. The total work monotonically increases
          as $N_\mubold$ increases for a given parametrization, which is
          expected due to the nested search spaces. For a fixed ID, the optimal
          total work for parametrization PII is larger than that for PI since
          the search space for PI is a subset of that of PII. The other
          integrated  quantities are included for completeness, but do not
          exhibit trends (except for converging to a fixed value as $N_\mubold$
          increases) since they were not included in the optimization problem.}
 \label{tab:pitch-summ}
\end{table}

\begin{figure}[!htbp]
 \centering
 \begin{subfigure}{0.3\textwidth}
  \centering
  \includegraphics[width=\textwidth]{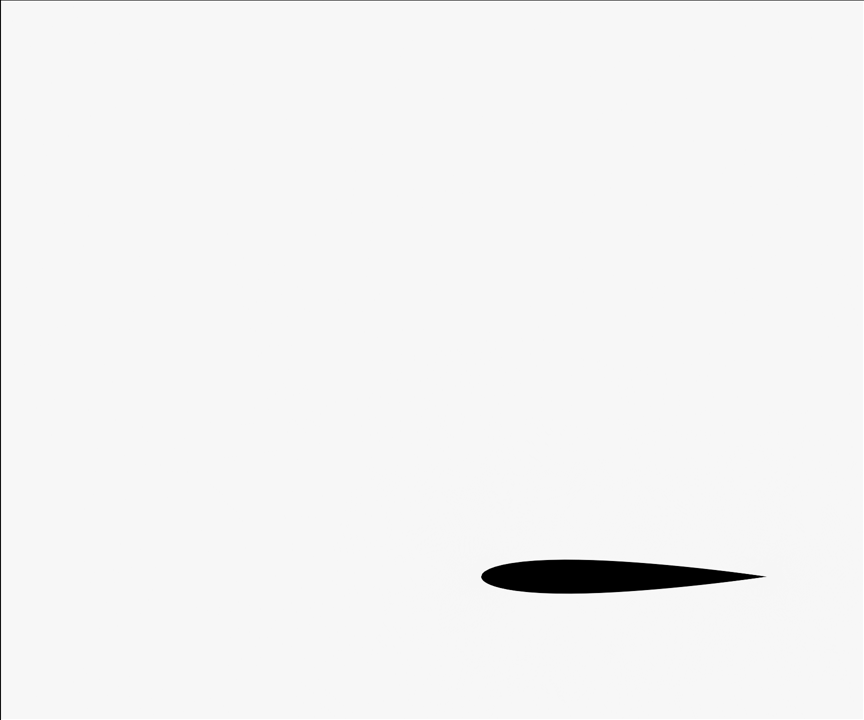}
 \end{subfigure} \hfill
 \begin{subfigure}{0.3\textwidth}
  \centering
  \includegraphics[width=\textwidth]{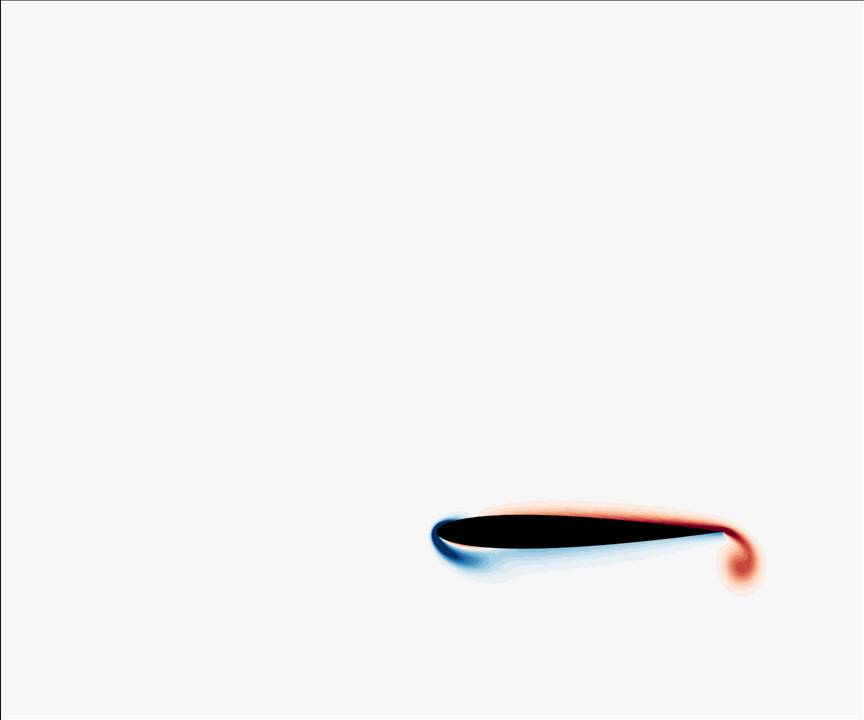}
 \end{subfigure} \hfill
 \begin{subfigure}{0.3\textwidth}
  \centering
  \includegraphics[width=\textwidth]{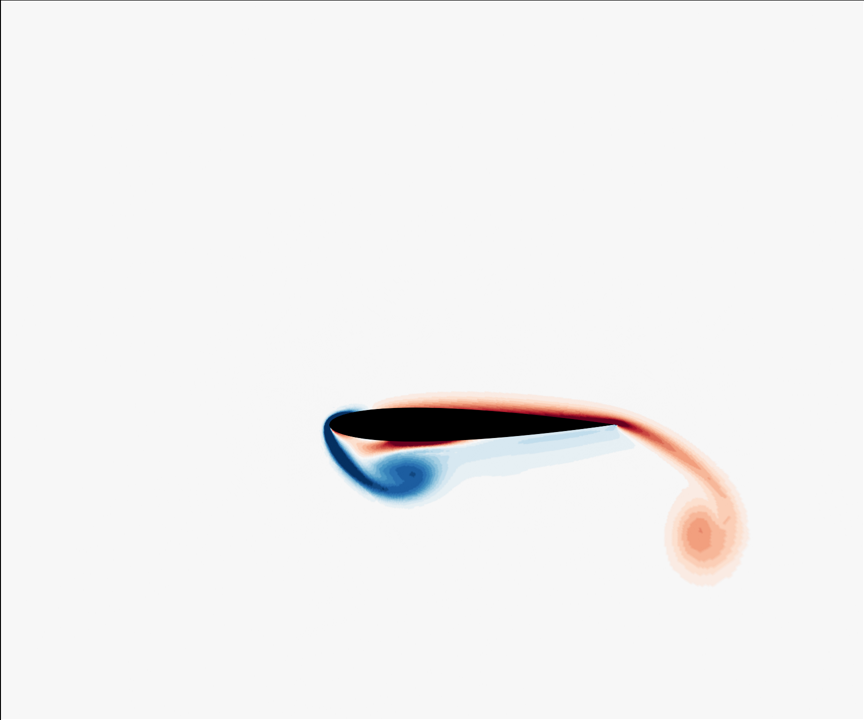}
 \end{subfigure} \\ \vspace{3mm}
 \begin{subfigure}{0.3\textwidth}
  \centering
  \includegraphics[width=\textwidth]{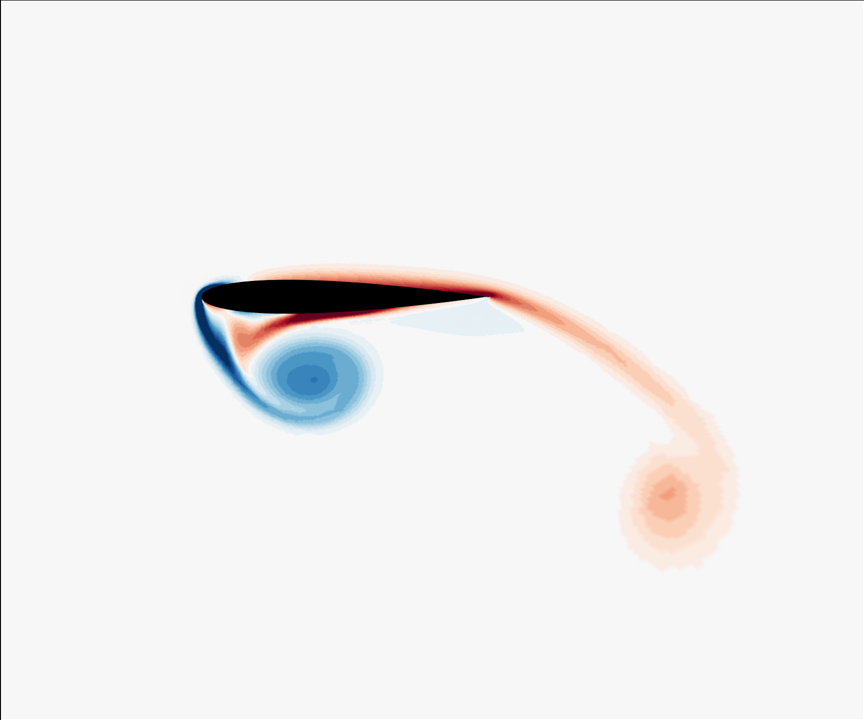}
 \end{subfigure} \hfill
 \begin{subfigure}{0.3\textwidth}
  \centering
  \includegraphics[width=\textwidth]{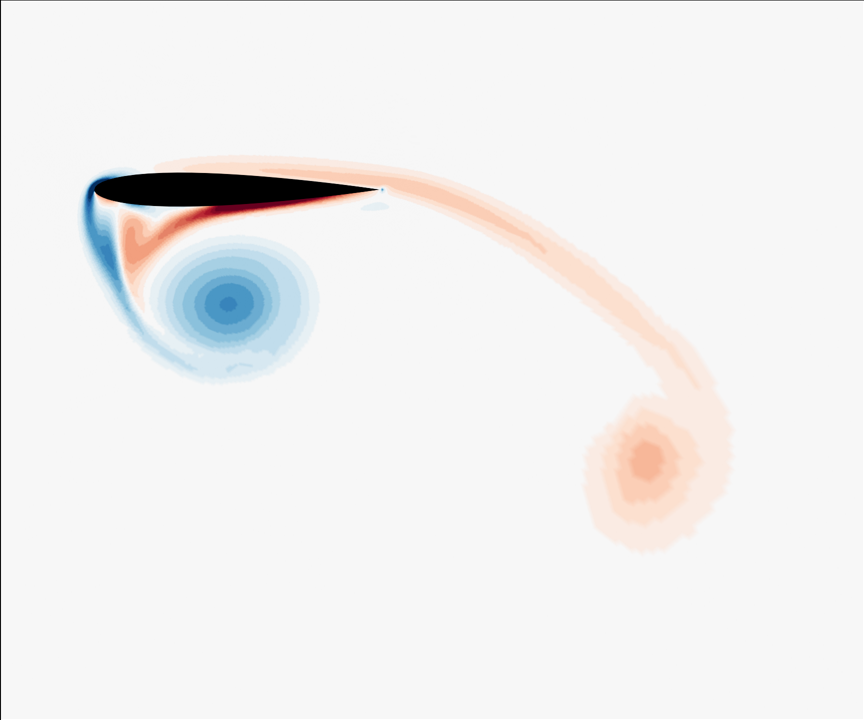}
 \end{subfigure} \hfill
 \begin{subfigure}{0.3\textwidth}
  \centering
  \includegraphics[width=\textwidth]{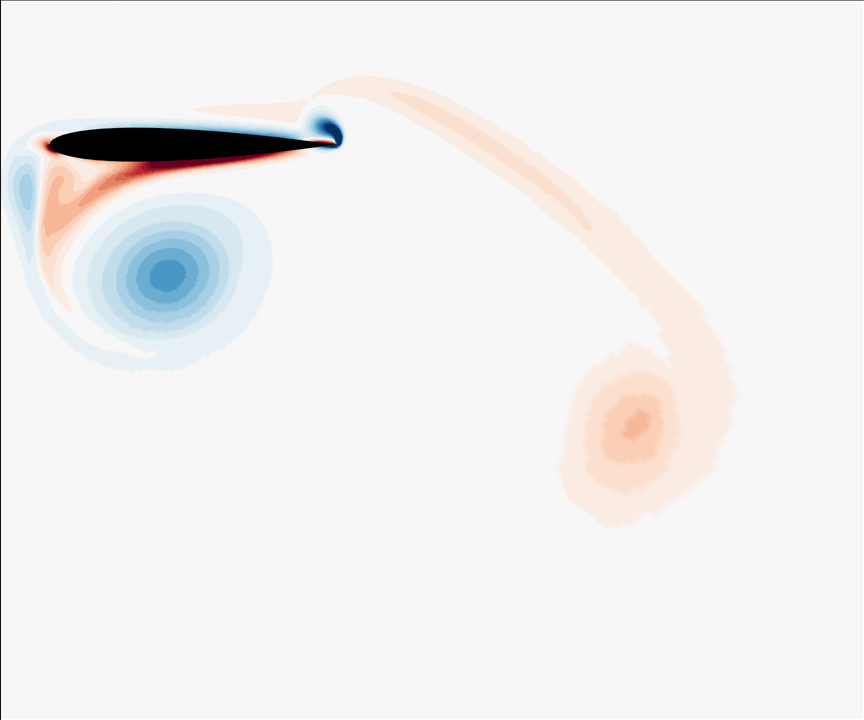}
 \end{subfigure}
 \caption{Flow vorticity around airfoil undergoing motion corresponding to
          initial guess for optimization, i.e., pure heaving
          (\ref{line:pitch-init}). Flow separation off leading edge implies a
          large amount of work required to complete mission.
          Snapshots taken at times $t = 0.0,~0.8,~1.6,~2.4,~3.2,~4.0$.}
 \label{fig:pitch-vort-init}
\end{figure}

\begin{figure}[!htbp]
 \centering
 \begin{subfigure}{0.3\textwidth}
  \centering
  \includegraphics[width=\textwidth]{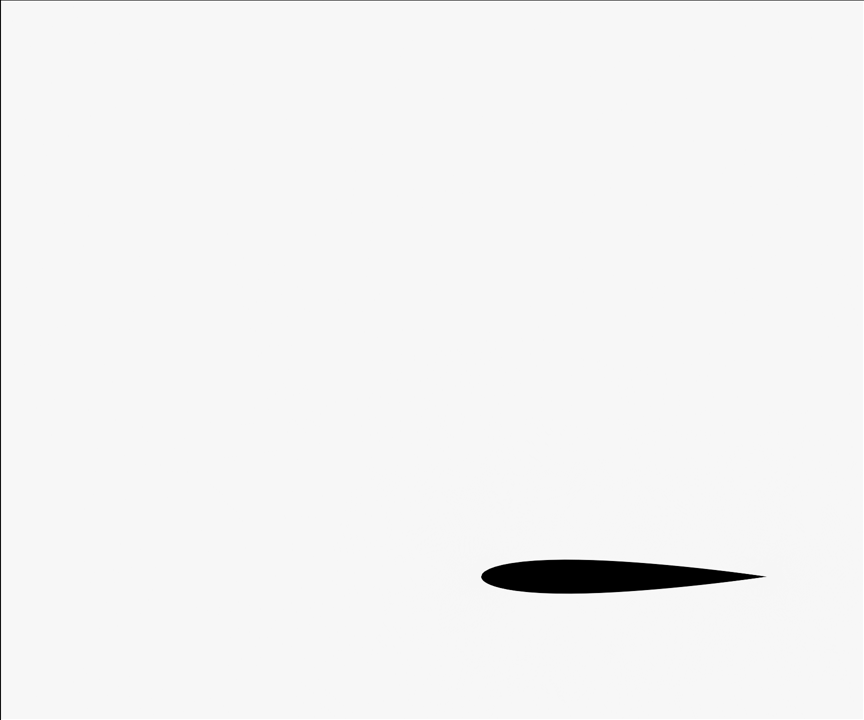}
 \end{subfigure} \hfill
 \begin{subfigure}{0.3\textwidth}
  \centering
  \includegraphics[width=\textwidth]{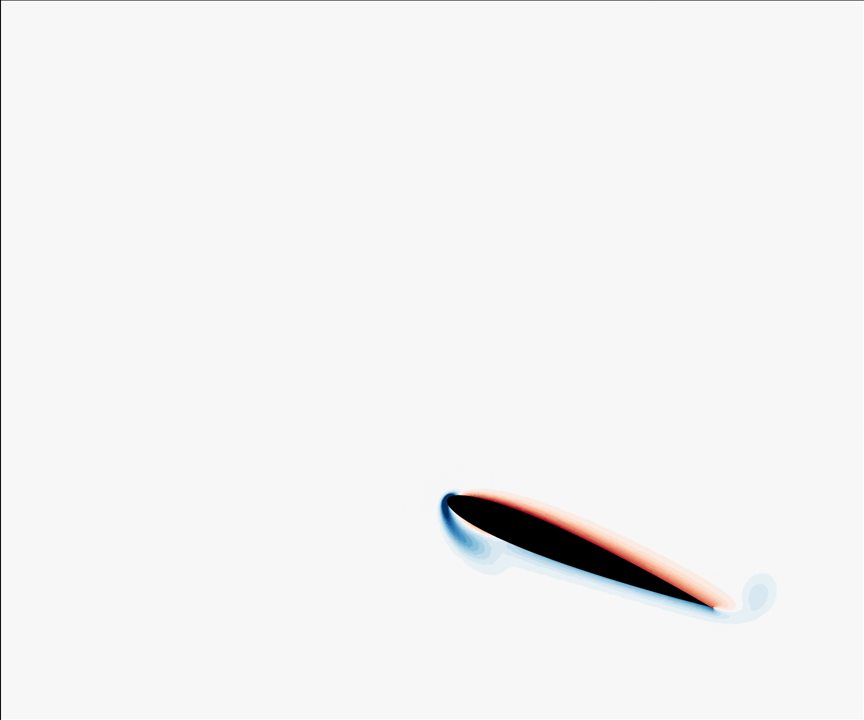}
 \end{subfigure} \hfill
 \begin{subfigure}{0.3\textwidth}
  \centering
  \includegraphics[width=\textwidth]{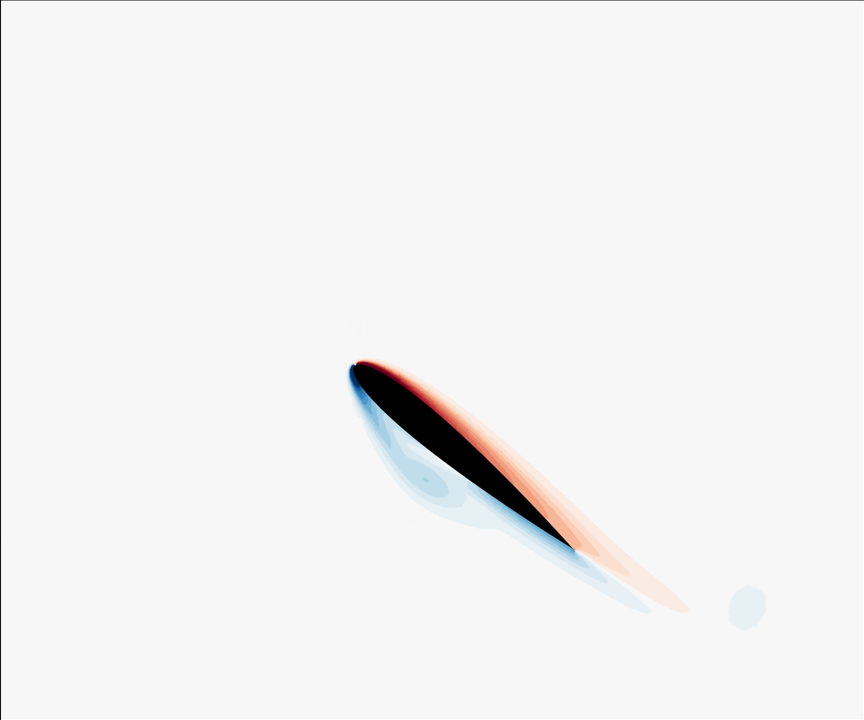}
 \end{subfigure} \\ \vspace{3mm}
 \begin{subfigure}{0.3\textwidth}
  \centering
  \includegraphics[width=\textwidth]{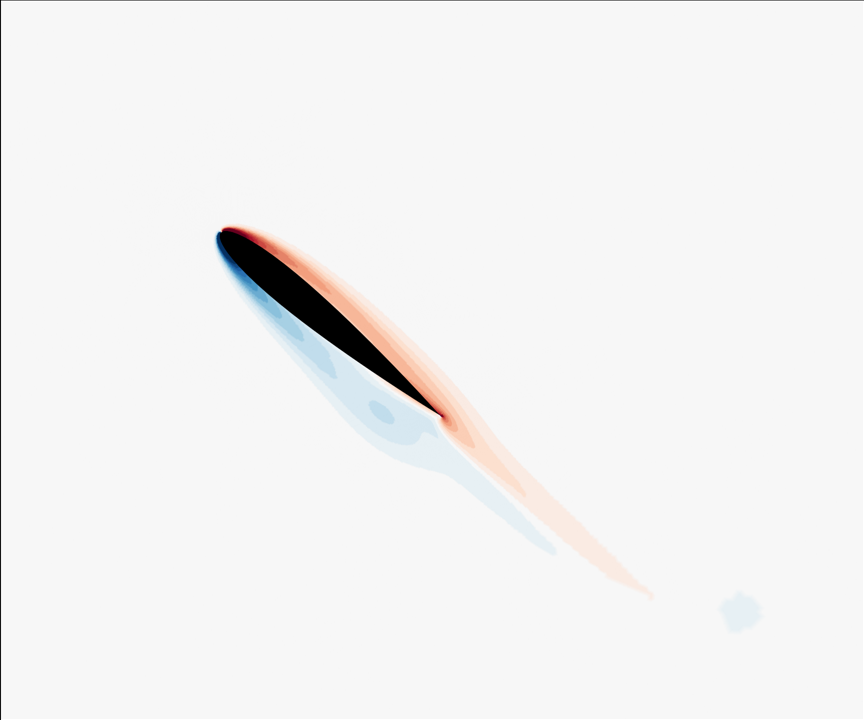}
 \end{subfigure} \hfill
 \begin{subfigure}{0.3\textwidth}
  \centering
  \includegraphics[width=\textwidth]{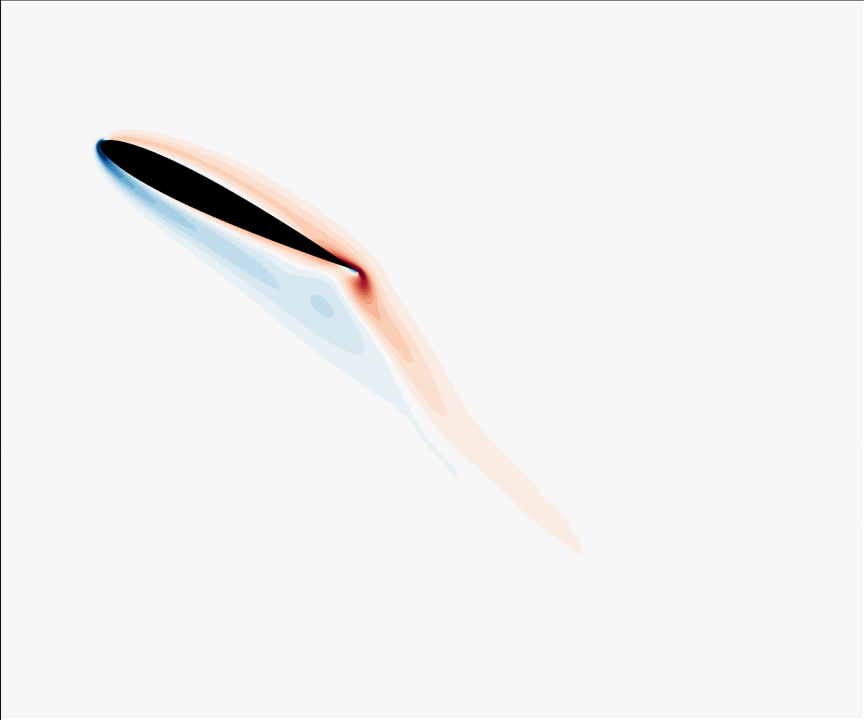}
 \end{subfigure} \hfill
 \begin{subfigure}{0.3\textwidth}
  \centering
  \includegraphics[width=\textwidth]{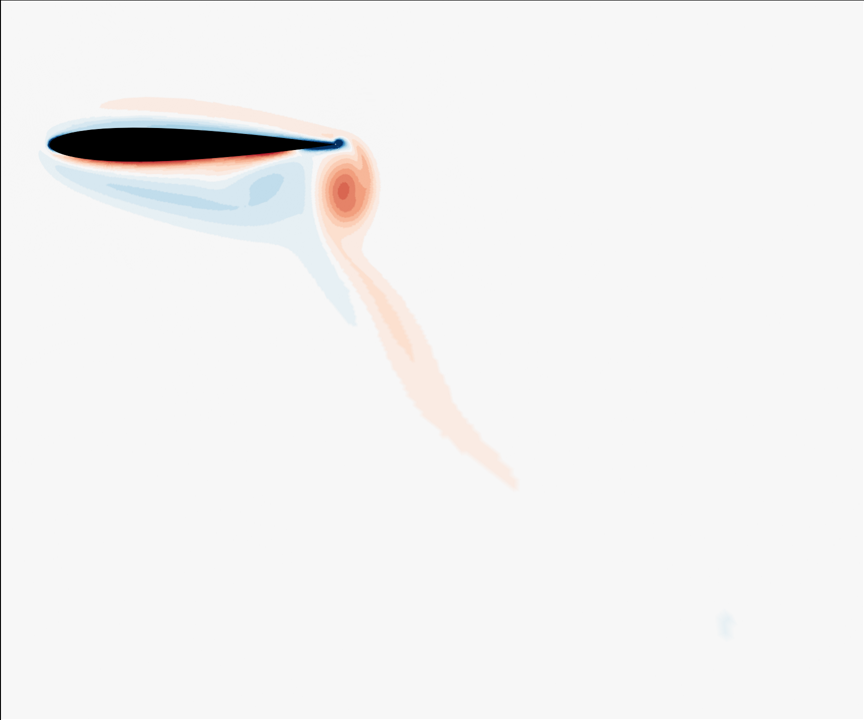}
 \end{subfigure}
 \caption{Flow vorticity around airfoil undergoing motion corresponding to
          optimal pitching motion for fixed translational motion, i.e., solution
          of (\ref{opt:pitch-disc}) under parametrization PI
          (\ref{line:pitch-opt1}). The pitching motion greatly reduces the
          degree of flow separation and vortex shedding compared to the initial
          guess, and requires less work to complete the mission.
          Snapshots taken at times $t = 0.0,~0.8,~1.6,~2.4,~3.2,~4.0$.}
 \label{fig:pitch-vort-opt1}
\end{figure}

\begin{figure}[!htbp]
 \centering
 \begin{subfigure}{0.3\textwidth}
  \centering
  \includegraphics[width=\textwidth]{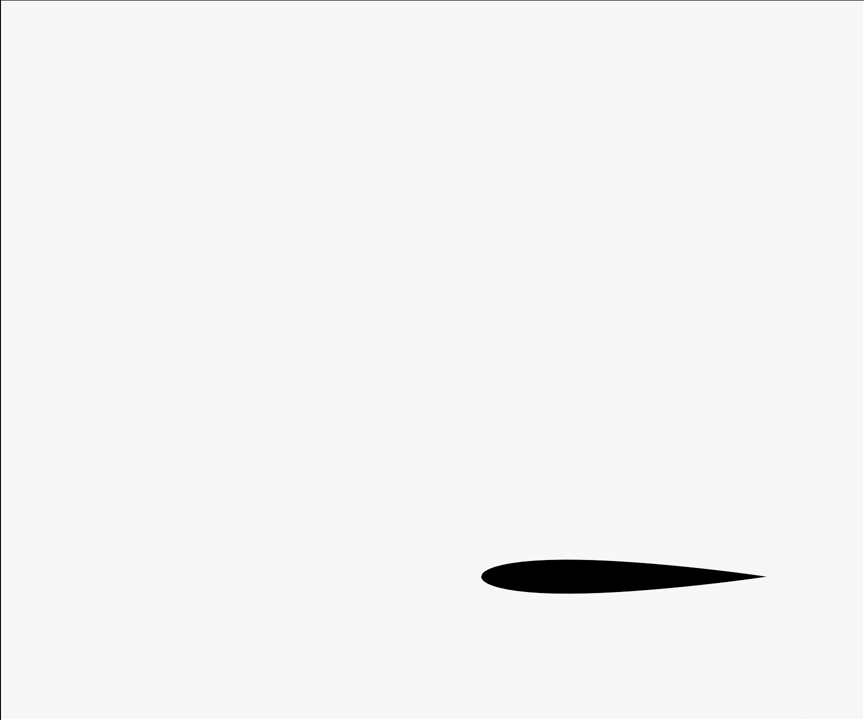}
 \end{subfigure} \hfill
 \begin{subfigure}{0.3\textwidth}
  \centering
  \includegraphics[width=\textwidth]{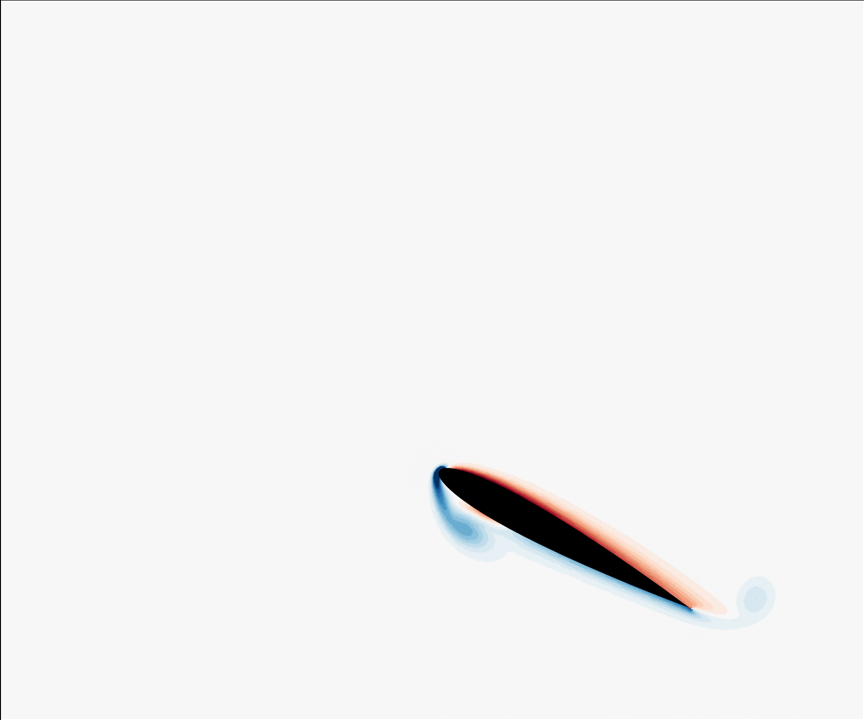}
 \end{subfigure} \hfill
 \begin{subfigure}{0.3\textwidth}
  \centering
  \includegraphics[width=\textwidth]{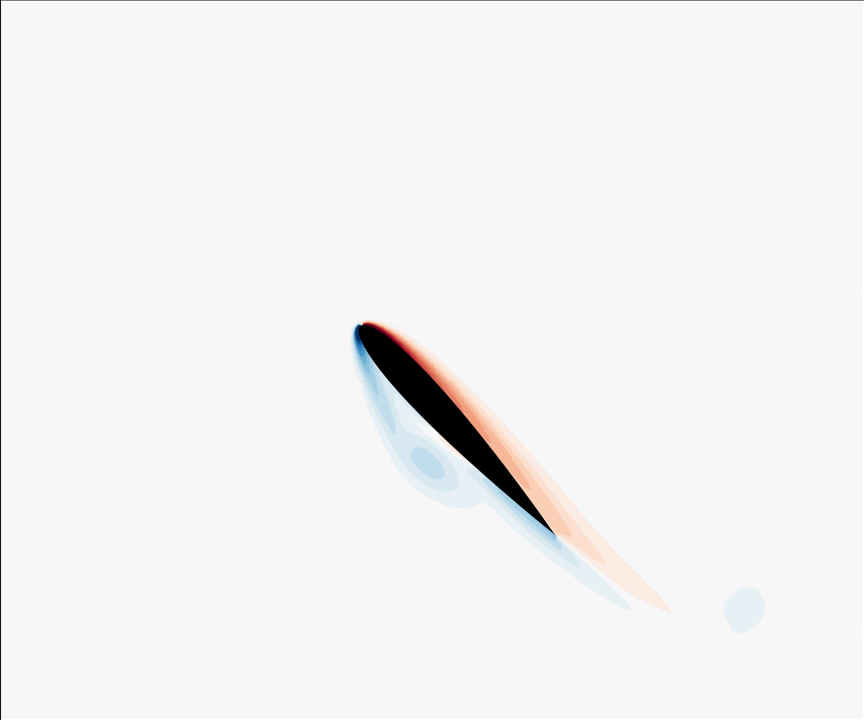}
 \end{subfigure} \\ \vspace{3mm}
 \begin{subfigure}{0.3\textwidth}
  \centering
  \includegraphics[width=\textwidth]{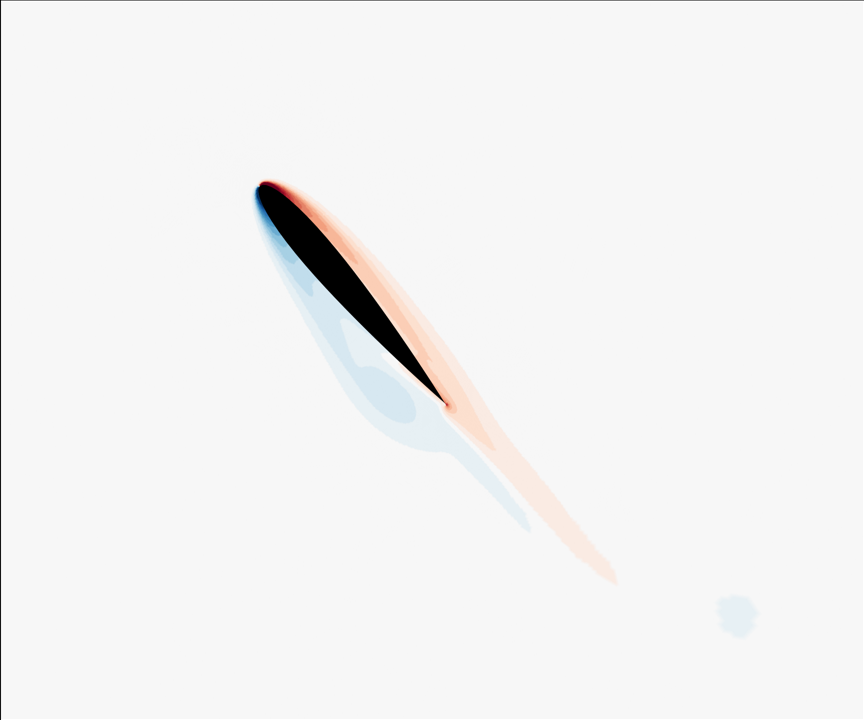}
 \end{subfigure} \hfill
 \begin{subfigure}{0.3\textwidth}
  \centering
  \includegraphics[width=\textwidth]{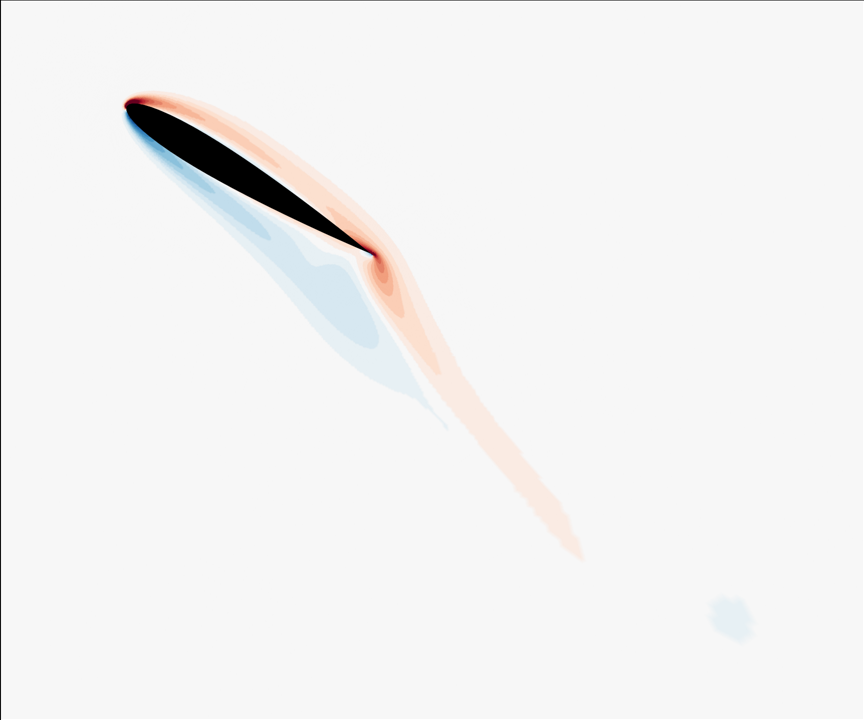}
 \end{subfigure} \hfill
 \begin{subfigure}{0.3\textwidth}
  \centering
  \includegraphics[width=\textwidth]{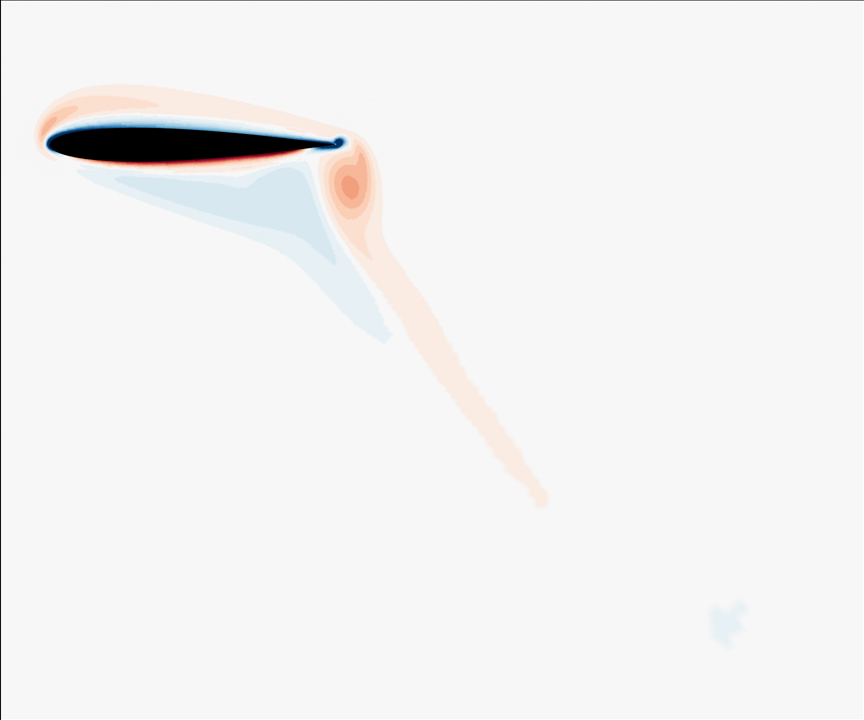}
 \end{subfigure}
 \caption{Flow vorticity around airfoil undergoing motion corresponding to
          optimal rigid body motion, i.e. solution of (\ref{opt:pitch-disc})
          under parametrization PII (\ref{line:pitch-opt2}). This rigid body
          motion further reduces the degree of flow separation and required work
          to complete the mission. This motion differs from the solution of PI
          as it has a larger pitch amplitude and slightly overshoots the final
          vertical position before settling to the required position.
          Snapshots taken at times $t = 0.0,~0.8,~1.6,~2.4,~3.2,~4.0$.}
 \label{fig:pitch-vort-opt2}
\end{figure}

The motion of the airfoil and vorticity of the surrounding flow are shown in
Figure~\ref{fig:pitch-vort-init} (nominal trajectory),
Figure~\ref{fig:pitch-vort-opt1} (optimal solution under parametrization PI),
and Figure~\ref{fig:pitch-vort-opt2} (optimal solution under parametrization
PII). The flow corresponding to the nominal configuration experiences flow
separation and vortex shedding, which results in the relatively large amount
of total energy to complete the mission. Fixing the translational motion and
optimizing over the rotation (PI) dramatically reduces the amount of shedding
and consequently reduces the amount of work required. Optimizing the entire
rigid body motion (PII) further reduces the shedding and required work.

\subsection{Energetically Optimal Shape and Flapping Motion of 2D Airfoil
            at Constant Impulse} \label{subsec:app-flap}
\begin{figure}[h]
  \centering
  \includegraphics[width=2.25in]{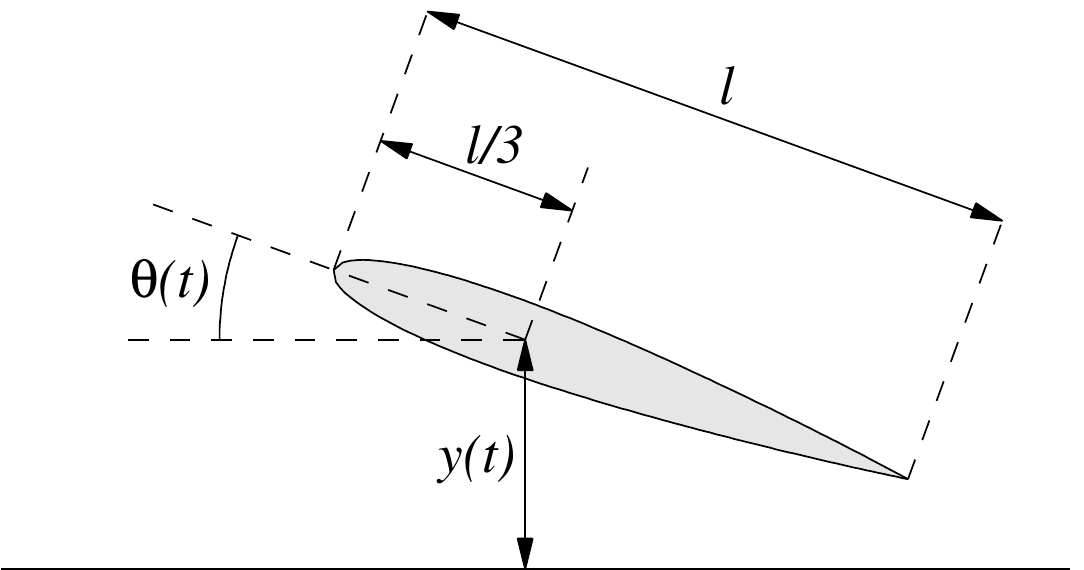} \hspace{2cm}
  \includegraphics[width=2.25in]{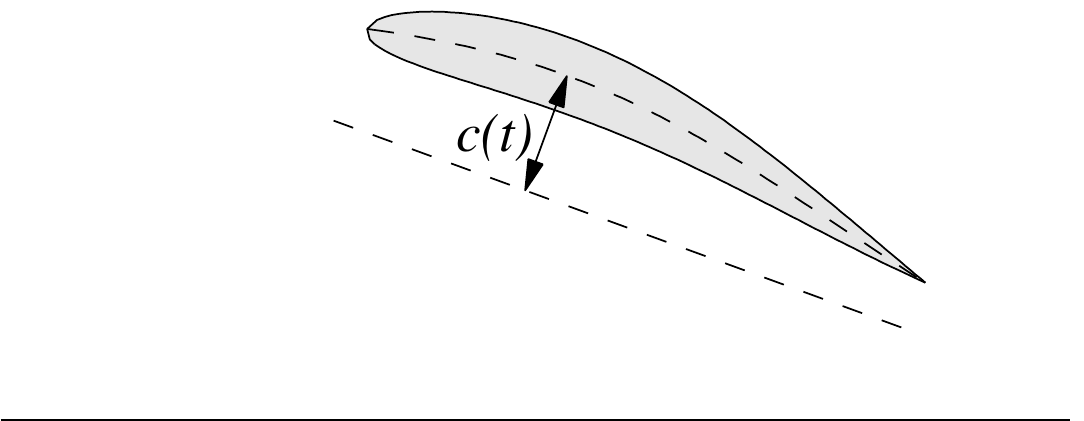}
  \caption{Airfoil}\label{fig:airfoil-traj-shape}
\end{figure}

In this section, the high-order, time-dependent PDE-constrained optimization
framework introduced in this document is applied to find the
energetically optimal flapping motion, under an impulse constraint, of a
2D NACA0012 airfoil (Figure~\ref{fig:airfoil-traj-shape}) with chord length
$l = 1$ and zero-thickness trailing edge. The governing equations are the
2D compressible, isentropic Navier-Stokes equations.

The goal is to determine the flapping motion -- $y(t)$ and $\theta(t)$ -- and
shape -- $c(t)$ -- of the airfoil that minimizes the total energy
such than a $x$-impulse of $q$ is achieved, i.e.,
\begin{equation} \label{opt:flap-cont}
  \begin{aligned}
    & \underset{\Ubm,~\mubold}{\text{minimize}}
    & & \Wcal(\Ubm, \mubold) \\
    & \text{subject to}
    & & \Jcal_x(\Ubm, \mubold) = q \\
    & & & \pder{\Ubm}{t} + \nabla \cdot \Fbm(\Ubm,~\nabla\Ubm) = 0 \quad
          \text{ in } v(\mubold, t).
  \end{aligned}
\end{equation}
The flapping frequency is fixed at $0.2$, which corresponds to a period
of $T = 5$. Proper initialization of the
flow is the initial condition that results in a \emph{time-periodic} flow
\citep{zahr2016periodic} to completely avoid non-physical transients and
simulate representative, \emph{in-flight} conditions; this work uses a crude
approximation that initializes the flow from the steady-state condition,
simulates $3$ periods of the flapping motion, and integrates the quantities
of interest over the last period only. The deformation of the domain is
determined from the value of $c(t)$ using the spatial blending
map of Section~\ref{subsec:implement-domain} with
\begin{equation}\label{eqn:morph-def}
 \varphibold(\Xbm, \mubold, t) =
 \begin{bmatrix}
   0 \\
   2c(t)e^{-[(\Xbm - \xbm_0) \cdot \ebm_1]^2}
 \end{bmatrix}.
\end{equation}
The deformation in (\ref{eqn:morph-def}) is volume preserving since
$\det(\Ibm + \nabla_\Xbm\varphibold) \equiv 1$, but it does not preserve
the surface area of the airfoil. A number of engineering situations exist
where time-dependent shape changes alter the surface area of an object
in a flow including control surface deployment on vehicles and extending
panels on wings \cite{supekar2007design}. Additionally, morphing in nature is
not necessarily surface area preserving \cite{ho2003unsteady}. Previous work
that has considered optimal time-morphed geometries, using a vortex-lattice
flow solver, also considered non-area-preserving deformations
\cite{stanford2010analytical, ghommem2012global}.
The time-dependent surface area is accounted for in the quantities of
interest (\ref{eqn:funcls}) through the spatial integral over the
time-dependent boundary, $\Gammabold$, i.e., $\int_\Gammabold \cdot\,dS$.

The trajectory of the airfoil -- $y(t)$, and $\theta(t)$ -- and its shape
-- $c(t)$ -- are discretized via cubic splines with
$m_y+1$, $m_\theta+1$, and $m_c+1$ knots, respectively, with
boundary conditions that enforce
\begin{equation} \label{eqn:mirror}
  y(t) = -y(t+T/2) \qquad \qquad
  \theta(t) = -\theta(t+T/2)  \qquad \qquad
  c(t) = -c(t+T/2).
\end{equation}
These boundary conditions\footnote{Periodic and mirrored cubic splines of this
form with $m+1$ knots only have $m$ degrees of freedom since the boundary
condition prescribes the value of the $m+1$ knot from the values of the
others $m$.} for $y(t)$, $\theta(t)$, and $c(t)$ correspond to a
mirroring of the trajectory at $t = T/2$ and implicitly enforces periodicity
with period $T$.
The knots are uniformly spaced between $0$ and $T$ in the $t$-dimension and
the knot values are optimization parameters. Since the unsteady simulation is
initialized from the steady-state flow, non-zero velocities of the airfoil at
$t = 0$ will result in non-physical transients. These transients are avoided
by blending the periodic cubic spline smoothly to the zero function at the
beginning of the time interval \citep{van2013adjoint}. Let $s_y(t; \mubold)$,
$s_\theta(t; \mubold)$, and $s_c(t; \mubold)$ denote the
periodic cubic spline approximations. Then, the flapping and shape
trajectories are defined as
\begin{equation}
   y(t) = b(t) s_y(t; \mubold) \qquad \qquad
   \theta(t) = b(t) s_\theta(t; \mubold) \qquad \qquad
   c(t) = b(t) s_c(t; \mubold),
\end{equation}
where $\displaystyle{b(t) = 1.0 - e^{-t^2}}$.
\begin{table}[!htbp]
 \centering
 \input{dat/flap-setup.tab}
 \caption{Summary of parametrizations considered in
          Section~\ref{subsec:app-flap}. The number of periodic cubic spline
          knots used to discretize $y(t)$, $\theta(t)$, and $\c(t)$
          are $m_y+1$, $m_\theta+1$, and $m_c+1$, respectively.
          FI freezes the airfoil shape and considers only rigid body motions
          ($m_y = m_\theta \neq 0, m_c = 0$). FII parametrizes both
          shape and kinematic motion ($m_y = m_\theta = m_c \neq 0$).}
\label{tab:flap-setup}
\end{table}
Table~\ref{tab:flap-setup} summarizes two parametrizations considered in this
section:
\begin{inparaenum}
 \item[(FI)]  rigid body motion parametrized via cubic splines and shape
              fixed at nominal value and
 \item[(FII)] rigid body motion and shape of airfoil parametrized via cubic
              splines.
\end{inparaenum}
With this parametrization of the airfoil kinematics and shape, spatial and
temporal discretization with the high-order scheme of Section~\ref{sec:govern}
leads to the  fully discrete version of the optimization problem in
(\ref{opt:flap-cont})
\begin{equation} \label{opt:flap-disc}
  \begin{aligned}
    & \underset{\substack{\ubm^{(0)},~\dots,~\ubm^{(N_t)} \in \Rbb^{N_\ubm},\\
                        \kbm_1^{(1)},~\dots,~\kbm_s^{(N_t)} \in \Rbb^{N_\ubm},\\
                        \mubold \in \Rbb^{N_\mubold}}}
                        {\text{minimize}}
    & & W(\ubm^{(0)},~\dots,~\ubm^{(N_t)},~\kbm_1^{(1)},~\dots,~\kbm_s^{(N_t)},~
          \mubold) \\
    & \text{subject to}
    & & J_x(\ubm^{(0)},~\dots,~\ubm^{(N_t)},~\kbm_1^{(1)},~\dots,
            ~\kbm_s^{(N_t)},~\mubold) = q \\
    & & &   \ubm^{(0)} = \ubm_0 \\
    & & & \ubm^{(n)} = \ubm^{(n-1)} + \sum_{i = 1}^s b_i\kbm^{(n)}_i \\
    & & & \Mbb\kbm^{(n)}_i = \Delta t_n\rbm\left(\ubm_i^{(n)},~\mubold,~
                                                 t_{n-1} + c_i\Delta t_n\right).
  \end{aligned}
\end{equation}

\begin{figure}[!htbp]
  \centering
  \hspace*{\fill}%
  \begin{subfigure}{0.49\textwidth}
    \input{tikz/flap-y-hist.tikz}
  \end{subfigure} \hfill
  \begin{subfigure}{0.49\textwidth}
   \input{tikz/flap-th-hist.tikz}
  \end{subfigure} \hspace*{\fill} \\
  \begin{subfigure}{0.49\textwidth}
   \input{tikz/flap-camber-hist.tikz}
  \end{subfigure}
  \caption{Trajectories of $y(t)$, $\theta(t)$, and $c(t)$ at initial
           guess (\ref{line:flap-init}), solution of (\ref{opt:flap-disc})
           under parametrization FI
           ($q = 0.0$: \ref{line:flap-opt1-0},
            $q = 1.0$: \ref{line:flap-opt1-1},
            $q = 2.5$: \ref{line:flap-opt1-2.5}), and solution of
           (\ref{opt:flap-disc}) under parametrization FII
           ($q = 0.0$: \ref{line:flap-opt2-0},
            $q = 1.0$: \ref{line:flap-opt2-1},
            $q = 2.5$: \ref{line:flap-opt2-2.5}) from
            Table~\ref{tab:flap-setup}.}
  \label{fig:flap-traj-kine}
\end{figure}
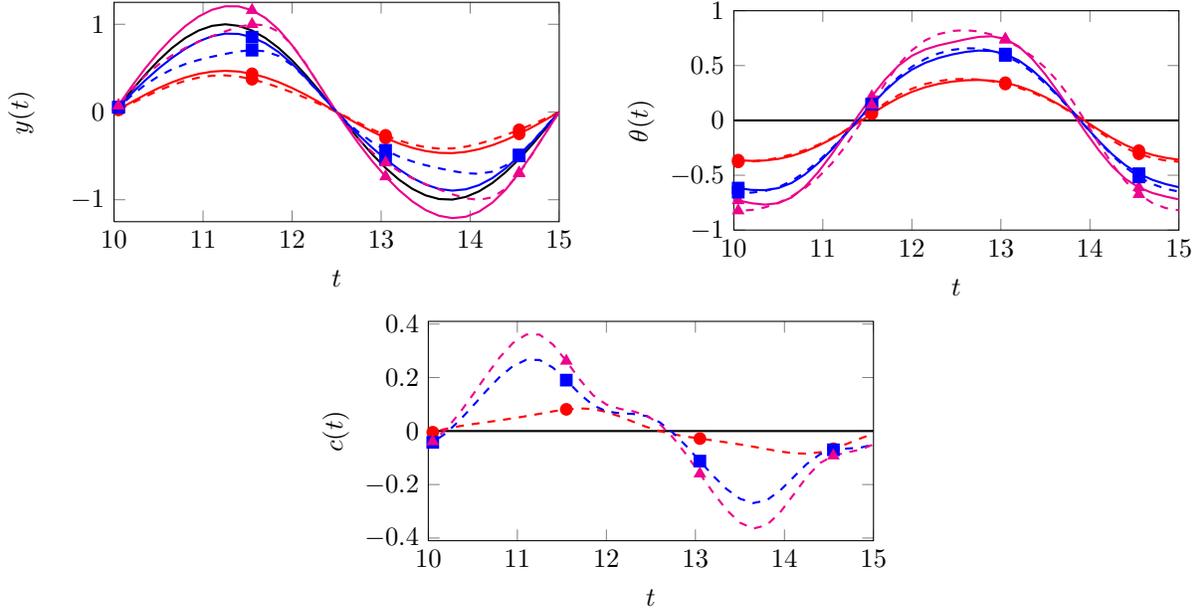

Given the gradient verification from the previous section, attention is turned
directly to the optimization problem in (\ref{opt:flap-disc}) for various
values of the impulse constraint, $q$. The optimization solver used in this
section is SNOPT \citep{gill2002snopt}, a nonlinearly constrained SQP method.
Figure~\ref{fig:flap-traj-kine} contains the initial guess for the optimization
problem in (\ref{opt:pitch-disc}) as well as its solution under both
parametrization, FI and FII. The initial guess for the optimization problem is
a pure heaving motion at a fixed shape, i.e., $c(t) = \theta(t) = 0$. The
solution under parametrization PI freezes the shape at its nominal
configuration (NACA0012) and modifies the rigid body motion. Pitch is
introduced for all values of the impulse constraint and the amplitude of the
heaving motion is decreased for $q = 0.0, 1.0$ and increased for $q = 2.5$.  The
solution under parametrization PII reduces the heaving amplitude and slightly
increases the pitch amplitude as compared to PI. It also introduces
non-trivial camber.

\begin{figure}[!htbp]
  \centering
  \input{tikz/flap-hist.tikz}
  \caption{Time history of total power, $\Pcal^h(\ubm, \mubold, t)$, and
           $x$-directed force, $\Fcal_x^h(\ubm, \mubold, t)$, imparted onto foil
           by fluid at initial guess (\ref{line:flap-init}), solution of
           (\ref{opt:flap-disc}) under parametrization FI
           ($q = 0.0$: \ref{line:flap-opt1-0},
            $q = 1.0$: \ref{line:flap-opt1-1},
            $q = 2.5$: \ref{line:flap-opt1-2.5}), and solution of
           (\ref{opt:flap-disc}) under parametrization FII
           ($q = 0.0$: \ref{line:flap-opt2-0},
            $q = 1.0$: \ref{line:flap-opt2-1},
            $q = 2.5$: \ref{line:flap-opt2-2.5}) from
            Table~\ref{tab:flap-setup}.}
  \label{fig:flap-traj-qoi}
\end{figure}
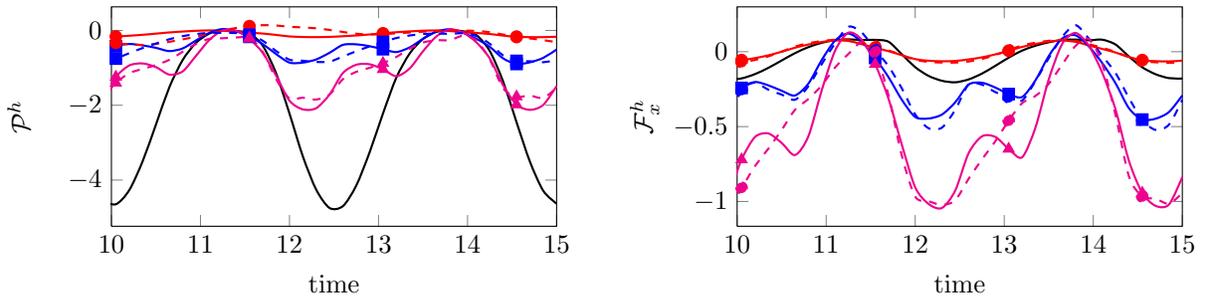

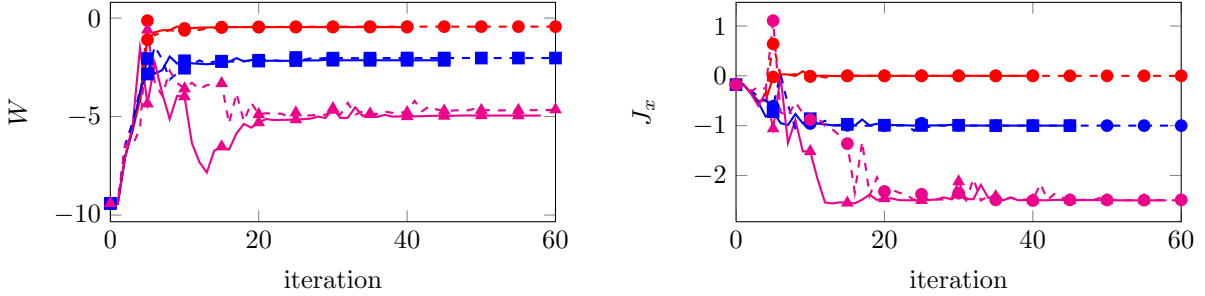
\begin{figure}[!htbp]
  \centering
  \input{tikz/flap-conv.tikz}
  \caption{Convergence of quantities of interest, $W$ and $J_x$, with
           optimization iteration for parametrization FI
           ($q = 0.0$: \ref{line:flap-opt1-0},
            $q = 1.0$: \ref{line:flap-opt1-1},
            $q = 2.5$: \ref{line:flap-opt1-2.5}) and FII
           ($q = 0.0$: \ref{line:flap-opt2-0},
            $q = 1.0$: \ref{line:flap-opt2-1},
            $q = 2.5$: \ref{line:flap-opt2-2.5}) from
           Table~\ref{tab:flap-setup}. Each optimization iteration requires
           a primal flow computation -- to evaluate quantities of interest --
           and its corresponding adjoint -- to evaluate the gradient of
           quantities of interest.}
  \label{fig:flap-conv}
\end{figure}

The instantaneous quantities of interest -- $W$ and $J_x$ in this case -- for
the nominal motion and shape and solution of (\ref{opt:flap-disc}) under
parametrizations PI and PII are included in Figure~\ref{fig:flap-traj-qoi}.
It is clear that the optimal solution under both parametrizations result in a
time history of the total power that is uniformly closer to $0$ than that at
the nominal trajectory, which is expected since $W$ is the objective function.
It is also clear that larger values of the impulse constraint require more
power to complete the flapping motion. While it may not be clear from
Figure~\ref{fig:flap-traj-qoi}, the integration of $\Fcal_x^h$ leads to an
impulse that exactly conforms to the specified value of $q$. This can be
seen more clearly in Figure~\ref{fig:flap-conv}. These observations can also
be verified in Figure~\ref{fig:flap-conv} and Table~\ref{tab:flap-summ} that
summarizes the optimal values of the integrated quantities of interest.

Figure~\ref{fig:flap-conv} shows the convergence of the integrated quantities
of interest with iterations in the optimization solver. The aforementioned
observations can be verified by inspection of the final iteration: all
impulse constraints are satisfied, larger values of $q$ require more work to
achieve, and morphing the shape of the airfoil allows for a slight reduction
in the required work. After $20$ iterations, the impulse constraint is
satisfied for $q = 0.0, 1.0$ and reduction of the work has essentially ceased,
implying the optimization could have been terminated at that point. The case
with $q = 2.5$ requires an additional 15 - 20 iterations to settle to a
converged solution.

\begin{figure}[!htbp]
 \centering
 \begin{subfigure}{0.3\textwidth}
  \centering
  \includegraphics[width=\textwidth]{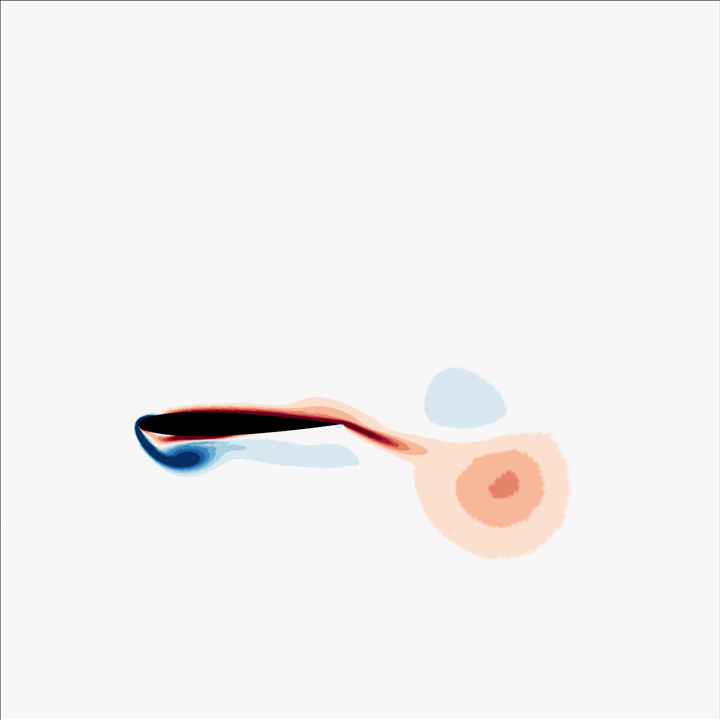}
 \end{subfigure} \hfill
 \begin{subfigure}{0.3\textwidth}
  \centering
  \includegraphics[width=\textwidth]{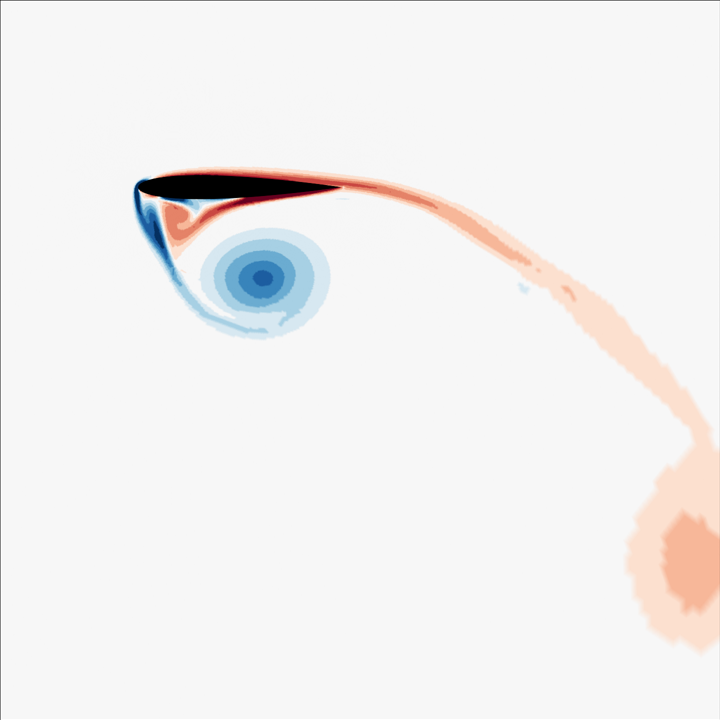}
 \end{subfigure} \hfill
 \begin{subfigure}{0.3\textwidth}
  \centering
  \includegraphics[width=\textwidth]{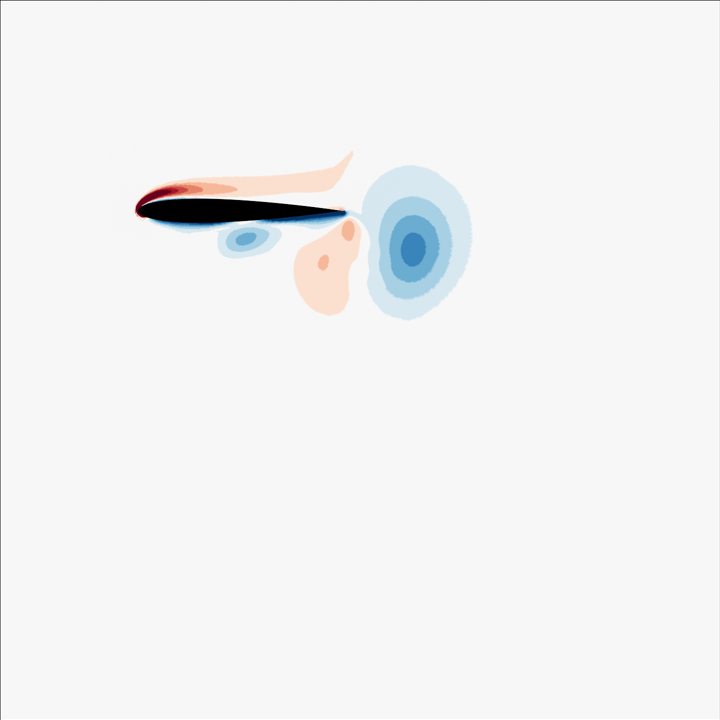}
 \end{subfigure} \\ \vspace{3mm}
 \begin{subfigure}{0.3\textwidth}
  \centering
  \includegraphics[width=\textwidth]{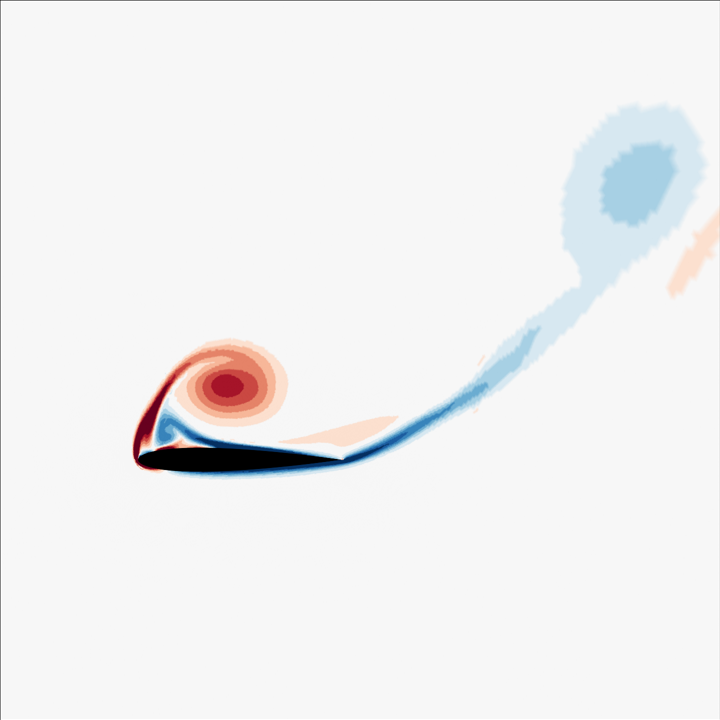}
 \end{subfigure} \hfill
 \begin{subfigure}{0.3\textwidth}
  \centering
  \includegraphics[width=\textwidth]{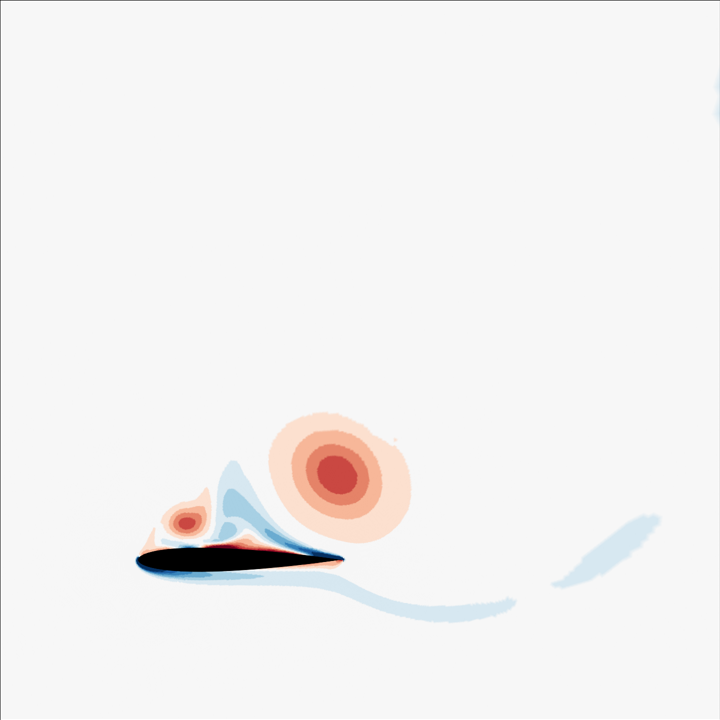}
 \end{subfigure} \hfill
 \begin{subfigure}{0.3\textwidth}
  \centering
  \includegraphics[width=\textwidth]{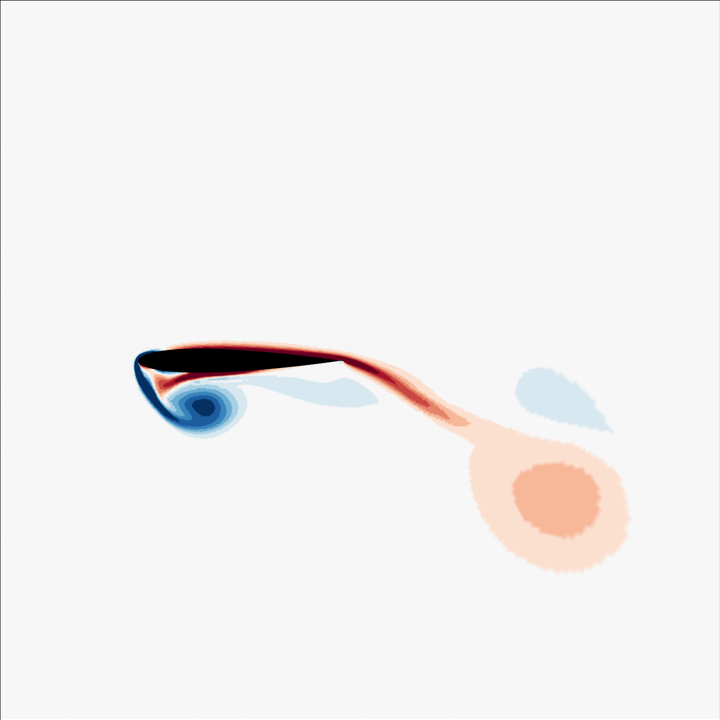}
 \end{subfigure}
 \caption{Flow vorticity around flapping airfoil undergoing motion corresponding
          to initial guess for optimization problem (\ref{opt:flap-disc}), i.e.,
          pure heaving (\ref{line:flap-init}). Flow separation off leading edge
          implies a large amount of work required for flapping motion.
          Snapshots taken at times $t = 9.75,~10.8,~11.85,~12.9,~13.95,~15.0$.}
 \label{fig:flap-vort-init}
\end{figure}

\begin{figure}[!htbp]
 \centering
 \begin{subfigure}{0.3\textwidth}
  \centering
  \includegraphics[width=\textwidth]{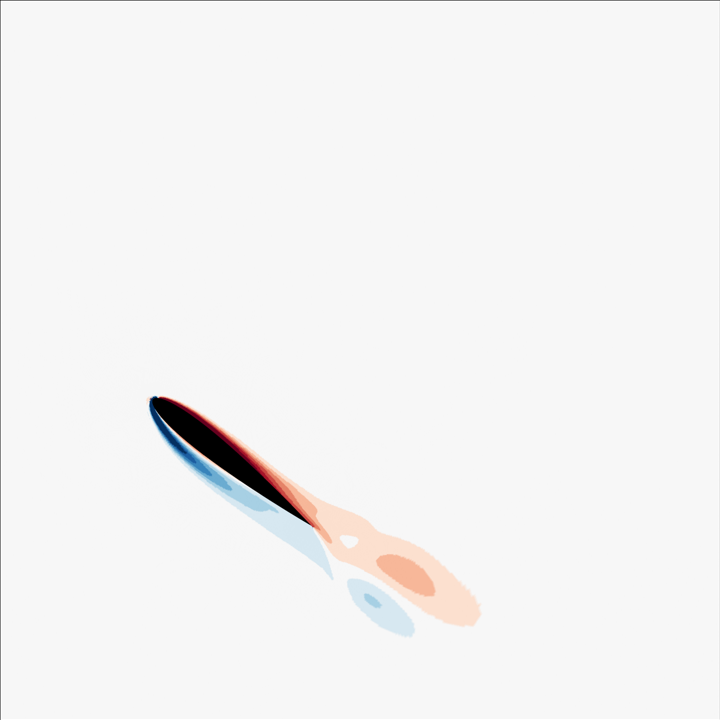}
 \end{subfigure} \hfill
 \begin{subfigure}{0.3\textwidth}
  \centering
  \includegraphics[width=\textwidth]{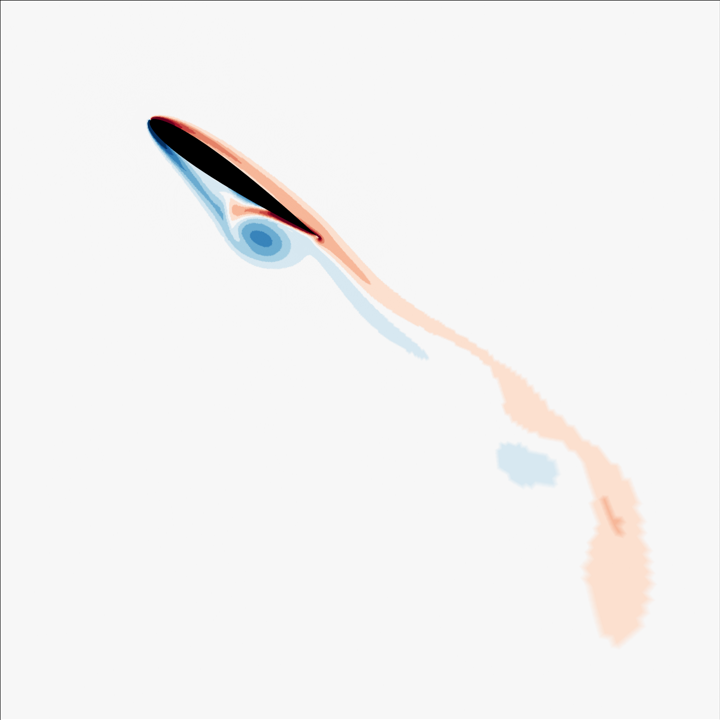}
 \end{subfigure} \hfill
 \begin{subfigure}{0.3\textwidth}
  \centering
  \includegraphics[width=\textwidth]{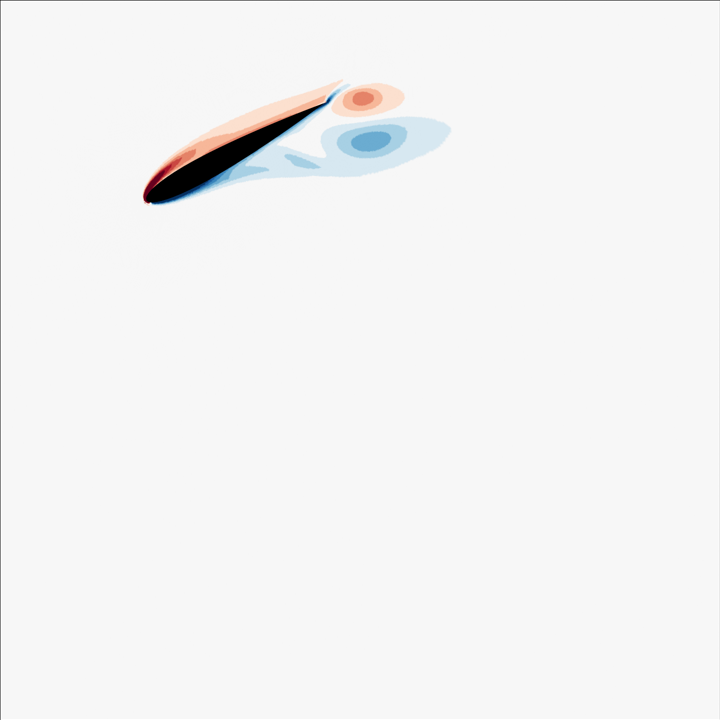}
 \end{subfigure} \\ \vspace{3mm}
 \begin{subfigure}{0.3\textwidth}
  \centering
  \includegraphics[width=\textwidth]{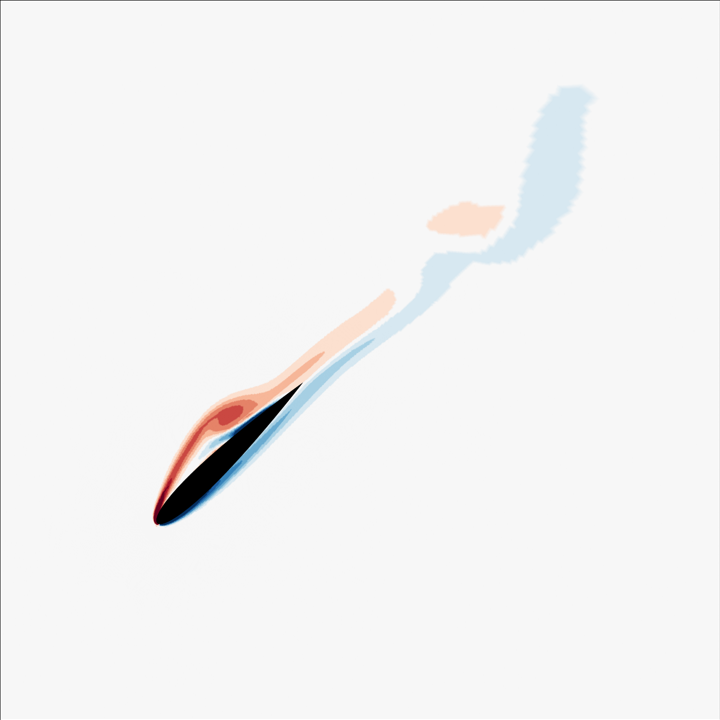}
 \end{subfigure} \hfill
 \begin{subfigure}{0.3\textwidth}
  \centering
  \includegraphics[width=\textwidth]{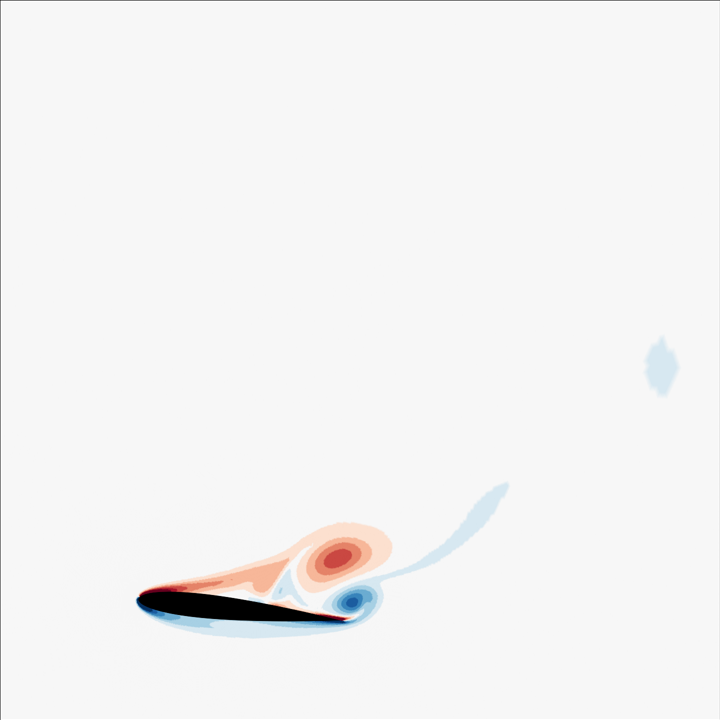}
 \end{subfigure} \hfill
 \begin{subfigure}{0.3\textwidth}
  \centering
  \includegraphics[width=\textwidth]{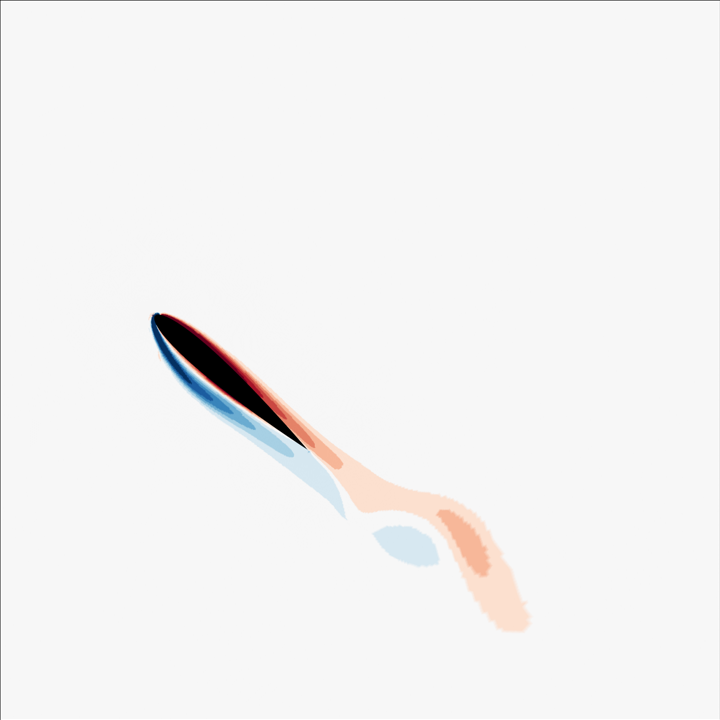}
 \end{subfigure}
 \caption{Flow vorticity around flapping airfoil undergoing optimal rigid
          body motion corresponding to the solution of (\ref{opt:flap-disc})
          under parametrization FI. The $x$-directed impulse is $J_x = 2.5$.
          The pitching motion greatly reduces the degree of flow separation and
          vortex shedding compared to the initial guess, and requires less work
          to complete the flapping motion and generate desired impulse.
          Snapshots taken at times $t = 9.75,~10.8,~11.85,~12.9,~13.95,~15.0$.}
 \label{fig:flap-vort-opt1}
\end{figure}

\begin{figure}[!htbp]
 \centering
 \begin{subfigure}{0.3\textwidth}
  \centering
  \includegraphics[width=\textwidth]{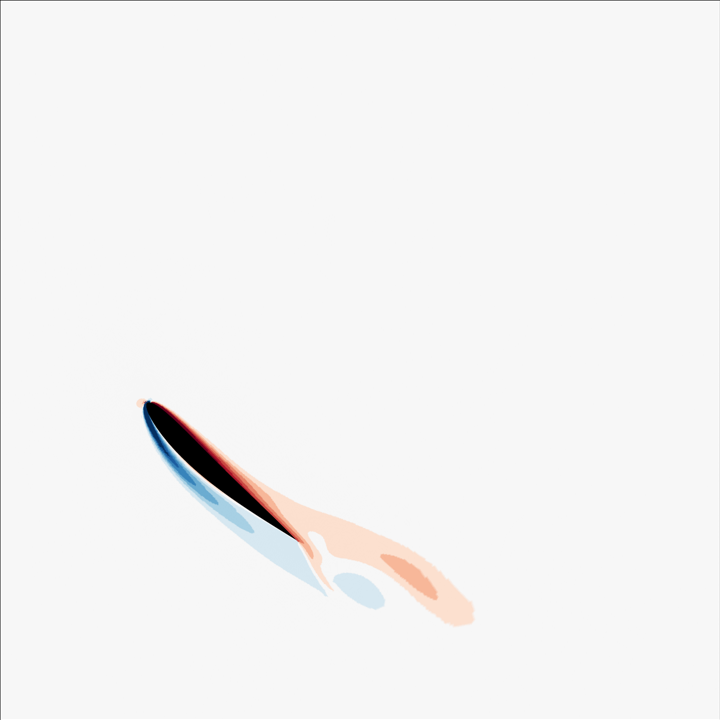}
 \end{subfigure} \hfill
 \begin{subfigure}{0.3\textwidth}
  \centering
  \includegraphics[width=\textwidth]{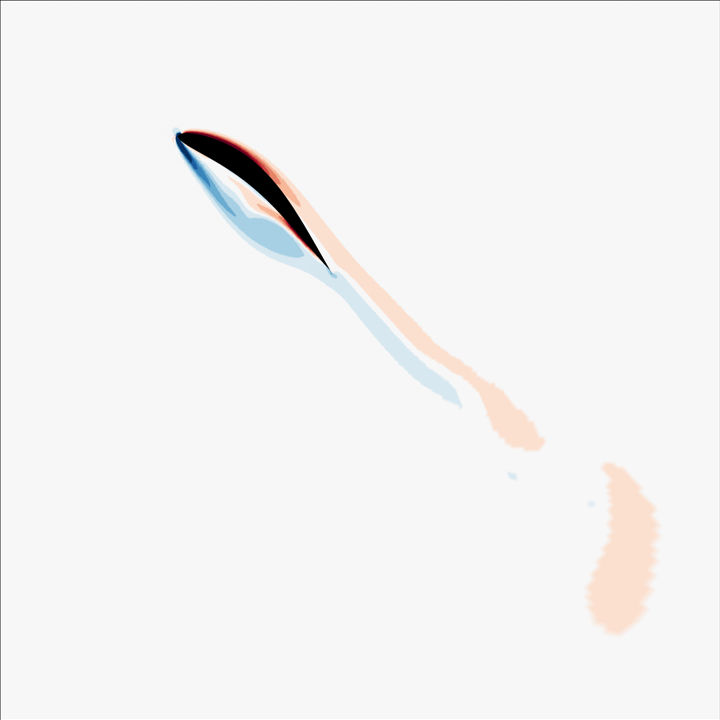}
 \end{subfigure} \hfill
 \begin{subfigure}{0.3\textwidth}
  \centering
  \includegraphics[width=\textwidth]{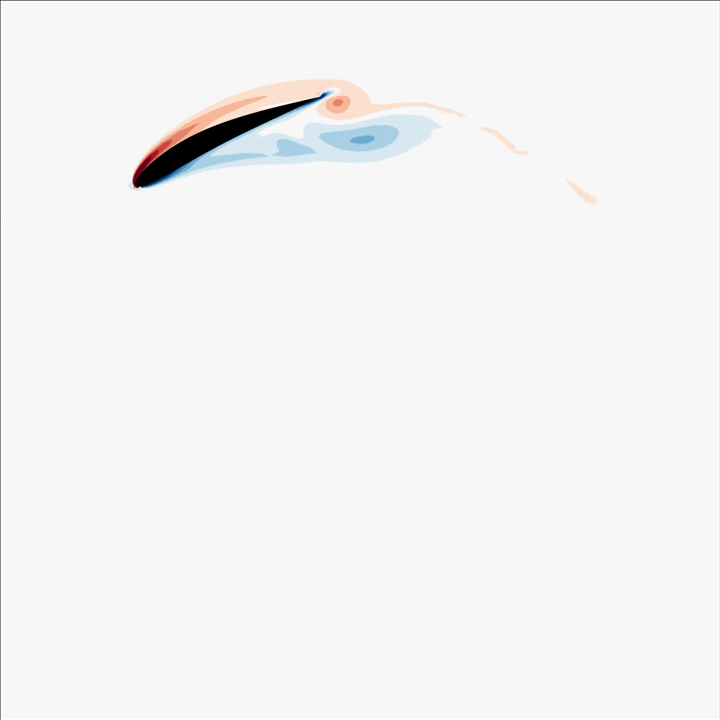}
 \end{subfigure} \\ \vspace{3mm}
 \begin{subfigure}{0.3\textwidth}
  \centering
  \includegraphics[width=\textwidth]{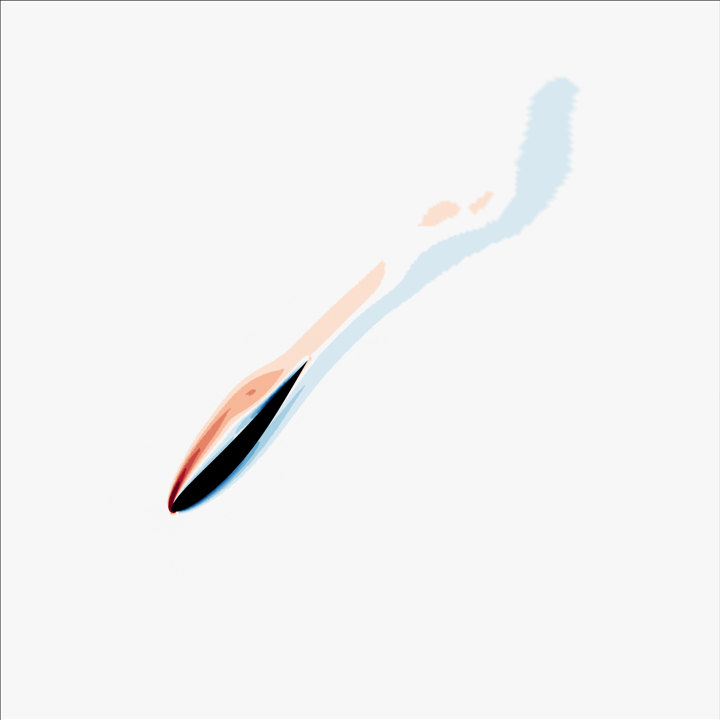}
 \end{subfigure} \hfill
 \begin{subfigure}{0.3\textwidth}
  \centering
  \includegraphics[width=\textwidth]{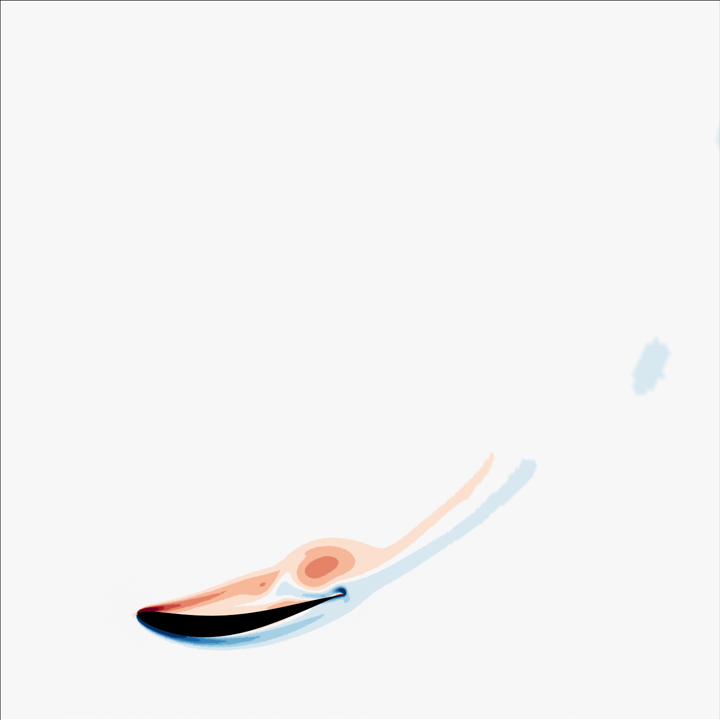}
 \end{subfigure} \hfill
 \begin{subfigure}{0.3\textwidth}
  \centering
  \includegraphics[width=\textwidth]{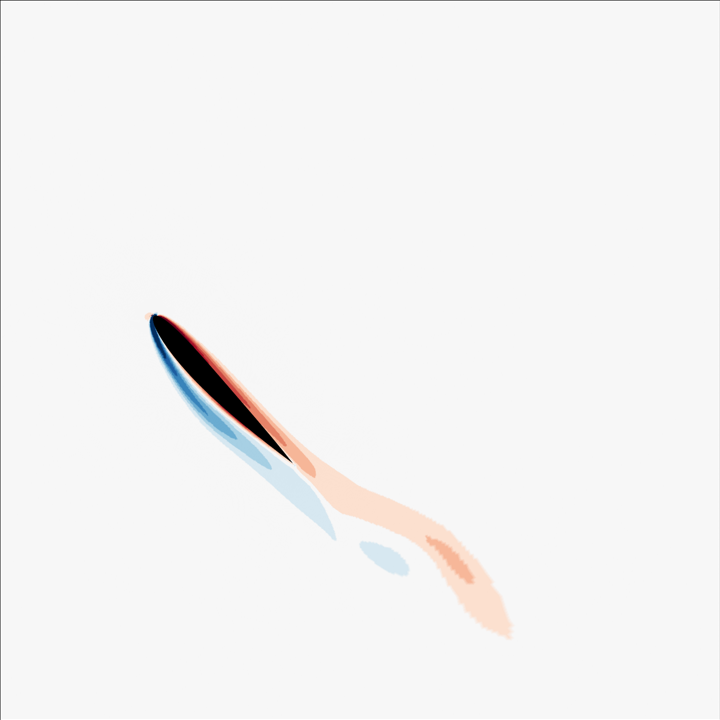}
 \end{subfigure}
 \caption{Flow vorticity around flapping airfoil undergoing optimal deformation
          and kinematic motion, corresponding to the solution of
          (\ref{opt:flap-disc}) under parametrization FII. The $x$-directed
          impulse is $J_x = 2.5$. The morphing further reduces the flow
          separation and work required to complete the flapping motion and
          generate desired impulse.
          Snapshots taken at times $t = 9.75,~10.8,~11.85,~12.9,~13.95,~15.0$.}
 \label{fig:flap-vort-opt2}
\end{figure}

The shape and motion of the airfoil and vorticity of the surrounding flow are
shown in Figure~\ref{fig:flap-vort-init} (nominal),
Figure~\ref{fig:flap-vort-opt1} (optimal solution under parametrization FI
for $q = 2.5$), and Figure~\ref{fig:flap-vort-opt2} (optimal solution under
parametrization FII for $q = 2.5$). The flow corresponding to the nominal
configuration experiences flow separation and vortex shedding, which results
in the relatively large amount of total energy to complete the flapping motion
and \emph{does not satisfy the impulse constraint}. Fixing the shape
and optimizing over the heaving and pitching motion (FI) dramatically reduces
the amount of shedding and consequently reduces the amount of work required.
Optimizing the shape in addition to the pitching and heaving motion (FII)
further reduces the shedding and required work. The solution of FI and FII
both satisfy the impulse constraint to greater than $8$ digits of accuracy.

\begin{table}[!htbp]
 \centering
 \input{dat/flap-summ.tab}
 \caption{Table summarizing integrated quantities of interest at optimal
          solution of each optimization problem for each impulse level. In all
          cases, the desired value of $J_x$ is achieved to greater than $4$
          digits of accuracy. The optimal solution for larger values of the
          impulse constraint require more total work to complete flapping
          motion, i.e., work monotonically increases in magnitude as value of
          impulse constraint increases. Smaller values of total work are
          achievable if airfoil is allowed to morph its shape in addition its
          rigid body motion. The other integrated quantities are included for
          completeness, but do not exhibit trends since they were not in the
          optimization problem.}
 \label{tab:flap-summ}
\end{table}

To conclude this section, a brief comparison of the optimal flapping motions
found in this work are compared to those found in the literature. From
Figure~\ref{fig:flap-traj-kine}, the pitch of the foil leads its plunge by
approximately $90^\circ$ in all optimal flapping motions, a result that was
found in several works that range from experimental to computational
\cite{tuncer2005optimization, ramamurti2001simulation, platzer2008flapping,
      oyama2009aerodynamic}.
The improved efficiency is largely due to a dramatic reduction in leading
edge vortex shedding characteristic of pure heaving motions
(Figure~\ref{fig:flap-vort-init})
\cite{tuncer2005optimization, platzer2008flapping}.
The specific pitching and heaving amplitudes were determined by the optimizer
such that the impulse constraint is satisfied; as the impulse requirement is
increased, the magnitude of the pitch and plunge increase and
shedding off the leading edge is induced (Figure~\ref{fig:flap-vort-opt1})
\cite{oyama2009aerodynamic}. The time-dependent shape deformation slightly
reduces the magnitude of the vortices shedding off the leading edge, which can
be seen by comparing
Figures~\ref{fig:flap-vort-opt1}~and~\ref{fig:flap-vort-opt2}.

\section{Conclusion}\label{sec:conc}
This document discussed an Arbitrary Lagrangian-Eulerian formulation of
conservation laws on deforming domains and a high-order spatial discretization
using a discontinuous Galerkin finite element method, previously introduced
in \citep{persson2009dgdeform}. A diagonally implicit Runge-Kutta scheme was
used for high-order temporal integration with the benefit of achieving
stability beyond second-order and not coupling all intermediate states like
general, implicit Runge-Kutta methods. High-order spatial and temporal
integration of quantities of interest was done in a solver-consistent manner,
i.e., DG shape functions used for spatial integration and DIRK used for
temporal integration. The result is a globally high-order accurate
discretization of deforming domain conservation laws and relevant quantities of
interest. The fully discrete adjoint equations for this globally high-order
discretization were derived and the corresponding adjoint method introduced
for computing gradients of quantities of interest along the manifold of
high-order solutions of the conservation law. Implementation details regarding
computation of necessary partial derivatives and domain deformation were
included. The adjoint method was used in the context of gradient-based
optimization to solve two optimal shape and control problems constrained by
the isentropic, compressible Navier-Stokes equations, one of which
incorporated a physics-based constraint. The first optimization problem was
determination of the energetically optimal trajectory of an airfoil; the
optimization solver reduced the energy required to complete the specified
mission nearly an order of magnitude. The second optimization problem was the
determination of the energetically optimal flapping motion and airfoil shape
such that the $x$-directed impulse on the airfoil was fixed at a prescribed
value; the solver found a configuration that satisfied the constraint to
greater than $8$ digits of precision and reduced the required energy.

This work is an initial step toward problems of engineering and scientific
relevance. While the work extends to multiphysics problems defined on a
single domain \emph{without modification}, its extension to multiple domain
problems -- a canonical example being fluid-structure interaction problems --
requires development, particularly if high-order accuracy is sought
\cite{froehle2013high, froehle2014high, mishra2015time}. The extension to
chaotic problems will also have significant engineering implications as this
will be required to solve optimization problems governed by turbulent flow.
Definition of a meaningful optimization problem in this context becomes a
challenge as flow initialization \citep{zahr2016periodic} and sensitivity
computation becomes difficult \cite{wang2014least}.
Low-order and high-order methods will both suffer from these difficulties as
they arise from the chaotic nature of the physical problem.
For this type of automated optimization approach to be competitive in
engineering practice, the large cost associated with repeatedly solving the
underlying partial differential equation must be addressed; two promising
approaches are the use of inexact gradients \cite{heinkenschloss1995analysis}
to accelerate adjoint computations and adaptive reduced-order models
\cite{zahr2015progressive} to reduce the cost of the primal and sensitivity
computations.

\appendix
\section{Discrete Adjoint Equation Derivation via Lagrange Multiplier
         Interpretation}\label{sec:app-a}
In this section, an alternate derivation of the adjoint equations is provided
using constrained optimization theory and the test variables introduced
in Section~\ref{subsec:adj-deriv} enter as \emph{Lagrange multipliers}
corresponding to constraints in an auxiliary optimization problem.

Fix $\bar\mubold \in \Rbb^{N_\mubold}$ and consider the following
\emph{auxiliary optimization problem}, introduced purely for the derivation
of the adjoint equation for a particular functional, $F$,
\begin{equation} \label{opt:uns-disc-adj-alt}
  \begin{aligned}
    & \underset{\substack{\ubm^{(0)},~\dots,~\ubm^{(N_t)} \in \Rbb^{N_\ubm},\\
                        \kbm_1^{(1)},~\dots,~\kbm_s^{(N_t)} \in \Rbb^{N_\ubm}}}
                        {\text{minimize}}
    & & F(\ubm^{(0)},~\dots,~\ubm^{(N_t)},~\kbm_1^{(1)},~\dots,~\kbm_s^{(N_t)},~
          \bar\mubold) \\
    & \text{subject to}
    & &   \tilde\rbm^{(0)}(\ubm^{(0)},~\bar\mubold) = 0 \\
    & & & \tilde\rbm^{(n)}(\ubm^{(n-1)},~\ubm^{(n)},
                           \kbm_1^{(n)}, \dots, \kbm_s^{(n)},~
                           \bar\mubold) = 0 \\
    & & & \Rbm_i^{(n)}(\ubm^{(n-1)},~
                     \kbm_1^{(n)}, \dots, \kbm_i^{(n)},~
                     \bar\mubold) = 0.
  \end{aligned}
\end{equation}
Similar to Section ~\ref{subsec:adj-deriv}, $F$ corresponds to any output
functional whose gradient is desired.  In the auxiliary optimization problem,
\emph{only the state vectors} are taken as optimization variables since
$\bar\mubold$ is fixed.  If the discrete PDE has a unique solution,
the solution of (\ref{opt:uns-disc-adj-alt}) is equivalent to the solution of
the fully discrete PDE in (\ref{eqn:abs-govern-disc}).  In this case,
the optimization problem (\ref{opt:uns-disc-adj-alt}) is simply a feasibility
problem since the optimal solution is independent of the objective function.

The Lagrangian \cite{nocedal2006numerical} corresponding to the optimization
problem in (\ref{opt:uns-disc-adj-alt}) takes the form
\begin{equation} \label{eqn:lag}
  \begin{aligned}
    \Lcal(\ubm^{(0)},~\dots,~\ubm^{(N_t)},&~
          \kbm_1^{(1)},~\dots,~\kbm_s^{(N_t)},
          \lambdabold^{(0)},~\dots,~\lambdabold^{(N_t)},~
          \kappabold_1^{(1)},~\dots,~\kappabold_s^{(N_t)}) = \\
    &F(\ubm^{(0)},~\dots,~\ubm^{(N_t)},~\kbm_1^{(1)},~\dots,~\kbm_s^{(N_t)},~
          \bar\mubold) -
     \sum_{n=0}^{N_t} {\lambdabold^{(n)}}^T\tilde\rbm^{(n)} -
     \sum_{n=1}^{N_t}\sum_{i=1}^s {\kappabold_i^{(n)}}^T\Rbm_i^{(n)}.
  \end{aligned}
\end{equation}
Then, $(\ubm^{(0)*},~\dots,~\ubm^{(N_t)*},~
        \kbm_1^{(1)*},~\dots,~\kbm_s^{(N_t)*},
        \lambdabold^{(0)*},~\dots,~\lambdabold^{(N_t)*},~
        \kappabold_1^{(1)*},~\dots,~\kappabold_s^{(N_t)*})$
is a critical point of (\ref{opt:uns-disc-adj-alt}) if it makes the Lagrangian
stationary, i.e.
\begin{equation} \label{eqn:kkt}
  \begin{aligned}
    \pder{\Lcal}{\ubm^{(n)}} &= 0 \qquad & n = 0, \dots, N_t\\
    \pder{\Lcal}{\kbm_i^{(n)}} &= 0 \qquad & i = 1, \dots, s,
                                            ~n = 1, \dots, N_t \\
    \pder{\Lcal}{\lambdabold^{(n)}} &= 0 \qquad & n = 0, \dots, N_t\\
    \pder{\Lcal}{\kappabold_i^{(n)}} &= 0 \qquad & i = 1, \dots, s,
                                                  ~n = 1, \dots, N_t.
  \end{aligned}
\end{equation}
The equations in (\ref{eqn:kkt}) are the first-order optimality conditions,
or the Karush-Kuhn-Tucker (KKT) conditions. The conditions in (\ref{eqn:kkt})
on the derivatives of the Lagrangian with respect to the Lagrange multipliers
recover the constraints in (\ref{opt:uns-disc-adj-alt}).  From
(\ref{eqn:kkt}), the condition $\displaystyle{\pder{\Lcal}{\ubm^{(N_t)}}} = 0$
gives rise to
\begin{equation} \label{eqn:uns-disc-adj-alt-1}
  \pder{\tilde\rbm^{(N_t)}}{\ubm^{(N_t)}}^T\lambdabold^{(N_t)} =
                                                       \pder{F}{\ubm^{(N_t)}}.
\end{equation}
The conditions on the derivatives of the Lagrangian with respect to the 
state variable $\ubm^{(n-1)}$ for $n = 1, \dots, N_t$, result in
\begin{equation} \label{eqn:uns-disc-adj-alt-2}
\pder{\tilde\rbm^{(n)}}{\ubm^{(n-1)}}^T\lambdabold^{(n)}
   + \pder{\tilde\rbm^{(n-1)}}{\ubm^{(n-1)}}^T\lambdabold^{(n-1)} =
     \pder{F}{\ubm^{(n-1)}}^T - \sum_{i=1}^s \pder{\Rbm_i^{(n)}}{\ubm^{(n-1)}}^T
                                                            \kappabold_i^{(n)}.
\end{equation}
Finally, the conditions on the derivatives of the Lagrangian with respect to
the stage variables, $\kbm_i^{(n)}$ for $n = 1, \dots, N_t$ and
$i = 1, \dots, s$, are expanded as
\begin{equation}\label{eqn:uns-disc-adj-alt-3}
   \sum_{j=i}^s \pder{\Rbm_j^{(n)}}{\kbm_i^{(n)}}^T\kappabold_j^{(n)} =
    \pder{F}{\kbm_i^{(n)}} - \pder{\tilde\rbm^{(n)}}{\kbm_i^{(n)}}^T
                                                              \lambdabold^{(n)}.
\end{equation}
Equations (\ref{eqn:uns-disc-adj-alt-1}) - (\ref{eqn:uns-disc-adj-alt-3}),
which correspond to the first two conditions in the KKT system (\ref{eqn:kkt})
of the optimization problem in (\ref{opt:uns-disc-adj-alt}), are exactly
the adjoint equations (\ref{eqn:uns-disc-adj}) derived in
Section~\ref{subsec:adj-deriv}.

Substitution of the adjoint equation in (\ref{eqn:uns-disc-adj-alt-1}) -
(\ref{eqn:uns-disc-adj-alt-3}) into the expression for $\oder{F}{\mubold}$
in (\ref{eqn:uns-disc-adj-derive-2}) recovers expression
(\ref{eqn:funcl-grad-nosens}) for $\oder{F}{\mubold}$, which is independent
of the state sensitivities.  Finally, substitution of the expressions in
(\ref{eqn:abs-govern-disc}) for $\rbm^{(n)}$ and $\Rbm_i^{(n)}$ in the
adjoint equations (\ref{eqn:uns-disc-adj-alt-1}) -
(\ref{eqn:uns-disc-adj-alt-3}) and gradient reconstruction formula
(\ref{eqn:funcl-grad-nosens-dirk}) recover their specialization to the
case of a DIRK temporal discretization, i.e. (\ref{eqn:uns-disc-adj-dirk})
and (\ref{eqn:funcl-grad-nosens-dirk}), respectively.

\section*{Acknowledgments}
This work was supported in part by the Department of Energy Computational
Science Graduate Fellowship Program of the Office of Science and National
Nuclear Security Administration in the Department of Energy under contract
DE-FG02-97ER25308 (MZ), and by the Director, Office of
Science, Computational and Technology Research, U.S. Department of
Energy under contract number DE-AC02-05CH11231 (PP).
The content of this publication does not necessarily
reflect the position or policy of any of these supporters, and no official
endorsement should be inferred.

\bibliographystyle{ieeetr}
\bibliography{biblio}
\end{document}

%% file: mycommands.tex
\usepackage{bm}

\newcommand{\beq}{\begin{equation}}
\newcommand{\eeq}{\end{equation}}
\newcommand{\esplit}{\end{split}}
\newcommand{\beqalign}{\begin{array}{rl}}
\newcommand{\eeqalign}{\end{array}}

\newtheorem*{Rem}{Remarks}

\newcommand{\norm}[1]{\ensuremath{|| #1 ||}}

\newcommand{\pder}[2]{\ensuremath{\frac{\partial #1}{\partial #2}}} 

\newcommand{\oder}[2]{\ensuremath{\frac{\mathrm{d} #1}{\mathrm{d} #2}}} 

\newcommand{\Fcal}{\ensuremath{\mathcal{F}}}
\newcommand{\Gcal}{\ensuremath{\mathcal{G}}}

\newcommand{\Jcal}{\ensuremath{\mathcal{J}}}

\newcommand{\Lcal}{\ensuremath{\mathcal{L}}}

\newcommand{\Pcal}{\ensuremath{\mathcal{P}}}

\newcommand{\Wcal}{\ensuremath{\mathcal{W}}}

\newcommand{\Mbb}{\ensuremath{\mathbb{M} }}

\newcommand{\Rbb}{\ensuremath{\mathbb{R} }}

\newcommand\Fbm{{\ensuremath{\bm{F}}}}
\newcommand\Gbm{{\ensuremath{\bm{G}}}}

\newcommand\Ibm{{\ensuremath{\bm{I}}}}

\newcommand\Qbm{{\ensuremath{\bm{Q}}}}
\newcommand\Rbm{{\ensuremath{\bm{R}}}}

\newcommand\Ubm{{\ensuremath{\bm{U}}}}

\newcommand\Xbm{{\ensuremath{\bm{X}}}}

\newcommand\cbm{{\ensuremath{\bm{c}}}}

\newcommand\ebm{{\ensuremath{\bm{e}}}}
\newcommand\fbm{{\ensuremath{\bm{f}}}}

\newcommand\kbm{{\ensuremath{\bm{k}}}}

\newcommand\rbm{{\ensuremath{\bm{r}}}}

\newcommand\ubm{{\ensuremath{\bm{u}}}}
\newcommand\vbm{{\ensuremath{\bm{v}}}}

\newcommand\xbm{{\ensuremath{\bm{x}}}}

\newcommand\bbold{\ensuremath{\mathbf{b}}}

\newcommand\gbold{\ensuremath{\mathbf{g}}}

\newcommand\xbold{\ensuremath{\mathbf{x}}}

\newcommand\lambdabold{{\ensuremath{\boldsymbol{\lambda}}}}

\newcommand\mubold{{\ensuremath{\boldsymbol{\mu}}}}
\newcommand\kappabold{{\ensuremath{\boldsymbol{\kappa}}}}

\newcommand\varphibold{{\ensuremath{\boldsymbol{\varphi}}}}

\newcommand\Gammabold{{\ensuremath{\boldsymbol{\Gamma}}}}

\newcommand\zerobold{\ensuremath{\mathbf{0}}}

%% file: dat/pitch-setup.tab
\begin{tabular}{l|cccc|cccc} \toprule
& \multicolumn{4}{c}{PI} &
  \multicolumn{4}{c}{PII} \\\midrule
ID & $m_x$ & $m_y$ & $m_\theta$ & $N_\mubold$
    & $m_x$ & $m_y$ & $m_\theta$ & $N_\mubold$ \\\midrule
\rowcolor{black!00} 1 & 0 & 0 &   2 &   3 &   2 &   2 &   2 &   9 \\
\rowcolor{black!10} 2 & 0 & 0 &   6 &   7 &   6 &   6 &   6 &  21 \\
\rowcolor{black!00} 3 & 0 & 0 &  10 &  11 &  10 &  10 &  10 &  33 \\
\rowcolor{black!10} 4 & 0 & 0 &  15 &  16 &  15 &  15 &  15 &  48 \\
\rowcolor{black!00} 5 & 0 & 0 &  25 &  26 &  25 &  25 &  25 &  78 \\
\rowcolor{black!10} 6 & 0 & 0 &  50 &  51 &  50 &  50 &  50 & 153 \\
\rowcolor{black!00} 7 & 0 & 0 & 100 & 101 & 100 & 100 & 100 & 303 \\\bottomrule
\end{tabular}

%% file: tikz/pitch-findiff-rev.tikz
\begin{tikzpicture}
\begin{axis}[
 width=10cm,
 height=5cm,
 xmode=log,
 ymode=log,
 ymin=1.0e-11,
 ymax=1.0e-5,
 xtick={1.0e-1,1.0e-3,1.0e-5,1.0e-7,1.0e-9},
 ytick={1.0e-5,1.0e-7,1.0e-9,1.0e-11},
 xlabel={$\tau$},
 ylabel={$\norm{\oder{W}{\mubold} - \frac{\Delta W}{\Delta \mubold}}/
          \norm{\frac{\Delta W}{\Delta \mubold}}$}
]

\addplot [black, solid, thick, mark=triangle*, mark repeat=1]  table[x index=0, y index=4] {foil2d/cfd/flap00/findiff_rev/ascii/app1_mesh-how1-0-3_cfd-1p4-0p2-1000-0p72-0-1-_1-0__tdisc2-0p2-10-_0p5-1_-3.findiff_rel.cd4}; \label{line:findiff-w}

\end{axis}
\end{tikzpicture}

%% file: tikz/pitch-spatial-conv-qoi.tikz
\begin{tikzpicture}
\begin{groupplot}[
    group style={
        group name=findiff plots,
        group size=2 by 2,
        horizontal sep=2.0cm
    },
    width=7.5cm,
    height=4.5cm,
    xmode=log,
    ymode=log,
]
\nextgroupplot[xticklabels={}, ylabel=$\left|J_x-J_x^*\right|$]
\addplot [black, solid, thick, mark=o]  table[x expr=\thisrowno{0}^0.5, y index = 1] {foil2d/cfd/pitch00/conv/post/pord1_nstage3_nstep1000.space}; \label{line:p1}
\addplot [blue, solid, thick, mark=square]  table[x expr=\thisrowno{0}^0.5, y index = 1] {foil2d/cfd/pitch00/conv/post/pord2_nstage3_nstep1000.space}; \label{line:p2}
\addplot [red, solid, thick, mark=triangle]  table[x expr=\thisrowno{0}^0.5, y index = 1] {foil2d/cfd/pitch00/conv/post/pord3_nstage3_nstep1000.space}; \label{line:p3}

\nextgroupplot[xticklabels={}]
\addplot [black, solid, thick, mark=o]  table[x expr=12*\thisrowno{0}, y index = 1] {foil2d/cfd/pitch00/conv/post/pord1_nstage3_nstep1000.space.timings};
\addplot [blue, solid, thick, mark=square]  table[x expr=12*\thisrowno{0}, y index = 1] {foil2d/cfd/pitch00/conv/post/pord2_nstage3_nstep1000.space.timings};
\addplot [red, solid, thick, mark=triangle]  table[x expr=12*\thisrowno{0}, y index = 1] {foil2d/cfd/pitch00/conv/post/pord3_nstage3_nstep1000.space.timings};

\nextgroupplot[xlabel={$\sqrt{N_\ubm}$}, ylabel=$\left|W-W^*\right|$]
\addplot [black, solid, thick, mark=o]  table[x expr=\thisrowno{0}^0.5, y index = 4] {foil2d/cfd/pitch00/conv/post/pord1_nstage3_nstep1000.space};
\addplot [blue, solid, thick, mark=square]  table[x expr=\thisrowno{0}^0.5, y index = 4] {foil2d/cfd/pitch00/conv/post/pord2_nstage3_nstep1000.space};
\addplot [red, solid, thick, mark=triangle]  table[x expr=\thisrowno{0}^0.5, y index = 4] {foil2d/cfd/pitch00/conv/post/pord3_nstage3_nstep1000.space};

\nextgroupplot[xlabel={CPU time (core-seconds)}]
\addplot [black, solid, thick, mark=o]  table[x expr=12*\thisrowno{0}, y index = 4] {foil2d/cfd/pitch00/conv/post/pord1_nstage3_nstep1000.space.timings};
\addplot [blue, solid, thick, mark=square]  table[x expr=12*\thisrowno{0}, y index = 4] {foil2d/cfd/pitch00/conv/post/pord2_nstage3_nstep1000.space.timings};
\addplot [red, solid, thick, mark=triangle]  table[x expr=12*\thisrowno{0}, y index = 4] {foil2d/cfd/pitch00/conv/post/pord3_nstage3_nstep1000.space.timings};

\end{groupplot}
\end{tikzpicture}

%% file: tikz/pitch-traj-x.tikz
\begin{tikzpicture}

\begin{axis}[
  width=7.5cm,
  height=4.5cm,
  xlabel=time,
  ylabel=$x(t)$]
\addplot [black, solid, thick, mark=o, mark repeat=20]  table[x expr=0.04*\coordindex, y index=1] {foil2d/cfd/pitch00/opt/freeze_x/testcase00/ascii/app1_mesh-how1-0-3_cfd-1p4-0p2-1000-0p72-0-1-_0-0__tdisc0-4-100-3_it000.xyth};
\addplot [red, solid, thick, mark=square, mark repeat=20]  table[x expr=0.04*\coordindex, y index=1] {foil2d/cfd/pitch00/opt/freeze_xy/testcase14/ascii/app1_mesh-how1-0-3_cfd-1p4-0p2-1000-0p72-0-1-_0-0__tdisc0-4-100-3_it064.xyth};
\addplot [blue, solid, thick, mark=triangle, mark repeat=20]  table[x expr=0.04*\coordindex, y index=1] {foil2d/cfd/pitch00/opt/freeze_x/testcase14/ascii/app1_mesh-how1-0-3_cfd-1p4-0p2-1000-0p72-0-1-_0-0__tdisc0-4-100-3_it045.xyth};
\end{axis}

\end{tikzpicture}

%% file: tikz/pitch-traj-y.tikz
\begin{tikzpicture}

\begin{axis}[
  width=8cm,
  height=4.5cm,
  xlabel=time,
  ylabel=$y(t)$]
\addplot [black, solid, thick, mark=o, mark repeat=20]  table[x expr=0.04*\coordindex, y index=2] {foil2d/cfd/pitch00/opt/freeze_x/testcase00/ascii/app1_mesh-how1-0-3_cfd-1p4-0p2-1000-0p72-0-1-_0-0__tdisc0-4-100-3_it000.xyth};
\addplot [red, solid, thick, mark=square, mark repeat=20]  table[x expr=0.04*\coordindex, y index=2] {foil2d/cfd/pitch00/opt/freeze_xy/testcase14/ascii/app1_mesh-how1-0-3_cfd-1p4-0p2-1000-0p72-0-1-_0-0__tdisc0-4-100-3_it064.xyth};
\addplot [blue, solid, thick, mark=triangle, mark repeat=20]  table[x expr=0.04*\coordindex, y index=2] {foil2d/cfd/pitch00/opt/freeze_x/testcase14/ascii/app1_mesh-how1-0-3_cfd-1p4-0p2-1000-0p72-0-1-_0-0__tdisc0-4-100-3_it045.xyth};
\end{axis}

\end{tikzpicture}

%% file: tikz/pitch-traj-th.tikz
\begin{tikzpicture}

\begin{axis}[
  width=8cm,
  height=4.5cm,
  xlabel=time,
  ylabel=$\theta(t)$]
\addplot [black, solid, thick, mark=o, mark repeat=20]  table[x expr=0.04*\coordindex, y index=3] {foil2d/cfd/pitch00/opt/freeze_x/testcase00/ascii/app1_mesh-how1-0-3_cfd-1p4-0p2-1000-0p72-0-1-_0-0__tdisc0-4-100-3_it000.xyth}; \label{line:pitch-init}
\addplot [red, solid, thick, mark=square, mark repeat=20]  table[x expr=0.04*\coordindex, y index=3] {foil2d/cfd/pitch00/opt/freeze_xy/testcase14/ascii/app1_mesh-how1-0-3_cfd-1p4-0p2-1000-0p72-0-1-_0-0__tdisc0-4-100-3_it064.xyth}; \label{line:pitch-opt1}
\addplot [blue, solid, thick, mark=triangle, mark repeat=20]  table[x expr=0.04*\coordindex, y index=3] {foil2d/cfd/pitch00/opt/freeze_x/testcase14/ascii/app1_mesh-how1-0-3_cfd-1p4-0p2-1000-0p72-0-1-_0-0__tdisc0-4-100-3_it045.xyth}; \label{line:pitch-opt2}
\end{axis}

\end{tikzpicture}

%% file: tikz/pitch-traj-qoi.tikz
\begin{tikzpicture}
\begin{groupplot}[
    group style={
        group name=findiff plots,
        group size=2 by 3,
        horizontal sep=2.4cm
    },
    width=7.5cm,
    height=4.5cm,
    yticklabel style={
        /pgf/number format/fixed
    },
    scaled y ticks=false
]
\nextgroupplot[xticklabels={}, ylabel=$\Fcal_x^h$]
\addplot [black, solid, thick, mark=o, mark repeat=20, smooth]  table[x expr=0.04*\coordindex, y index=0] {foil2d/cfd/pitch00/opt/freeze_xy/testcase00/ascii/app1_mesh-how1-0-3_cfd-1p4-0p2-1000-0p72-0-1-_0-0__tdisc0-4-100-3_it000.allpost.instant};
\addplot [red, solid, thick, mark=square, mark repeat=20, smooth]  table[x expr=0.04*\coordindex, y index=0] {foil2d/cfd/pitch00/opt/freeze_xy/testcase14/ascii/app1_mesh-how1-0-3_cfd-1p4-0p2-1000-0p72-0-1-_0-0__tdisc0-4-100-3_it063.allpost.instant};
\addplot [blue, solid, thick, mark=triangle, mark repeat=20, smooth]  table[x expr=0.04*\coordindex, y index=0] {foil2d/cfd/pitch00/opt/freeze_x/testcase14/ascii/app1_mesh-how1-0-3_cfd-1p4-0p2-1000-0p72-0-1-_0-0__tdisc0-4-100-3_it045.allpost.instant};

\nextgroupplot[xticklabels={}, ylabel=$\Fcal_y^h$]
\addplot [black, solid, thick, mark=o, mark repeat=20, smooth]  table[x expr=0.04*\coordindex, y index=1] {foil2d/cfd/pitch00/opt/freeze_xy/testcase00/ascii/app1_mesh-how1-0-3_cfd-1p4-0p2-1000-0p72-0-1-_0-0__tdisc0-4-100-3_it000.allpost.instant};
\addplot [red, solid, thick, mark=square, mark repeat=20, smooth]  table[x expr=0.04*\coordindex, y index=1] {foil2d/cfd/pitch00/opt/freeze_xy/testcase14/ascii/app1_mesh-how1-0-3_cfd-1p4-0p2-1000-0p72-0-1-_0-0__tdisc0-4-100-3_it063.allpost.instant};
\addplot [blue, solid, thick, mark=triangle, mark repeat=20, smooth]  table[x expr=0.04*\coordindex, y index=1] {foil2d/cfd/pitch00/opt/freeze_x/testcase14/ascii/app1_mesh-how1-0-3_cfd-1p4-0p2-1000-0p72-0-1-_0-0__tdisc0-4-100-3_it045.allpost.instant};

\nextgroupplot[xticklabels={}, ylabel=$\Pcal^h$]
\addplot [black, solid, thick, mark=o, mark repeat=20, smooth]  table[x expr=0.04*\coordindex, y index=3] {foil2d/cfd/pitch00/opt/freeze_xy/testcase00/ascii/app1_mesh-how1-0-3_cfd-1p4-0p2-1000-0p72-0-1-_0-0__tdisc0-4-100-3_it000.allpost.instant};
\addplot [red, solid, thick, mark=square, mark repeat=20, smooth]  table[x expr=0.04*\coordindex, y index=3] {foil2d/cfd/pitch00/opt/freeze_xy/testcase14/ascii/app1_mesh-how1-0-3_cfd-1p4-0p2-1000-0p72-0-1-_0-0__tdisc0-4-100-3_it063.allpost.instant};
\addplot [blue, solid, thick, mark=triangle, mark repeat=20, smooth]  table[x expr=0.04*\coordindex, y index=3] {foil2d/cfd/pitch00/opt/freeze_x/testcase14/ascii/app1_mesh-how1-0-3_cfd-1p4-0p2-1000-0p72-0-1-_0-0__tdisc0-4-100-3_it045.allpost.instant};

\nextgroupplot[xticklabels={}, ylabel=$\Pcal_x^h$]
\addplot [black, solid, thick, mark=o, mark repeat=20, smooth]  table[x expr=0.04*\coordindex, y index=4] {foil2d/cfd/pitch00/opt/freeze_xy/testcase00/ascii/app1_mesh-how1-0-3_cfd-1p4-0p2-1000-0p72-0-1-_0-0__tdisc0-4-100-3_it000.allpost.instant};
\addplot [red, solid, thick, mark=square, mark repeat=20, smooth]  table[x expr=0.04*\coordindex, y index=4] {foil2d/cfd/pitch00/opt/freeze_xy/testcase14/ascii/app1_mesh-how1-0-3_cfd-1p4-0p2-1000-0p72-0-1-_0-0__tdisc0-4-100-3_it063.allpost.instant};
\addplot [blue, solid, thick, mark=triangle, mark repeat=20, smooth]  table[x expr=0.04*\coordindex, y index=4] {foil2d/cfd/pitch00/opt/freeze_x/testcase14/ascii/app1_mesh-how1-0-3_cfd-1p4-0p2-1000-0p72-0-1-_0-0__tdisc0-4-100-3_it045.allpost.instant};

\nextgroupplot[xlabel={time}, ylabel=$\Pcal_y^h$]
\addplot [black, solid, thick, mark=o, mark repeat=20, smooth]  table[x expr=0.04*\coordindex, y index=5] {foil2d/cfd/pitch00/opt/freeze_xy/testcase00/ascii/app1_mesh-how1-0-3_cfd-1p4-0p2-1000-0p72-0-1-_0-0__tdisc0-4-100-3_it000.allpost.instant};
\addplot [red, solid, thick, mark=square, mark repeat=20, smooth]  table[x expr=0.04*\coordindex, y index=5] {foil2d/cfd/pitch00/opt/freeze_xy/testcase14/ascii/app1_mesh-how1-0-3_cfd-1p4-0p2-1000-0p72-0-1-_0-0__tdisc0-4-100-3_it063.allpost.instant};
\addplot [blue, solid, thick, mark=triangle, mark repeat=20, smooth]  table[x expr=0.04*\coordindex, y index=5] {foil2d/cfd/pitch00/opt/freeze_x/testcase14/ascii/app1_mesh-how1-0-3_cfd-1p4-0p2-1000-0p72-0-1-_0-0__tdisc0-4-100-3_it045.allpost.instant};

\nextgroupplot[xlabel={time}, ylabel=$\Pcal_\theta^h$]
\addplot [black, solid, thick, mark=o, mark repeat=20, smooth]  table[x expr=0.04*\coordindex, y index=6] {foil2d/cfd/pitch00/opt/freeze_xy/testcase00/ascii/app1_mesh-how1-0-3_cfd-1p4-0p2-1000-0p72-0-1-_0-0__tdisc0-4-100-3_it000.allpost.instant};
\addplot [red, solid, thick, mark=square, mark repeat=20, smooth]  table[x expr=0.04*\coordindex, y index=6] {foil2d/cfd/pitch00/opt/freeze_xy/testcase14/ascii/app1_mesh-how1-0-3_cfd-1p4-0p2-1000-0p72-0-1-_0-0__tdisc0-4-100-3_it063.allpost.instant};
\addplot [blue, solid, thick, mark=triangle, mark repeat=20, smooth]  table[x expr=0.04*\coordindex, y index=6] {foil2d/cfd/pitch00/opt/freeze_x/testcase14/ascii/app1_mesh-how1-0-3_cfd-1p4-0p2-1000-0p72-0-1-_0-0__tdisc0-4-100-3_it045.allpost.instant};
\end{groupplot}
\end{tikzpicture}

%% file: tikz/pitch-conv.tikz
\begin{tikzpicture}
\begin{groupplot}[
    group style={
        group name=findiff plots,
        group size=2 by 1,
        horizontal sep=2.4cm
    },
    width=7.5cm,
    height=4.5cm,
]

\nextgroupplot[xlabel=iteration, ylabel={$W$}]
\addplot [red, solid, thick, mark=square, mark repeat=5]  table[x expr=\coordindex, y index=4]  {foil2d/cfd/pitch00/opt/freeze_xy/app1_mesh-how1-0-3_cfd-1p4-0p2-1000-0p72-0-1-_0-0__tdisc0-4-100-3.allpost.integrated.allconv};
\addplot [blue, solid, thick, mark=triangle, mark repeat=5]  table[x expr=\coordindex, y index=4]  {foil2d/cfd/pitch00/opt/freeze_x/app1_mesh-how1-0-3_cfd-1p4-0p2-1000-0p72-0-1-_0-0__tdisc0-4-100-3.allpost.integrated.allconv};

\nextgroupplot[xlabel=$N_\mubold$, ylabel={$W$}, xmode=log]
\addplot [red, solid, thick, mark=square]  table[x expr=\thisrowno{0}+1, y index=4] {foil2d/cfd/pitch00/opt/freeze_xy/app1_mesh-how1-0-3_cfd-1p4-0p2-1000-0p72-0-1-_0-0__tdisc0-4-100-3.allpost.integrated.paramconv};
\addplot [blue, solid, thick, mark=triangle]  table[x expr=3*(\thisrowno{0}+1), y index=4] {foil2d/cfd/pitch00/opt/freeze_x/app1_mesh-how1-0-3_cfd-1p4-0p2-1000-0p72-0-1-_0-0__tdisc0-4-100-3.allpost.integrated.paramconv};

\end{groupplot}
\end{tikzpicture}

%% file: dat/pitch-summ.tab
\begin{tabular}{lccccccc} \toprule
& \multicolumn{7}{c}{PI} \\\midrule
ID & 1 & 2 & 3 & 4 & 5 & 6 & 7 \\\midrule
\rowcolor{black!00} $W$   & -2.1951e-01 & -1.5881e-01 & -1.5358e-01 & -1.5128e-01 & -1.5026e-01 & -1.4950e-01 & -1.4924e-01 \\
\rowcolor{black!10} $W_x$ &  8.1329e-02 &  5.6090e-02 &  4.9543e-02 &  4.5924e-02 &  4.5085e-02 &  4.4712e-02 &  4.4707e-02 \\
\rowcolor{black!00} $W_y$ & -2.3460e-01 & -1.8153e-01 & -1.7122e-01 & -1.6544e-01 & -1.6374e-01 & -1.6298e-01 & -1.6294e-01 \\
\rowcolor{black!10} $W_\theta$ &  6.6234e-02 &  3.3370e-02 &  3.1906e-02 &  3.1768e-02 &  3.1604e-02 &  3.1223e-02 &  3.1010e-02 \\
\rowcolor{black!00} $F_x$ & -1.9234e-01 & -1.3123e-01 & -1.1886e-01 & -1.1136e-01 & -1.0912e-01 & -1.0810e-01 & -1.0800e-01 \\
\rowcolor{black!10} $F_y$ & -5.1539e-01 & -3.1711e-01 & -3.1816e-01 & -3.1877e-01 & -3.2551e-01 & -3.2959e-01 & -3.3063e-01 \\\bottomrule
& \multicolumn{7}{c}{PII} \\\midrule
ID & 1 & 2 & 3 & 4 & 5 & 6 & 7 \\\midrule
\rowcolor{black!00} $W$   & -1.7357e-01 &  -1.2095e-01 & -1.1733e-01
                          & -1.1629e-01 & -1.1603e-01 & -1.1557e-01
                          & -1.1502e-01 \\
\rowcolor{black!10} $W_x$ & 9.6487e-03 & -1.4123e-02 & -1.4328e-02
                          & -1.4967e-02 & -1.5021e-02 & -1.5061e-02
                          & -1.5027e-02 \\
\rowcolor{black!00} $W_y$ & -1.1041e-01 & -6.2238e-02 & -6.1036e-02
                          & -6.0425e-02 & -6.0032e-02 & -5.9489e-02
                          & -5.9245e-02 \\
\rowcolor{black!10} $W_\theta$ & 7.2807e-02 & 4.4585e-02 & 4.1963e-02 & 4.0895e-02 & 4.0980e-02 & 4.1023e-02 & 4.0749e-02 \\
\rowcolor{black!00} $F_x$ & -4.1265e-02 & 2.8091e-02 & 2.7677e-02 & 2.9596e-02
                          & 2.9848e-02 & 3.0231e-02 & 3.0221e-02 \\
\rowcolor{black!10} $F_y$ & -3.2231e-01 & -1.064e-01 & -1.0806e-01 & -1.0890e-01
                          & -1.1343e-01 & -1.1626e-01 & -1.1764e-01 \\\bottomrule
\end{tabular}

%% file: dat/flap-setup.tab
\begin{tabular}{cccc|cccc} \toprule
\multicolumn{4}{c}{FI} &
\multicolumn{4}{c}{FII} \\\midrule
   $m_y$ & $m_\theta$ & $m_c$ & $N_\mubold$ &
   $m_y$ & $m_\theta$ & $m_c$ & $N_\mubold$ \\\midrule
\rowcolor{black!00}   4 & 4 & 0 & 6 & 4 & 4 & 4 & 9 \\\bottomrule
\end{tabular}

%% file: tikz/flap-y-hist.tikz
\begin{tikzpicture}

\begin{axis}[
name=flap-y-hist,
width=7.5cm,
height=4.5cm,
xmin=10,
xmax=15,
ymin=-1.25,
ymax=1.25,
xlabel=$t$,
ylabel=$y(t)$]

\addplot [black, solid, thick]  table[x index=0, y index=2]  {foil2d/cfd/flap01/opt/thrust0/freeze_x/testcase00/ascii/app1_mesh-how1-0-3_cfd-1p4-0p2-1000-0p72-0-1-_1-0__tdisc2-0p2-100-_3-1_-3_it000.xyth};
\addplot [red, solid, thick, mark=*, mark repeat=10]  table[x index=0, y index=2] {foil2d/cfd/flap01/opt/thrust0/freeze_x/testcase00/ascii/app1_mesh-how1-0-3_cfd-1p4-0p2-1000-0p72-0-1-_1-0__tdisc2-0p2-100-_3-1_-3_it042.xyth};
\addplot [blue, solid, thick, mark=square*, mark repeat=10]  table[x index=0, y index=2] {foil2d/cfd/flap01/opt/thrust1/freeze_x/testcase00/ascii/app1_mesh-how1-0-3_cfd-1p4-0p2-1000-0p72-0-1-_1-0__tdisc2-0p2-100-_3-1_-3_it047.xyth};
\addplot [magenta, solid, thick, mark=triangle*, mark repeat=10] table[x index=0, y index=2] {foil2d/cfd/flap01/opt/thrust2p5/freeze_x/testcase00/ascii/app1_mesh-how1-0-3_cfd-1p4-0p2-1000-0p72-0-1-_1-0__tdisc2-0p2-100-_3-1_-3_it062.xyth};

\addplot [red, dashed, thick, mark=*, mark repeat=10, mark options={solid}]  table[x index=0, y index=2] {foil2d/cfd/flap01morph01/opt/thrust0/freeze_x/testcase00/ascii/app1_mesh-how1-0-3_cfd-1p4-0p2-1000-0p72-0-1-_1-0__tdisc2-0p2-100-_3-1_-3_it091.xyth};
\addplot [blue, dashed, thick, mark=square*, mark repeat=10, mark options={solid}]  table[x index=0, y index=2] {foil2d/cfd/flap01morph01/opt/thrust1/freeze_x/testcase00/ascii/app1_mesh-how1-0-3_cfd-1p4-0p2-1000-0p72-0-1-_1-0__tdisc2-0p2-100-_3-1_-3_it075.xyth};
\addplot [magenta, dashed, thick, mark=triangle*, mark repeat=10, mark options={solid}]  table[x index=0, y index=2] {foil2d/cfd/flap01morph01/opt/thrust2p5/freeze_x/testcase00/ascii/app1_mesh-how1-0-3_cfd-1p4-0p2-1000-0p72-0-1-_1-0__tdisc2-0p2-100-_3-1_-3_it100.xyth};

\end{axis}

\end{tikzpicture}

%% file: tikz/flap-th-hist.tikz
\begin{tikzpicture}

\begin{axis}[
name=flap-th-hist,
width=7.5cm,
height=4.5cm,
xmin=10,
xmax=15,
ymin=-1,
ymax=1,
xlabel=$t$,
ylabel=$\theta(t)$]
\addplot [black, solid, thick]  table[x index=0, y index=3] {foil2d/cfd/flap01/opt/thrust0/freeze_x/testcase00/ascii/app1_mesh-how1-0-3_cfd-1p4-0p2-1000-0p72-0-1-_1-0__tdisc2-0p2-100-_3-1_-3_it000.xyth}; \label{line:flap-init}

\addplot [red, solid, thick, mark=*, mark repeat=10]  table[x index=0, y index=3] {foil2d/cfd/flap01/opt/thrust0/freeze_x/testcase00/ascii/app1_mesh-how1-0-3_cfd-1p4-0p2-1000-0p72-0-1-_1-0__tdisc2-0p2-100-_3-1_-3_it042.xyth}; \label{line:flap-opt1-0}
\addplot [blue, solid, thick, mark=square*, mark repeat=10]  table[x index=0, y index=3] {foil2d/cfd/flap01/opt/thrust1/freeze_x/testcase00/ascii/app1_mesh-how1-0-3_cfd-1p4-0p2-1000-0p72-0-1-_1-0__tdisc2-0p2-100-_3-1_-3_it047.xyth}; \label{line:flap-opt1-1}
\addplot [magenta, solid, thick, mark=triangle*, mark repeat=10]  table[x index=0, y index=3] {foil2d/cfd/flap01/opt/thrust2p5/freeze_x/testcase00/ascii/app1_mesh-how1-0-3_cfd-1p4-0p2-1000-0p72-0-1-_1-0__tdisc2-0p2-100-_3-1_-3_it062.xyth}; \label{line:flap-opt1-2.5}

\addplot [red, dashed, thick, mark=*, mark repeat=10, mark options={solid}]  table[x index=0, y index=3] {foil2d/cfd/flap01morph01/opt/thrust0/freeze_x/testcase00/ascii/app1_mesh-how1-0-3_cfd-1p4-0p2-1000-0p72-0-1-_1-0__tdisc2-0p2-100-_3-1_-3_it091.xyth}; \label{line:flap-opt2-0}
\addplot [blue, dashed, thick, mark=square*, mark repeat=10, mark options={solid}]  table[x index=0, y index=3] {foil2d/cfd/flap01morph01/opt/thrust1/freeze_x/testcase00/ascii/app1_mesh-how1-0-3_cfd-1p4-0p2-1000-0p72-0-1-_1-0__tdisc2-0p2-100-_3-1_-3_it075.xyth}; \label{line:flap-opt2-1}
\addplot [magenta, dashed, thick, mark=triangle*, mark repeat=10, mark options={solid}]  table[x index=0, y index=3] {foil2d/cfd/flap01morph01/opt/thrust2p5/freeze_x/testcase00/ascii/app1_mesh-how1-0-3_cfd-1p4-0p2-1000-0p72-0-1-_1-0__tdisc2-0p2-100-_3-1_-3_it100.xyth}; \label{line:flap-opt2-2.5}

\end{axis}
\end{tikzpicture}

%% file: tikz/flap-camber-hist.tikz
\begin{tikzpicture}

\begin{axis}[
name=flap-th-hist,
width=7.5cm,
height=4.5cm,
xmin=10,
xmax=15,
ymin=-0.41,
ymax=0.41,
xlabel=$t$,
ylabel=$c(t)$]
\addplot [black, solid, thick]  table[x index=0, y index=3] {foil2d/cfd/flap01morph01/opt/thrust0/freeze_x/testcase00/ascii/app1_mesh-how1-0-3_cfd-1p4-0p2-1000-0p72-0-1-_1-0__tdisc2-0p2-100-_3-1_-3_it000.lhc};

\addplot [red, dashed, thick, mark=*, mark repeat=10, mark options={solid}]  table[x index=0, y index=3] {foil2d/cfd/flap01morph01/opt/thrust0/freeze_x/testcase00/ascii/app1_mesh-how1-0-3_cfd-1p4-0p2-1000-0p72-0-1-_1-0__tdisc2-0p2-100-_3-1_-3_it091.lhc};
\addplot [blue, dashed, thick, mark=square*, mark repeat=10, mark options={solid}]  table[x index=0, y index=3] {foil2d/cfd/flap01morph01/opt/thrust1/freeze_x/testcase00/ascii/app1_mesh-how1-0-3_cfd-1p4-0p2-1000-0p72-0-1-_1-0__tdisc2-0p2-100-_3-1_-3_it075.lhc};
\addplot [magenta, dashed, thick, mark=triangle*, mark repeat=10, mark options={solid}]  table[x index=0, y index=3] {foil2d/cfd/flap01morph01/opt/thrust2p5/freeze_x/testcase00/ascii/app1_mesh-how1-0-3_cfd-1p4-0p2-1000-0p72-0-1-_1-0__tdisc2-0p2-100-_3-1_-3_it100.lhc};

\end{axis}
\end{tikzpicture}

%% file: tikz/flap-hist.tikz
\begin{tikzpicture}
\begin{groupplot}[
    group style={
        group name=findiff plots,
        group size=2 by 1,
        horizontal sep=2.4cm
    },
    width=7.5cm,
    height=4.5cm,
] 

\nextgroupplot[xmin=10, xmax=15, xlabel={time}, ylabel={$\Pcal^h$}]
\addplot [black, solid, thick, smooth]  table[x index=0, y index=4]  {foil2d/cfd/flap01/opt/thrust0/freeze_x/testcase00/ascii/app1_mesh-how1-0-3_cfd-1p4-0p2-1000-0p72-0-1-_1-0__tdisc2-0p2-100-_3-1_-3_it000.allpost.instant.wit};

\addplot [red, solid, thick, mark=*, mark repeat=10, smooth]  table[x index=0, y index=4]  {foil2d/cfd/flap01/opt/thrust0/freeze_x/testcase00/ascii/app1_mesh-how1-0-3_cfd-1p4-0p2-1000-0p72-0-1-_1-0__tdisc2-0p2-100-_3-1_-3_it042.allpost.instant.wit};
\addplot [blue, solid, thick, mark=square*, mark repeat=10, smooth]  table[x index=0, y index=4]  {foil2d/cfd/flap01/opt/thrust1/freeze_x/testcase00/ascii/app1_mesh-how1-0-3_cfd-1p4-0p2-1000-0p72-0-1-_1-0__tdisc2-0p2-100-_3-1_-3_it047.allpost.instant.wit};
\addplot [magenta, solid, thick, mark=triangle*, mark repeat=10, smooth]  table[x index=0, y index=4]  {foil2d/cfd/flap01/opt/thrust2p5/freeze_x/testcase00/ascii/app1_mesh-how1-0-3_cfd-1p4-0p2-1000-0p72-0-1-_1-0__tdisc2-0p2-100-_3-1_-3_it062.allpost.instant.wit};

\addplot [red, dashed, thick, mark=*, mark repeat=10, smooth, mark options={solid}]  table[x index=0, y index=4]  {foil2d/cfd/flap01morph01/opt/thrust0/freeze_x/testcase00/ascii/app1_mesh-how1-0-3_cfd-1p4-0p2-1000-0p72-0-1-_1-0__tdisc2-0p2-100-_3-1_-3_it091.allpost.instant.wit};
\addplot [blue, dashed, thick, mark=square*, mark repeat=10, smooth, mark options={solid}]  table[x index=0, y index=4]  {foil2d/cfd/flap01morph01/opt/thrust1/freeze_x/testcase00/ascii/app1_mesh-how1-0-3_cfd-1p4-0p2-1000-0p72-0-1-_1-0__tdisc2-0p2-100-_3-1_-3_it075.allpost.instant.wit};
\addplot [magenta, dashed, thick, mark=triangle*, mark repeat=10, smooth, mark options={solid}]  table[x index=0, y index=4]  {foil2d/cfd/flap01morph01/opt/thrust2p5/freeze_x/testcase00/ascii/app1_mesh-how1-0-3_cfd-1p4-0p2-1000-0p72-0-1-_1-0__tdisc2-0p2-100-_3-1_-3_it100.allpost.instant.wit};

\nextgroupplot[xmin=10, xmax=15, xlabel={time}, ylabel={$\Fcal_x^h$}]
\addplot [black, solid, thick, smooth]  table[x index=0, y index=1]  {foil2d/cfd/flap01/opt/thrust0/freeze_x/testcase00/ascii/app1_mesh-how1-0-3_cfd-1p4-0p2-1000-0p72-0-1-_1-0__tdisc2-0p2-100-_3-1_-3_it000.allpost.instant.wit};

\addplot [red, solid, thick, mark=*, mark repeat=10, smooth]  table[x index=0, y index=1]  {foil2d/cfd/flap01/opt/thrust0/freeze_x/testcase00/ascii/app1_mesh-how1-0-3_cfd-1p4-0p2-1000-0p72-0-1-_1-0__tdisc2-0p2-100-_3-1_-3_it042.allpost.instant.wit};
\addplot [blue, solid, thick, mark=square*, mark repeat=10, smooth]  table[x index=0, y index=1]  {foil2d/cfd/flap01/opt/thrust1/freeze_x/testcase00/ascii/app1_mesh-how1-0-3_cfd-1p4-0p2-1000-0p72-0-1-_1-0__tdisc2-0p2-100-_3-1_-3_it047.allpost.instant.wit};
\addplot [magenta, solid, thick, mark=triangle*, mark repeat=10, smooth]  table[x index=0, y index=1]  {foil2d/cfd/flap01/opt/thrust2p5/freeze_x/testcase00/ascii/app1_mesh-how1-0-3_cfd-1p4-0p2-1000-0p72-0-1-_1-0__tdisc2-0p2-100-_3-1_-3_it062.allpost.instant.wit};

\addplot [red, dashed, thick, mark=*, mark repeat=10, smooth]  table[x index=0, y index=1]  {foil2d/cfd/flap01morph01/opt/thrust0/freeze_x/testcase00/ascii/app1_mesh-how1-0-3_cfd-1p4-0p2-1000-0p72-0-1-_1-0__tdisc2-0p2-100-_3-1_-3_it061.allpost.instant.wit};
\addplot [blue, dashed, thick, mark=*, mark repeat=10, smooth]  table[x index=0, y index=1]  {foil2d/cfd/flap01morph01/opt/thrust1/freeze_x/testcase00/ascii/app1_mesh-how1-0-3_cfd-1p4-0p2-1000-0p72-0-1-_1-0__tdisc2-0p2-100-_3-1_-3_it024.allpost.instant.wit};
\addplot [magenta, dashed, thick, mark=*, mark repeat=10, smooth]  table[x index=0, y index=1]  {foil2d/cfd/flap01morph01/opt/thrust2p5/freeze_x/testcase00/ascii/app1_mesh-how1-0-3_cfd-1p4-0p2-1000-0p72-0-1-_1-0__tdisc2-0p2-100-_3-1_-3_it076.allpost.instant.wit};

\end{groupplot}
\end{tikzpicture}

%% file: tikz/flap-conv.tikz
\begin{tikzpicture}
\begin{groupplot}[
    group style={
        group name=findiff plots,
        group size=2 by 1,
        horizontal sep=2.4cm
    },
    width=7.5cm,
    height=4.5cm,
] 

\nextgroupplot[xmin=0, xmax=60, xlabel={iteration}, ylabel=$W$]
\addplot [red, solid, thick, mark=*, mark repeat=5]  table[x index=0, y index=4]  {foil2d/cfd/flap01/opt/thrust0/freeze_x/testcase00/ascii/app1_mesh-how1-0-3_cfd-1p4-0p2-1000-0p72-0-1-_1-0__tdisc2-0p2-100-_3-1_-3.allpost.instant.conv};
\addplot [blue, solid, thick, mark=square*, mark repeat=5]  table[x index=0, y index=4]  {foil2d/cfd/flap01/opt/thrust1/freeze_x/testcase00/ascii/app1_mesh-how1-0-3_cfd-1p4-0p2-1000-0p72-0-1-_1-0__tdisc2-0p2-100-_3-1_-3.allpost.instant.conv};
\addplot [magenta, solid, thick, mark=triangle*, mark repeat=5]  table[x index=0, y index=4]  {foil2d/cfd/flap01/opt/thrust2p5/freeze_x/testcase00/ascii/app1_mesh-how1-0-3_cfd-1p4-0p2-1000-0p72-0-1-_1-0__tdisc2-0p2-100-_3-1_-3.allpost.instant.conv};

\addplot [red, dashed, thick, mark=*, mark repeat=5, mark options={solid}]  table[x index=0, y index=4]  {foil2d/cfd/flap01morph01/opt/thrust0/freeze_x/testcase00/ascii/app1_mesh-how1-0-3_cfd-1p4-0p2-1000-0p72-0-1-_1-0__tdisc2-0p2-100-_3-1_-3.allpost.instant.conv};
\addplot [blue, dashed, thick, mark=square*, mark repeat=5, mark options={solid}]  table[x index=0, y index=4]  {foil2d/cfd/flap01morph01/opt/thrust1/freeze_x/testcase00/ascii/app1_mesh-how1-0-3_cfd-1p4-0p2-1000-0p72-0-1-_1-0__tdisc2-0p2-100-_3-1_-3.allpost.instant.conv};
\addplot [magenta, dashed, thick, mark=triangle*, mark repeat=5, mark options={solid}]  table[x index=0, y index=4]  {foil2d/cfd/flap01morph01/opt/thrust2p5/freeze_x/testcase00/ascii/app1_mesh-how1-0-3_cfd-1p4-0p2-1000-0p72-0-1-_1-0__tdisc2-0p2-100-_3-1_-3.allpost.instant.conv};

\nextgroupplot[xmin=0, xmax=60, xlabel={iteration}, ylabel=$J_x$]
\addplot [red, solid, thick, mark=*, mark repeat=5]  table[x index=0, y index=1]  {foil2d/cfd/flap01/opt/thrust0/freeze_x/testcase00/ascii/app1_mesh-how1-0-3_cfd-1p4-0p2-1000-0p72-0-1-_1-0__tdisc2-0p2-100-_3-1_-3.allpost.instant.conv};
\addplot [blue, solid, thick, mark=square*, mark repeat=5]  table[x index=0, y index=1]  {foil2d/cfd/flap01/opt/thrust1/freeze_x/testcase00/ascii/app1_mesh-how1-0-3_cfd-1p4-0p2-1000-0p72-0-1-_1-0__tdisc2-0p2-100-_3-1_-3.allpost.instant.conv};
\addplot [magenta, solid, thick, mark=triangle*, mark repeat=5]  table[x index=0, y index=1]  {foil2d/cfd/flap01/opt/thrust2p5/freeze_x/testcase00/ascii/app1_mesh-how1-0-3_cfd-1p4-0p2-1000-0p72-0-1-_1-0__tdisc2-0p2-100-_3-1_-3.allpost.instant.conv};

\addplot [red, dashed, thick, mark=*, mark repeat=5, mark options={solid}]  table[x index=0, y index=1]  {foil2d/cfd/flap01morph01/opt/thrust0/freeze_x/testcase00/ascii/app1_mesh-how1-0-3_cfd-1p4-0p2-1000-0p72-0-1-_1-0__tdisc2-0p2-100-_3-1_-3.allpost.instant.conv};
\addplot [blue, dashed, thick, mark=*, mark repeat=5, mark options={solid}]  table[x index=0, y index=1]  {foil2d/cfd/flap01morph01/opt/thrust1/freeze_x/testcase00/ascii/app1_mesh-how1-0-3_cfd-1p4-0p2-1000-0p72-0-1-_1-0__tdisc2-0p2-100-_3-1_-3.allpost.instant.conv};
\addplot [magenta, dashed, thick, mark=*, mark repeat=5, mark options={solid}]  table[x index=0, y index=1]  {foil2d/cfd/flap01morph01/opt/thrust2p5/freeze_x/testcase00/ascii/app1_mesh-how1-0-3_cfd-1p4-0p2-1000-0p72-0-1-_1-0__tdisc2-0p2-100-_3-1_-3.allpost.instant.conv};

\end{groupplot}
\end{tikzpicture}

%% file: dat/flap-summ.tab
\begin{tabular}{l|ccccccc} \toprule
& Initial & \multicolumn{3}{c}{FI}
          & \multicolumn{3}{c}{FII} \\\midrule
$q$ & & $0.0$ & $1.0$ & $2.5$ & $0.0$ & $1.0$ & $2.5$ \\\midrule 
\rowcolor{black!00} $W$        & -9.4096e+00 & -4.5695e-01 & -2.1419e+00 & -4.9476e+00 & -4.3252e-01 & -2.0271e+00 & -4.6110e+00 \\
\rowcolor{black!10} $W_x$      &  0.0000e+00 &  0.0000e+00 &  0.0000e+00 &  0.0000e+00 &  0.0000e+00 &  0.0000e+00 &  0.0000e+00 \\
\rowcolor{black!00} $W_y$      & -9.4096e+00 & -4.2807e-01 & -2.0642e+00 & -4.7967e+00 & -3.7363e-01 & -1.5413e+00 & -3.3712e+00 \\
\rowcolor{black!10} $W_\theta$ &  0.0000e+00 &  2.8883e-02 &  7.7694e-02 &  1.5083e-01 &  1.7101e-02 &  6.7900e-03 & 1.8744e-01 \\
\rowcolor{black!00} $F_x$      & -1.7660e-01 & -4.0490e-11 & -1.0000e+00 & -2.5000e+00 &  1.6937e-10 & -1.0000e+00 & -2.5000e+00 \\
\rowcolor{black!10} $F_y$      &  3.5413e-01 &  1.5989e-02 &  5.0480e-02 &  9.7240e-02 & -1.5292e-02 &  4.8657e-02 &  9.6440e-02 \\\bottomrule
\end{tabular}